\documentclass[11pt]{article}
\usepackage{amsmath,amsthm,amssymb,amsfonts,amscd}
\usepackage{amsxtra}
\usepackage{eucal}
\usepackage{mathrsfs}
\usepackage{color}
\usepackage{xcolor}
\usepackage{hyperref}
\usepackage{enumerate}
\usepackage{blkarray}
\usepackage{multirow}
\usepackage{tikz}
\usetikzlibrary{topaths,calc}
\usepackage{enumitem}
\usepackage{tcolorbox}
\usepackage{graphicx}
\usepackage{pdfpages}
\usepackage[numbers,sort]{natbib}

\tikzstyle{edge} = [fill opacity=.5,line cap=round, line join=round, line width=20pt]
\tikzstyle{thickedge} = [fill opacity=.5,line cap=round, line join=round, line width=30pt]

\pgfdeclarelayer{background}
\pgfsetlayers{background,main}

\usepackage[a4paper,margin=2cm]{geometry}
\definecolor{forestgreen(traditional)}{rgb}{0.0, 0.27, 0.13}
\definecolor{forestgreen(web)}{rgb}{0.13, 0.55, 0.13}
\definecolor{airforceblue}{rgb}{0.36, 0.54, 0.66}

\newcommand{\rnk}{\mathrm{rank}}
\newcommand{\GF}{\mathrm{GF}}
\newcommand{\supp}{\mathrm{supp}}
\newcommand{\innk}[2]{\langle #1, #2 \rangle_K}
\newcommand{\inn}[2]{\langle #1, #2 \rangle}
\newcommand{\mb}[1]{\mathbf{#1}}
\newcommand{\mbt}{\mathbf{t}}

\newcommand\abs[1]{\lvert #1\rvert}

\hypersetup{
  colorlinks,
  citecolor=forestgreen(web),
  linkcolor=airforceblue,
  urlcolor=magenta,
  pdftitle={Prime vertex-minors of a prime graph},
  pdfauthor={Donggyu Kim and Sang-il Oum}}

\newtheorem{theorem}{Theorem}[section]
\newtheorem{proposition}[theorem]{Proposition}
\newtheorem{lemma}[theorem]{Lemma}

\newtheorem{corollary}[theorem]{Corollary}

\theoremstyle{definition}

\theoremstyle{remark}

\usepackage{authblk}
\usepackage{lineno}
\title{Prime vertex-minors of a prime graph}
\author[2,1]{Donggyu Kim\thanks{Supported by the Institute for Basic Science (IBS-R029-C1).}}
\author[$*$1,2]{Sang-il Oum}
\affil[1]{Discrete Mathematics Group,
Institute for Basic Science (IBS),
Daejeon,~South~Korea}
\affil[2]{Department of Mathematical Sciences, KAIST, Daejeon, South~Korea}
\affil[ ]{Email: \texttt{donggyu@kaist.ac.kr}, \texttt{sangil@ibs.re.kr}}
\date{\today}	

\begin{document}
\maketitle
\begin{abstract}
  A graph is prime if it does not admit a partition $(A,B)$ of its vertex set such that $\min\{|A|,|B|\} \geq 2$ and 
  the rank of the $A\times B$ submatrix of its adjacency matrix is at most $1$.
  A vertex $v$ of a graph is non-essential if at least two of the three kinds of vertex-minor reductions at~$v$ result in prime graphs.
  
  In 1994, Allys proved that every prime graph with at least four vertices has a non-essential vertex unless it is locally equivalent to a cycle graph.
  We prove that every prime graph with at least four vertices has at least two non-essential vertices unless it is locally equivalent to a cycle graph.
  As a corollary, we show that for a prime graph $G$ with at least six vertices and a vertex $x$, there is a vertex $v \ne x$ such that $G \setminus v$ or $G * v \setminus v$ is prime, unless $x$ is adjacent to all other vertices and $G$ is isomorphic to a particular graph on odd number of vertices.

  Furthermore, we show that a prime graph with at least four vertices has at least three non-essential vertices, unless it is locally equivalent to a graph consisting of at least two internally-disjoint paths 
  between two fixed distinct vertices having no common neighbors.
  We also prove analogous results for pivot-minors.
\end{abstract}

\section{Introduction}
An edge $e$ of a graph $G$ is \emph{non-essential} if the deletion of $e$ in $G$, denoted by $G \setminus e$, or the contraction of $e$ in $G$, denoted by $G / e$, is simple and $3$-connected.
Tutte's wheel theorem~{\cite[(4.1)]{Tutte1961}} states that every simple $3$-connected graph has a non-essential edge unless it is isomorphic to a wheel graph.
As a generalization, Oxley and Wu~{\cite{Oxley2000a}} showed that every simple $3$-connected graph has at least two non-essential edges unless it is isomorphic to a wheel graph, and in~\cite{Oxley2000b} determined all simple $3$-connected graphs having exactly two non-essential edges.
Moreover, they~\cite{Oxley2004} investigated all simple $3$-connected graphs having exactly three non-essential edges.
Indeed, all of these results except the last have corresponding results for matroids; see~\cite{Tutte1966cm,Oxley2000a,Oxley2000b}.

We aim to prove analogous theorems for the vertex-minor relation.
In this paper, except for the first paragraph, all graphs are assumed to be simple, meaning that they have neither loops nor parallel edges.
For a vertex $v$ of a graph $G$, let $G*v$ be the graph obtained from $G$ by deleting all edges joining two neighbors of $v$ and adding edges joining non-adjacent pairs of two neighbors of $v$.
This operation is called the \emph{local complementation} at $v$ to $G$.
Two graphs are \emph{locally equivalent} if one can be obtained from the other by applying a sequence of local complementations.
A graph $H$ is a \emph{vertex-minor} of a graph $G$ if $H$ is an induced subgraph of a graph locally equivalent to~$G$.
For every edge $vw$, $G*v*w*v = G*w*v*w$ by Bouchet~{\cite[(8.2)]{Bouchet1988iso}}.
For an edge $vw$, let $G \wedge vw := G*v*w*v$.
This operation is called the \emph{pivoting} $vw$ to $G$.
Two graphs are \emph{pivot-equivalent} if one can be obtained from the other by applying a sequence of pivotings.
A graph $H$ is a \emph{pivot-minor} of a graph $G$ if $H$ is an induced subgraph of a graph pivot-equivalent to~$G$.
Graphs with the pivot-minor relation are closely related to binary matroids with the minor relation.
For example, fundamental graphs of a binary matroid $M$ are pairwise pivot-equivalent and fundamental graphs of minors of $M$ are pivot-minors of a fundamental graph of $M$; see~{\cite[Section~3]{Oum2005}} for details.

A \emph{split} of a graph is a partition $(A,B)$ of its vertex set such that $\min\{|A|,|B|\} \geq 2$ and for some $A'\subseteq A$ and $B'\subseteq B$, two vertices $x\in A$ and $y\in B$ are adjacent if and only if $x\in A'$ and $y\in B'$.
Equivalently, the partition $(A,B)$ is a split if and only if $\min\{|A|,|B|\} \geq 2$ and the $A\times B$ submatrix of the adjacency matrix of the given graph has rank at most $1$.
A graph is \emph{prime} if it has no split.
Bouchet~{\cite[Corollary~3]{Bouchet1989iso}} showed that locally equivalent graphs have the same set of splits.
Thus, if a graph $G$ is prime, then every graph locally equivalent to $G$ is prime.

For a vertex $v$ and distinct neighbors $w$ and $w'$ of $v$, $G\wedge vw\wedge ww' = G\wedge vw'$; see Oum~{\cite[Proposition~2.5]{Oum2005}}.
Hence $G\wedge vw \setminus v$ and $G\wedge vw' \setminus v$ are pivot-equivalent (so locally equivalent) because $G \wedge vw \setminus v \wedge ww' = G \wedge vw' \setminus v$.
Let $G/v$ denote  $G \wedge vw \setminus v$ for an arbitrary neighbor $w$ of~$v$ if $v$ has a neighbor 
and $G\setminus v$ otherwise.
Note that $G/v$ is well defined up to pivot equivalence (and up to local equivalence).
Bouchet~{\cite[(9.2)]{Bouchet1988iso}} proved that for a graph $G$ and a vertex $v$, every vertex-minor of $G$ on $V(G) - \{v\}$ is locally equivalent to $G\setminus v$, $G*v\setminus v$, or $G / v$.

A vertex $v$ of a graph $G$ is \emph{non-essential} if at least two of $G\setminus v$, $G*v\setminus v$, and $G/v$ are prime.
Allys~{\cite[Theorem~4.3]{Allys1994iso}} proved that every prime graph with more than $4$ vertices has a non-essential vertex unless it is locally equivalent to a cycle.
We prove that, indeed, such graphs have at least two non-essential vertices.

\begin{theorem}\label{thm:noness2}
  Every prime graph with at least four vertices has at least two non-essential vertices unless it is locally equivalent to a cycle.
\end{theorem}

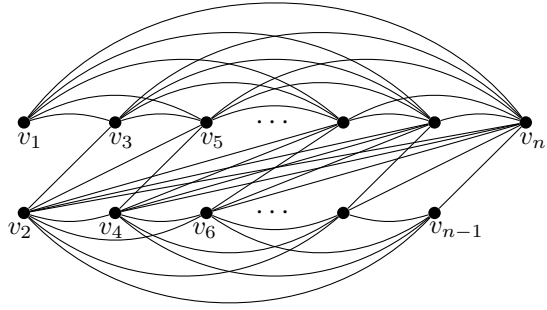
\begin{figure}
  \centering
  \begin{tikzpicture}
    \node[shape=circle,fill=black, scale=0.45] (v1) at (0,1.2) {};
    \node[shape=circle,fill=black, scale=0.45] (v3) at (1.2,1.2) {};
    \node[shape=circle,fill=black, scale=0.45] (v5) at (2.4,1.2) {};

    \node () at (3.3,1.2) {$\cdots$};

    \node[shape=circle,fill=black, scale=0.45] (vn4) at (4.2,1.2) {};
    \node[shape=circle,fill=black, scale=0.45] (vn2) at (5.4,1.2) {};
    \node[shape=circle,fill=black, scale=0.45] (vn) at (6.6,1.2) {};

    \node[shape=circle,fill=black, scale=0.45] (v2) at (0,0) {};
    \node[shape=circle,fill=black, scale=0.45] (v4) at (1.2,0) {};
    \node[shape=circle,fill=black, scale=0.45] (v6) at (2.4,0) {};

    \node () at (3.3,0) {$\cdots$};

    \node[shape=circle,fill=black, scale=0.45] (vn3) at (4.2,0) {};
    \node[shape=circle,fill=black, scale=0.45] (vn1) at (5.4,0) {};

    \draw (v2) -- (v3);
    \draw (v2) -- (v5);
    \draw (v2) -- (vn4);
    \draw (v2) -- (vn2);

    \draw (v4) -- (v5);
    \draw (v4) -- (vn4);
    \draw (v4) -- (vn2);

    \draw (v6) -- (vn4);
    \draw (v6) -- (vn2);

    \draw (vn3) -- (vn2);

    \draw (vn) -- (vn1);
    \draw (vn) -- (vn3);
    \draw (vn) -- (v6);
    \draw (vn) -- (v4);
    \draw (vn) -- (v2);

    \draw (vn) to [bend right = 17] (vn2);
    \draw (vn) to [bend right = 29] (vn4);
    \draw (vn) to [bend right = 41] (v5);
    \draw (vn) to [bend right = 47] (v3);
    \draw (vn) to [bend right = 53] (v1);

    \node () at (0.08,1.265-0.3) {\small{$v_1$}};
    \node () at (1.28,1.265-0.3) {\small{$v_3$}};
    \node () at (2.48,1.265-0.3) {\small{$v_5$}};
    \node () at (6.70,1.265-0.3) {\small{$v_{n}$}};

    \node () at (-0.05,-0.235) {\small{$v_2$}};
    \node () at (1.15,-0.235) {\small{$v_4$}};
    \node () at (2.38,-0.235) {\small{$v_6$}};
    \node () at (5.69,-0.235) {\small{$v_{n-1}$}};

    \draw (v1) to [bend left = 17] (v3);
    \draw (v3) to [bend left = 17] (v5);
    \draw (v5) to [bend left = 23] (vn4);
    \draw (vn4) to [bend left = 17] (vn2);

    \draw (v1) to [bend left = 29] (v5);
    \draw (v3) to [bend left = 35] (vn4);
    \draw (v5) to [bend left = 35] (vn2);

    \draw (v1) to [bend left = 41] (vn4);
    \draw (v3) to [bend left = 41] (vn2);

    \draw (v1) to [bend left = 47] (vn2);

    \draw (v2) to [bend left = -17] (v4);
    \draw (v4) to [bend left = -17] (v6);
    \draw (v6) to [bend left = -23] (vn3);
    \draw (vn3) to [bend left = -17] (vn1);

    \draw (v2) to [bend left = -29] (v6);
    \draw (v4) to [bend left = -35] (vn3);
    \draw (v6) to [bend left = -35] (vn1);

    \draw (v2) to [bend left = -41] (vn3);
    \draw (v4) to [bend left = -41] (vn1);

    \draw (v2) to [bend left = -47] (vn1);
  \end{tikzpicture}
  \caption{$H_n$ with odd $n$.}
  \label{fig:Gvn}
\end{figure}

As a corollary, we deduce the following strengthening of Allys~{\cite[Theorem~5.3]{Allys1994iso}}.
For a positive integer $n$, let $H_{n}$ be the graph on $\{v_1,v_2,\ldots,v_n\}$ such that for every $i$ and $j$ with $1 \le i < j \le n$, $v_i$ is adjacent to $v_j$ if and only if $i$ is even or $j$ is odd; see Figure~\ref{fig:Gvn}.
We note that if $n$ is odd, then $v_n$ is a unique vertex of degree $n-1$ in $H_n$.

\begin{corollary}\label{cor:rooted prime reduction0}
  Let $G$ be a prime graph with at least six vertices and let $x$ be a vertex of $G$.
  Then
  $G$ has a vertex $v \ne x$ such that $G \setminus v$ or $G * v \setminus v$ is prime, unless
  $|V(G)|$ is odd,
  $x$ is adjacent to all other vertices, and $G$ is isomorphic to $H_{|V(G)|}$.
\end{corollary}

\begin{figure}
  \centering
  \begin{tikzpicture}
    \node[shape=circle,fill=black, scale=0.5] (v1) at (0,0) {};
    \node[shape=circle,fill=black, scale=0.5] (v2) at (0.6,0.8) {};
    \node[shape=circle,fill=black, scale=0.5] (v3) at (1.5,1.2) {};
    \node[shape=circle,fill=black, scale=0.5] (v4) at (2.4,0.8) {};
    \node[shape=circle,fill=black, scale=0.5] (v5) at (3,0) {};
    \node[shape=circle,fill=black, scale=0.5] (v6) at (2.4,-0.8) {};
    \node[shape=circle,fill=black, scale=0.5] (v7) at (0.6,-0.8) {};
    \node[shape=circle,fill=black, scale=0.5] (v8) at (1.5,0) {};
    
    \draw (-0.28,-0.02) node {$v_1$};
    \draw (0.3,0.85) node {$v_2$};
    \draw (1.5,1.44) node {$v_3$};
    \draw (2.67,0.85) node {$v_4$};
    \draw (3.28,-0.02) node {$v_5$};
    \draw (2.67,-0.88) node {$v_6$};
    \draw (0.3,-0.88) node {$v_7$};
    \draw (1.5,-0.29) node {$v_8$};
    
    \draw (v1) -- (v2) -- (v3) -- (v4) -- (v5) -- (v6) -- (v7) -- (v1) -- (v8) -- (v5);
    \draw (v2) -- (v4);

    \draw (1.5,-1.3) node {$G$};
  \end{tikzpicture}

  \vspace{0.7cm}
  \begin{tikzpicture}
    \node[shape=circle,fill=black, scale=0.5] (v1) at (0,0) {};
    \node[shape=circle,fill=black, scale=0.5] (v2) at (0.6,0.8) {};
    \node[shape=circle,fill=black, scale=0.5, opacity=0.35] (v3) at (1.5,1.2) {};
    \node[shape=circle,fill=black, scale=0.5] (v4) at (2.4,0.8) {};
    \node[shape=circle,fill=black, scale=0.5] (v5) at (3,0) {};
    \node[shape=circle,fill=black, scale=0.5] (v6) at (2.4,-0.8) {};
    \node[shape=circle,fill=black, scale=0.5] (v7) at (0.6,-0.8) {};
    \node[shape=circle,fill=black, scale=0.5] (v8) at (1.5,0) {};

    \draw (v1) -- (v2) -- (v4) -- (v5) -- (v6) -- (v7) -- (v1) -- (v8) -- (v5);
    \draw[opacity=0.35,dashed] (v2) -- (v3) -- (v4);

    \draw (1.5,-1.3) node {$G \setminus v_3$};
  \end{tikzpicture}
  \hspace{1cm}
  \begin{tikzpicture}
    \node[shape=circle,fill=black, scale=0.5] (v1) at (0,0) {};
    \node[shape=circle,fill=black, scale=0.5] (v2) at (0.6,0.8) {};
    \node[shape=circle,fill=black, scale=0.5, opacity=0.35] (v3) at (1.5,1.2) {};
    \node[shape=circle,fill=black, scale=0.5] (v4) at (2.4,0.8) {};
    \node[shape=circle,fill=black, scale=0.5] (v5) at (3,0) {};
    \node[shape=circle,fill=black, scale=0.5] (v6) at (2.4,-0.8) {};
    \node[shape=circle,fill=black, scale=0.5] (v7) at (0.6,-0.8) {};
    \node[shape=circle,fill=black, scale=0.5] (v8) at (1.5,0) {};

    \draw (v1) -- (v2);
    \draw (v4) -- (v5) -- (v6) -- (v7) -- (v1) -- (v8) -- (v5);
    \draw[opacity=0.35,dashed] (v2) -- (v3) -- (v4);

    \draw (1.5,-1.3) node {$G*v_3\setminus v_3$};
  \end{tikzpicture}
  \hspace{1cm}
  \begin{tikzpicture}
    \node[shape=circle,fill=black, scale=0.5] (v1) at (0,0) {};
    \node[shape=circle,fill=black, scale=0.5] (v2) at (0.6,0.8) {};
    \node[shape=circle,fill=black, scale=0.5, opacity=0.35] (v3) at (1.5,1.2) {};
    \node[shape=circle,fill=black, scale=0.5] (v4) at (2.4,0.8) {};
    \node[shape=circle,fill=black, scale=0.5] (v5) at (3,0) {};
    \node[shape=circle,fill=black, scale=0.5] (v6) at (2.4,-0.8) {};
    \node[shape=circle,fill=black, scale=0.5] (v7) at (0.6,-0.8) {};
    \node[shape=circle,fill=black, scale=0.5] (v8) at (1.5,0) {};

    \draw (v1) to [bend right=5] (v4);
    \draw (v2) -- (v4) -- (v5) -- (v6) -- (v7) -- (v1) -- (v8) -- (v5);
    \draw[opacity=0.35,dashed] (v3) to [bend left=15] (v1);
    \draw[opacity=0.35,dashed] (v3) -- (v2);
    \draw[opacity=0.35,dashed] (v3) -- (v4);

    \draw (1.5,-1.3) node {$G\wedge v_2v_3 \setminus v_3$};
  \end{tikzpicture}
  \caption{The set of non-essential vertices of $G$ is $\{v_1,v_5,v_8\}$ and the set of non-pivotal vertices of $G$ is $\{v_1,v_2,v_3,v_4,v_5,v_8\}$. For instance, $G\setminus v_3$ is prime and neither $G*v_3 \setminus v_3$ nor $G\wedge v_2v_3 \setminus v_3$ is prime.}
  \label{fig:nonpiv_ess}
\end{figure}

A vertex $v$ of a graph $G$ is \emph{non-pivotal} if $G \setminus v$ or $G/v$ is prime.
Obviously, every non-essential vertex is non-pivotal.
However, a non-pivotal vertex is not necessarily non-essential; see Figure~\ref{fig:nonpiv_ess}.
As a corollary of Theorem~\ref{thm:noness2}, every prime graph with at least $5$ vertices has at least $2$ non-pivotal vertices unless it is locally equivalent to a cycle.
We extend this observation for graphs not pivot-equivalent to cycles.

\begin{theorem}
  \label{thm:nonpiv2}
  Every prime graph with at least four vertices has at least two non-pivotal vertices unless it is pivot-equivalent to a cycle.
\end{theorem}

Theorems~\ref{thm:noness2} and~\ref{thm:nonpiv2} do not imply each other immediately.
As mentioned earlier, a non-pivotal vertex need not be non-essential
and therefore Theorem~\ref{thm:nonpiv2} does not seem to imply Theorem~\ref{thm:noness2} directly.
Also, it is nontrivial to deduce Theorem~\ref{thm:nonpiv2} from Theorem~\ref{thm:noness2}, because 
locally equivalent graphs need not be pivot-equivalent.
For instance, for a cycle $C$ of length at least $6$ and its vertex~$v$, two graphs $C$ and $C*v$ are locally equivalent but not pivot-equivalent.
We remark that if two bipartite graphs are locally equivalent, then they are pivot-equivalent; shown by Fon-Der-Flaass~\cite{FonDerFlaass1989,FonDerFlaass1996}.
In Section~\ref{sec:nonpiv}, we prove Theorem~\ref{thm:nonpiv2} using Theorem~\ref{thm:noness2} together with nontrivial lemmas.

One can easily check that every cycle of length at least $5$ has no non-pivotal vertex (and no non-essential vertex).
By limiting our focus to bipartite graphs, we obtain the following corollary.
Note that every pivot-minor of a bipartite graph is bipartite; see Oum~{\cite[Corollary~2.3]{Oum2005}}.

\begin{corollary}\label{cor:bippiv}
  Every prime bipartite graph with at least four vertices has at least two non-pivotal vertices unless it is pivot-equivalent to an even cycle.
\end{corollary}

We remark that Corollary~\ref{cor:bippiv} is equivalent to Corollary~3.5 of Oxley and Wu~\cite{Oxley2000a} restricted to binary matroids.
We will explain this equivalence together with the relation between pivot-minors of bipartite graphs and minors of binary matroids in Appendix~\ref{app:binary matroid}.
It is easy to deduce from Corollary~\ref{cor:bippiv} a result of Oxley and Wu~{\cite[Corollary~3.5]{Oxley2000a}} on graphs, which states that every %
$3$-connected graph has at least $2$ non-essential edges unless it is isomorphic to a wheel graph.

We also characterize prime graphs with at least $3$ non-essential (or non-pivotal) vertices as follows.
Let $\Theta$ be the set of graphs consisting of at least two internally-disjoint paths between two fixed distinct vertices having no common neighbor.

\begin{theorem}\label{thm:noness3}
  A prime graph with at least four vertices has at least three non-essential vertices if and only if it is not locally equivalent to any graph in $\Theta$.
\end{theorem}

\begin{theorem}\label{thm:nonpiv3}
  A prime graph with at least four vertices has at least three non-pivotal vertices if and only if it is not pivot-equivalent to any graph in $\Theta$.
\end{theorem}

As we deduced Corollary~\ref{cor:bippiv} from Theorem~\ref{thm:nonpiv2}, we obtain the next corollary from Theorem~\ref{thm:nonpiv3}.

\begin{corollary}\label{cor:bippiv3}
  A bipartite prime graph with at least four vertices has at least three non-pivotal vertices if and only if it is not pivot-equivalent to any bipartite graph in $\Theta$.
\end{corollary}

We remark that similarly, Corollary~\ref{cor:bippiv3} implies Theorems~1.3 and~1.4 of Oxley and Wu~\cite{Oxley2000b} restricted to binary matroids.
This implication will be explained in Appendix~\ref{app:binary matroid}.

Bouchet~\cite{Bouchet1987cir} showed that 
every prime graph with at least six vertices has a prime vertex-minor with one fewer vertex, which was used in a recognition algorithm and the proof of obstructions for circle graphs of Bouchet~\cite{Bouchet1987cir, Bouchet1994circle}.
From Theorem~\ref{thm:noness3}, we obtain the following strengthening of Bouchet's result.

\begin{corollary}\label{cor:rooted prime reduction1}
  Let $G$ be a prime graph with at least six vertices and let $x$ and $y$ be vertices of $G$.
  Then
  there is a prime vertex-minor $H$ of $G$ such that $|V(H)| = |V(G)|-1$ and $x,y \in V(H)$.
\end{corollary}

Similarly, we deduce the following for pivot-minors from Theorems~\ref{thm:nonpiv2} and~\ref{thm:nonpiv3}. %

\begin{corollary}\label{cor:rooted prime reduction2}
  Let $G$ be a prime graph with at least four vertices and let $x$ and $y$ be vertices of $G$.
  Then
  there is a prime pivot-minor $H$ of $G$ such that $|V(H)| = |V(G)|-1$ and $x,y \in V(H)$, unless
  \begin{enumerate}[label=\rm(\roman*)]
    \item\label{item:rpr2i} $G$ is pivot-equivalent to a cycle, or
    \item\label{item:rpr2ii} $x \neq y$ and $G$ is pivot-equivalent to a graph consisting of at least three internally-disjoint paths between $x$ and $y$ that have no common neighbor.
  \end{enumerate}
\end{corollary}

Our results heavily depend on the following theorem providing non-essential vertices in a prime graph.
The \emph{adjacency matrix} $A_G = (a_{vw})$ of a graph $G = (V,E)$ is a $V\times V$ matrix over the binary field $\GF(2)$ such that $a_{vw}=1$ if and only if $vw \in E$.
The \emph{cut-rank function}~$\rho_G$ of~$G$ is a function from $2^{V}$ to~$\mathbb{Z}$ such that $\rho_G(X) := \rnk(A_G[X,V-X])$, where $A_G[X,V-X]$ is the $X \times (V-X)$ submatrix of $A_G$.
A \emph{tight path} $P$ in a $3$-uniform hypergraph $H$ is a hypergraph such that $V(P) \subseteq V(H)$, $E(P) \subseteq E(H)$, and
$V(P)$ admits an ordering $v_0,v_1,\ldots,v_{k+1}$ for which $k \geq 1$ and $E(P) = \{\{v_{i-1},v_{i},v_{i+1}\} : 1\le i \le k\}$. 
An \emph{end} of a tight path $P$ is a vertex incident with exactly one edge of $P$.
Note that if $k\ge 2$, $P$ has exactly two ends $v_0$ and $v_{k+1}$, and if $k=1$, then $P$ has exactly three ends $v_0$, $v_1$, and $v_2$.

\begin{theorem}\label{thm:maxpathG}
  Let $G$ be a prime graph with at least five vertices.
  Let $H$ be a $3$-uniform hypergraph on $V(G)$ such that 
  $E(H) = \{X \subseteq V(G) : |X| = 3 \text{ and } \rho_G(X) = 2\}$.
  If $G$ is not locally equivalent to a cycle, then at least two ends of a maximal tight path in $H$ are non-essential in $G$.
\end{theorem}

This paper is organized as follows.
In Section~\ref{sec:pre} we review graph-theoretic notions.
In Section~\ref{sec:isosys} we review isotropic systems defined by Bouchet~\cite{Bouchet1987iso,Bouchet1988iso,Bouchet1989iso}.
Isotropic systems capture the local equivalence of graphs in terms of linear algebra.
In Section~\ref{sec:tri} we define triangles in an isotropic system and describe properties of paths formed by triangles, and we prove Theorems~\ref{thm:maxpathG} and~\ref{thm:noness2}.
In Section~\ref{sec:main} we prove half of Theorem~\ref{thm:noness3}.
In order to complete the proof of Theorem~\ref{thm:noness3}, in Section~\ref{sec:theta} we classify graphs consisting of internally-disjoint paths sharing their ends according to the number of non-essential vertices.
In Section~\ref{sec:nonpiv} we prove Theorem~\ref{thm:nonpiv2}, Theorem~~\ref{thm:nonpiv3}, Corollary~\ref{cor:bippiv}, and Corollary~\ref{cor:bippiv3}.
In Section~\ref{sec:applications} we present proofs of Corollaries~\ref{cor:rooted prime reduction0},~\ref{cor:rooted prime reduction1}, and~\ref{cor:rooted prime reduction2}.

\section{Preliminaries}\label{sec:pre}

For a set $X$, we write $2^X$ to denote the set of subsets of $X$.
For sets $X$ and~$Y$, we write $X\triangle Y := (X - Y) \cup (Y-X)$.
For an $R\times C$ matrix $M$ over a field $\mathbb{F}$ and subsets $X \subseteq R$ and $Y \subseteq C$, let $M[X,Y]$ be the submatrix of $M$ whose rows are indexed by $X$ and columns are indexed $Y$.

For a graph $G = (V,E)$ and $v\in V$, let $G \setminus v$ denote a graph obtained from $G$ by deleting $v$.
For $X \subseteq V$,
let $G[X]$ be the induced subgraph of $G$ whose vertex set is $X$ and let $G-X := G[V - X]$.
Let $N_G(v)$ be the set of neighbors of a vertex $v$.
Subdividing an edge $e=vw$ is an operation replacing $e$ with a path of length two whose ends are $v$ and $w$.
An \emph{isolated vertex} of a graph is a vertex without a neighbor, and a \emph{pendant vertex} is a vertex with a unique neighbor.
Two vertices are \emph{twins} if none of the other vertices are adjacent to exactly one of them.
For an integer $k \geq 3$, we write $C_k$ to denote a cycle graph of length~$k$.
For positive integers $\ell_1, \ell_2, \dots \ell_m$ with $|\{i:\ell_i=1\}| \leq 1$, let $\theta(\ell_1, \dots, \ell_m)$ be the graph consisting of $m$ internally-disjoint paths between two fixed vertices of lengths $\ell_1, \dots, \ell_m$, respectively.
Recall that $\Theta$ is the set of graphs $\theta(\ell_1,\dots,\ell_m)$ such that $m\geq 2$ and $\ell_i \neq 2$ for all $i$.

\paragraph*{Prime graphs}

Note that a partition $(X,V-X)$ of $V$ is a split of a graph $G$ if and only if $\min\{|X|,|V-X|\}\geq 2$ and $\rho_G(X) \leq 1$.
It is easy to observe the following.
\begin{lemma}\label{lem:prime}
  Every prime graph with at least $4$ vertices is $2$-connected and has no twins.\qed
\end{lemma}
Thus one may observe that every prime graph with at least $4$ vertices has neither isolated vertices nor pendant vertices.
It also follows that there is no prime graph with exactly $4$ vertices.

\begin{figure}
  \centering
  \begin{tikzpicture}
    \node[shape=circle,fill=black, scale=0.5] (v1) at (30:1.2) {};
    \node[shape=circle,fill=black, scale=0.5] (v2) at (150:1.2) {};
    \node[shape=circle,fill=black, scale=0.5] (v3) at (270:1.2) {};
    \node[shape=circle,fill=black, scale=0.5] (w1) at (90:2.4) {};
    \node[shape=circle,fill=black, scale=0.5] (w2) at (210:2.4) {};
    \node[shape=circle,fill=black, scale=0.5] (w3) at (330:2.4) {};
    
    \draw (150:1.5) node {$v$};
    \draw (30:1.5) node {$w$};
    \draw (v1) -- (v2);
    \draw (v2) -- (v3);
    \draw (v3) -- (v1);
    \draw (w1) -- (v1);
    \draw (w1) -- (v2);
    \draw (w2) -- (v2);
    \draw (w2) -- (v3);
    \draw (w3) -- (v3);
    \draw (w3) -- (v1);

    \draw (270:2) node {$G$};
  \end{tikzpicture}
  \hspace{1cm}
  \begin{tikzpicture}
    \node[shape=circle,fill=black, scale=0.5] (v1) at (30:1.2) {};
    \node[shape=circle,fill=black, scale=0.5] (v2) at (150:1.2) {};
    \node[shape=circle,fill=black, scale=0.5] (v3) at (270:1.2) {};
    \node[shape=circle,fill=black, scale=0.5] (w1) at (90:2.4) {};
    \node[shape=circle,fill=black, scale=0.5] (w2) at (210:2.4) {};
    \node[shape=circle,fill=black, scale=0.5] (w3) at (330:2.4) {};
    
    \draw (150:1.5) node {$v$};
    \draw (v1) -- (v2);
    \draw (v2) -- (v3);
    \draw (w1) -- (v2);
    \draw (w2) -- (v2);
    \draw (w3) -- (v3);
    \draw (w3) -- (v1);

    \draw (v1) -- (w2);
    \draw (v3) -- (w1);
    \draw (w1) to [bend right=30] (w2);

    \draw (270:2) node {$G*v$};
  \end{tikzpicture}
  \hspace{1cm}
  \begin{tikzpicture}
    \node[shape=circle,fill=black, scale=0.5] (v1) at (30:1.2) {};
    \node[shape=circle,fill=black, scale=0.5] (v2) at (150:1.2) {};
    \node[shape=circle,fill=black, scale=0.5] (v3) at (270:1.2) {};
    \node[shape=circle,fill=black, scale=0.5] (w1) at (90:2.4) {};
    \node[shape=circle,fill=black, scale=0.5] (w2) at (210:2.4) {};
    \node[shape=circle,fill=black, scale=0.5] (w3) at (330:2.4) {};
    
    \draw (150:1.5) node {$w$};
    \draw (30:1.5) node {$v$};
    \draw (v1) -- (v2);
    \draw (v2) -- (v3);
    \draw (v3) -- (v1);
    \draw (w1) -- (v1);
    \draw (w1) -- (v2);
    \draw (w2) -- (v2);
    \draw (w3) -- (v1);

    \draw (w1) to [bend right=30] (w2);
    \draw (w2) to [bend right=15] (w3);
    \draw (w3) to [bend right=30] (w1);

    \draw (270:2) node {$G\wedge vw$};
  \end{tikzpicture}
  \caption{Local complementation and pivoting.}
  \label{fig:locpiv}
\end{figure}
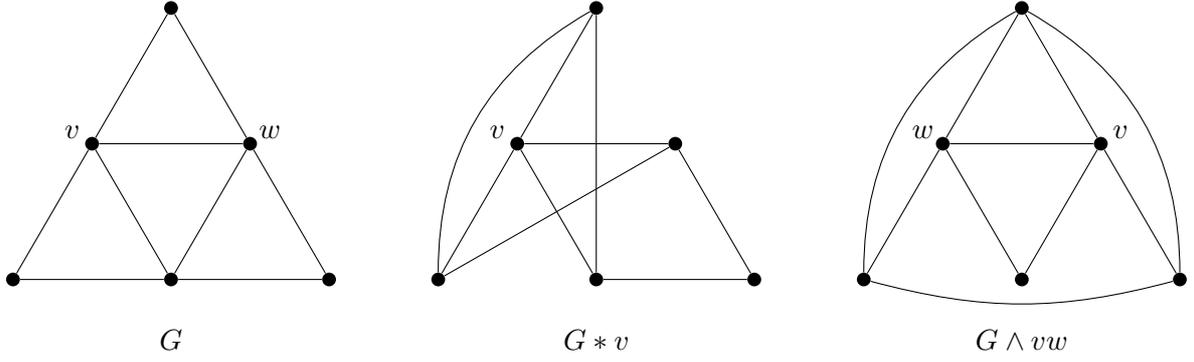

\paragraph*{Vertex-minors}
Recall that for a graph $G$ and a vertex $v$ of $G$,
\[
  G*v = (V(G), E(G) \triangle \{xy:\text{$x$ and $y$ are two distinct neighbors of $v$}\}).
\]
\begin{proposition}[Oum~{\cite[Proposition~2.1]{Oum2005}}]
  \label{prop:pivot1}
  For a graph $G$ and an edge $vw$ of $G$, let $G'$ be a graph on $V(G)$ such that
  \begin{align*}
    E(G')
    =
    E(G)
    &\triangle \{xy: x \in N_G(v)-(N_G(w) \cup \{w\}),\; y \in N_G(w)-(N_G(v) \cup \{v\}\}) \\
    &\triangle \{xy : x \in N_G(v)-(N_G(w) \cup \{w\}),\; y \in N_G(v) \cap N_G(w)\} \\
    &\triangle \{xy : x \in N_G(w)-(N_G(v) \cup \{v\}),\; y \in N_G(v) \cap N_G(w)\}.
  \end{align*}
  Then $G \wedge vw$ is equal to the graph obtained from $G'$ by exchanging the labels of $v$ and $w$; see Figure~\ref{fig:locpiv}.
\end{proposition}

A vertex~$v$ of a graph~$G$ is \emph{essential} if at most one of $G\setminus v$, $G * v \setminus v$, and $G / v$ is prime, and $v$ is \emph{pivotal} if neither $G \setminus v$ nor $G / v$ is prime.
Note that every pivotal vertex is essential.

The following proposition can be seen easily from the theory of isotropic systems~\cite{Bouchet1988iso}, and Geelen and Oum~\cite{Geelen2009circle} presented a short graph-theoretic proof.

\begin{proposition}[Geelen and Oum~{\cite[Lemma~3.1]{Geelen2009circle}}]
  \label{prop:vertexminor_localequiv}
  Let $G$ be a graph and $v,w$ be its vertices.
  \begin{enumerate}[label=\rm(\roman*)]
    \item If $v\neq w$ and $vw$ is not an edge of $G$, then $G*w \setminus v$, $G*w*v \setminus v$, $G*w/v$ are locally equivalent to $G \setminus v$, $G*v \setminus v$, $G/v$, respectively.
    \item If $v\neq w$ and $vw$ is an edge of $G$, then $G*w \setminus v$, $G*w*v \setminus v$, $G*w/v$ are locally equivalent to $G \setminus v$, $G/v$, $G*v \setminus v$, respectively.
    \item If $v=w$, then $G*w \setminus v$, $G*w*v \setminus v$, $G*w/v$ are locally equivalent to $G*v \setminus v$, $G \setminus v$, $G/v$, respectively.
  \end{enumerate}
\end{proposition}

As a corollary, we deduce the following result.

\begin{corollary}\label{cor:ess2}
  Locally equivalent graphs have the same set of non-essential vertices.
  \qed
\end{corollary}

It is known from {\cite[Proposition~2.5]{Oum2005}}
that 
$G\wedge xz=G\wedge xy\wedge yz$ for two edges $xy$ and $xz$ in a graph $G$,
and therefore we deduce the following corollary.
\begin{corollary}\label{cor:ess3}
  Pivot-equivalent graphs have the same set of non-pivotal vertices.
  \qed
\end{corollary}
\paragraph*{3-uniform hypergraphs and tight paths}
A pair $H =(V,E)$ is a \emph{hypergraph} if $V$ is a finite set and $E$ is a set of nonempty subsets of $V$, and we denote the vertex set of $H$ by $V(H) = V$ and denote the edge set of $H$ by $E(H) = E$.
A hypergraph is \emph{3-uniform} if every edge has cardinality $3$.
A hypergraph $H'$ is a \emph{partial hypergraph} of a hypergraph $H$ if $V(H') \subseteq V(H)$ and $E(H') \subseteq E(H)$.

Recall that
for %
a $3$-uniform hypegraph~$H$, a tight path $P$ %
in $H$ is a partial hypergraph that admits an ordering $v_0,v_1,\dots,v_{k+1}$ of $V(P)$ where $k\ge 1$ and $E(P) = \{\{v_{i-1},v_i,v_{i+1} \} : 1\leq i\leq k\}$.
We usually denote~$P$ by a sequence $v_0v_1v_2\dots v_{k+1}$ of distinct vertices.
The \emph{length} of a tight path is its number of edges.
A tight path $P$ in $H$ is \emph{maximal} if there is no tight path $Q$ in $H$ such that $E(P) \subsetneq E(Q)$.

\begin{lemma}
  A tight path $P = v_0v_1\dots v_{k+1}$ of length $k\geq 1$ in a $3$-uniform hypergraph $H$ is not maximal if and only if at least one of the following holds:
  \begin{enumerate}[label=\rm(\roman*)]
    \item There is a vertex $w$ such that $v_0v_1\dots v_{k+1}w$ or $wv_0v_1\dots v_{k+1}$ is a tight path in $H$.
    \item $k=1$ and there is a vertex $w$ such that $v_1v_0v_2w$ is a tight path in $H$.
    \item $k=2$ and there is a vertex $w$ such that $v_0v_2v_1v_3w$ or $wv_0v_2v_1v_3$ is a tight path in $H$.
  \end{enumerate}
\end{lemma}
\begin{proof}
  If $P$ is not maximal, then there is a tight path $Q$ of length $k+1$ containing $P$.
  In a vertex ordering of $Q$ certifying that $Q$ is a tight path, a vertex ordering of $P$ certifying that $P$ is a tight path can be obtained by deleting the first or the last vertex which is not a vertex of $P$.
  It remains to enumerate all vertex ordering of $P$ guaranteeing that $P$ is a tight path.

  If $k\geq 3$, then $P$ admits a unique vertex ordering $v_0v_1\dots v_{k+1}$ up to reversing.
  If $k=1$, then $P$ admits three vertex ordering $v_0v_1v_2$, $v_0v_2v_1$, and $v_1v_0v_2$ up to reversing.
  If $k=2$, then $P$ admits two vertex orderings $v_0v_1v_2v_3$ and $v_0v_2v_1v_3$ up to reversing.
\end{proof}

Recall that 
an end of a tight path $P$ is a vertex incident with exactly one edge of $P$.
An \emph{internal vertex} of~$P$ is a vertex incident with at least two edges of $P$.
Equivalently, an internal vertex is a vertex of $P$ that is not an end.
It is readily shown that the set of ends of $P = v_0v_1\dots v_{k+1}$ is $\{v_0,v_{k+1}\}$ if $k\geq 2$, and $V(P)$ if $k=1$.

\section{Isotropic systems}\label{sec:isosys}

We review isotropic systems defined by Bouchet~\cite{Bouchet1987iso,Bouchet1988iso,Bouchet1989iso}.
We follow notations in~\cite{Oum2005thesis,Oum2008wqo}.

Let $K = \{0, \alpha, \beta, \gamma\}$ be a $2$-dimensional vector space over the binary field $\GF(2)$, and let $\innk{\cdot}{\cdot} : K \times K \rightarrow \GF(2)$ be a bilinear form such that $\innk{x}{y} = 1$ if and only if $0\neq x \neq y \neq 0$.
For a finite set~$V$, let $K^V$ be the set of functions from $V$ to $K$, regarded as a $2|V|$-dimensional vector space over $\GF(2)$.
Let $\inn{\cdot}{\cdot} : K^V \times K^V \rightarrow \GF(2)$ be a bilinear form such that $\inn{\mb{a}}{\mb{b}} = \sum_{v\in V} \innk{\mb{a}(v)}{\mb{b}(v)}$.
For a subspace $L$ of $K^V$, let $L^\perp := \{\mb{a}\in K^V : \text{$\inn{\mb{a}}{\mb{b}} = 0$ for all $\mb{b} \in L$}\}$.
The \emph{support} of a vector $\mb{a} \in K^V$, denoted by $\supp(\mb{a})$, is the set of elements $v$ in $V$ such that $\mb{a}(v) \neq 0$.
A vector $\mb{a} \in K^V$ is \emph{complete} if $\supp(\mb{a}) = V$.
Two vectors $\mb{a}$ and $\mb{b}$ in $K^V$ are \emph{supplementary} if they are complete and $\mb{a}(v) \neq \mb{b}(v)$ for every $v\in V$.

A subspace $L$ of $K^V$ is \emph{totally isotropic} if $\innk{\mb{a}}{\mb{b}} = 0$ for all vectors $\mb{a}$ and $\mb{b}$ in $L$, equivalently, $L \subseteq L^\perp$.
Note that for every subspace $L$ of $K^V$, we have $\dim(L) + \dim(L^\perp) = \dim(K^V) = 2|V|$; see Lang~{\cite[Theorem~6.4]{Lang1987}}.
Hence for a totally isotropic subspace $L$ of $K^V$, we have $\dim(L) \leq |V|$ where the equality holds if and only if $L^\perp = L$.
An \emph{isotropic system} is a pair $(V,L)$ consisting of a finite set $V$ and a subspace $L$ of $K^V$ such that $L$ is totally isotropic and $\dim(L) = |V|$.
For an isotropic system $S=(V,L)$, we call each element $v\in V$ a \emph{vertex} of $S$.

\subsection{Minors}

For a subset $X$ of $V$, let $p_X$ be a map from $K^V$ to $K^X$ such that $(p_X(\mb{a}))(v) = \mb{a}(v)$.
For $\mb{a} \in K^V$ and $X\subseteq V$, let $\mb{a}[X]$ be a vector in $K^V$ such that
\[
  \mb{a}[X](v) =
  \begin{cases}
    \mb{a}(v) & \text{if $v \in X$}, \\
    0             & \text{otherwise}.
  \end{cases}
\]
For a subspace $L$ of $K^V$, $v\in V$, and $x \in K - \{0\}$, let
\[
  L|_{x}^v := \{p_{V-\{v\}}(\mb{a}) \in K^{V-\{v\}} : \text{$\mb{a}\in L$ and $\mb{a}(v) \in \{0,x\}$}\}.
\]
For an isotropic system $S = (V,L)$, let $S|_{x}^v := (V - \{v\}, L|_{x}^v)$ be the \emph{elementary minor} of $S$ at $v\in V$ with respect to $x \in K-\{0\}$.
The isotropic system $S$ has three elementary minors $S|_\alpha^v$, $S|_\beta^v$, $S|_\gamma^v$ at~$v$.
Bouchet~{\cite[(8.1)]{Bouchet1987iso}} proved that every elementary minor of an isotropic system is an isotropic system.
An isotropic system $S$ is a \emph{minor} of an isotropic system $S'$ if $S = S'|^{v_1}_{x_1}\dots|^{v_t}_{x_t}$ for some vertices $v_1,\dots,v_t$ of $S'$ and $x_1,\dots,x_t \in K-\{0\}$.

\subsection{Connectivity}

For a subspace $L$ of $K^V$ and a subset $X$ of $V$, let
\begin{align*}
  L|_{\subseteq X} 
  &:= \{p_X(\mb{a}) : \text{$\mb{a} \in L$ and $\supp(\mb{a}) \subseteq X$}\}
  \ \text{and} \\
  L|_{X} 
  &:= \{p_X(\mb{a}) : \text{$\mb{a} \in L$}\}.
\end{align*}
The \emph{connectivity function} of an isotropic system $S = (V,L)$ is a function $c_S : 2^{V} \rightarrow \mathbb{Z}$ such that $c_S(X) = |X| - \dim(L|_{\subseteq X})$.
We omit the subscript $S$ in $c_S$ if it is clear from the context.

\begin{lemma}[Oum~{\cite[Lemma~5.1]{Oum2008wqo}}]
  \label{lem:ortho}
  Let $L$ be a totally isotropic subspace of $K^V$ and $X$ be a subset of $V$.
  Then $(L|_{\subseteq X})^\perp = L^\perp|_{X}$.
\end{lemma}

For an isotropic system $S = (V,L)$ and $X\subseteq V$, by Lemma~\ref{lem:ortho}, we have $(L|_{\subseteq X})^\perp = L^\perp|_{X} = L|_X$.
Hence $\dim(L|_{\subseteq X}) + \dim(L|_{X}) = \dim(K^X) = 2|X|$ and $c_S(X) = \dim(L|_X) - |X|$.

\begin{proposition}[Bouchet~\cite{Bouchet1989iso}; see Allys~{\cite[Proposition~2.3.1]{Allys1994iso}}]\label{prop:conn0}
  Let $S = (V,L)$ be an isotropic system with the connectivity function $c$.
  Then for all subsets $X$, $Y \subseteq V$, the following hold.
  \begin{enumerate}[label=\rm(\roman*)]
    \item $0\leq c(X) \leq |X|$.
    \item\label{item:c02} $c(X) = c(V-X)$.
    \item\label{item:c03} $c(X) + c(Y) \geq c(X\cup Y) + c(X\cap Y)$. %
  \end{enumerate}
\end{proposition}
The following two lemmas display handy properties of the connectivity function.

\begin{lemma}[Allys~{\cite[Lemma~3.1]{Allys1994iso}}]\label{lem:conn1}
  Let $S = (V,L)$ be an isotropic system with the connectivity function $c$.
  For a subset $X \subseteq V$ and a vertex $v\in V-X$, the following hold.
  \begin{enumerate}[label=\rm(\roman*)]
    \item\label{item:c11} $c(X)-1 \leq c(X\cup \{v\}) \leq c(X)+1$.
    \item\label{item:c12} $c(X \cup \{v\}) \leq c(X)$ if and only if $L|_{\subseteq X \cup \{v\}}$ has a vector $\mb{a}$ such that $\mb{a}(v) \neq 0$.
    \item\label{item:c13} $c(X \cup \{v\}) = c(X)-1$ if and only if $L|_{\subseteq X \cup \{v\}}$ has vectors $\mb{a}$, $\mb{b}$ such that $0\neq \mb{a}(v) \neq \mb{b}(v) \neq 0$.
  \end{enumerate}
\end{lemma}

\begin{lemma}[Allys~{\cite[Proposition~3.2]{Allys1994iso}}]\label{lem:conn2}
  Let $S = (V,L)$ be an isotropic system and $S|_x^v$ be its elementary minor such that $L$ has no vector whose support is $\{v\}$.
  Let $c$ and $c'$ be the connectivity functions of $S$ and $S|_x^v$, respectively.
  Then for a subset $X$ of $V - \{v\}$, the following hold.
  \begin{enumerate}[label=\rm(\roman*)]
    \item\label{item:c21} $c(X) - 1 \leq c'(X) \leq c(X)$ and $c(X \cup \{v\}) - 1 \leq c'(X) \leq c(X\cup \{v\})$.
    \item\label{item:c22} $c'(X) = c(X) - 1$ if and only if $L|_{\subseteq X \cup \{v\}}$ has a vector $\mb{a}$ such that $\mb{a}(v) = x$.
  \end{enumerate}
\end{lemma}

For an isotropic system $S$ and a positive integer $k$, a partition $(X,Y)$ of the vertex set of $S$ is a \emph{$k$-separation} of $S$ if $\min\{|X|, |Y|\} \geq k$ and $c_S(X) < k$.
An isotropic system is \emph{$k$-connected} if it has no $k'$-separation with $1\leq k' < k$.

The following lemma is straightforward from the definition.

\begin{lemma}
  \label{lem:pt}
  If $S = (V,L)$ is a $3$-connected isotropic system with $|V| \geq 4$, then $|\supp(\mb{a})| \geq 3$ for every nonzero vector $\mb{a} \in L$.
\end{lemma}
\begin{proof}
  Suppose that there is a nonzero vector $\mb{a} \in L$ with $|\supp(\mb{a})| \leq 2$.
  Let $X = \supp(\mb{a})$.
  Then $\dim(L|_{\subseteq X}) \geq 1$ because of $\mb{a}$, so $c_S(X) = |X| - \dim(L|_{\subseteq X}) \leq |X|-1$.
  Therefore, $S$ has an $|X|$-separation $(X,V-X)$, which contradicts the assumption that $S$ is $3$-connected.
\end{proof}

For an isotropic system $S$, a vertex $v$ is \emph{non-essential} if at least two of $S|^v_\alpha$, $S|^v_\beta$, $S|^v_\gamma$ are $3$-connected, and \emph{essential} otherwise.

\subsection{Fundamental graphs}\label{ssec:gph}

For a graph $G=(V,E)$ and two supplementary vectors $\mb{a}$ and $\mb{b}$ in $K^V$, let $L_G$ be the subspace of $K^V$ spanned by $\{\mb{a}[N_G(v)] + \mb{b}[\{v\}] : v \in V\}$.
Bouchet~{\cite[(3.1)]{Bouchet1988iso}} proved that $S = (V,L_G)$ is an isotropic system.
We call a triple $(G,\mb{a},\mb{b})$ a \emph{graphic presentation} of $S$.

A vector $\mb{a} \in K^V$ is an \emph{Eulerian} vector of an isotropic system $S = (V,L)$ if $\mb{a}$ is complete and $\mb{a}[X] \not\in L$ for every nonempty subset $X$ of $V$.
\begin{lemma}[Bouchet~{\cite[(4.1)]{Bouchet1988iso}}]
  \label{lem:eur}
  Let $S$ be an isotropic system.
  For every complete vector $\mb{c}$, there is an Eulerian vector $\mb{a}$ of $S$ supplementary to $\mb{c}$.
\end{lemma}

\begin{proposition}[Bouchet~{\cite[(4.3) and~(4.4)]{Bouchet1988iso}}]
  \label{prop:basis}
  Let $\mb{a}$ be an Eulerian vector of an isotropic system $S = (V,L)$.
  Then for each $v\in V$, there is a unique vector $\mb{b}_v \in L$ such that
  \begin{enumerate}[label=\rm(\roman*)]
    \item\label{item:basis1} $\innk{\mb{b}_v(v)}{\mb{a}(v)} = 1$, and
    \item\label{item:basis2} $\innk{\mb{b}_v(w)}{\mb{a}(w)} = 0$ for all $w \in V - \{v\}$.
  \end{enumerate}
  Moreover, $\mb{b}_v(w) \neq 0$ if and only if $\mb{b}_w(v) \neq 0$ for all distinct $v,w\in V$, and $\{ \mb{b}_v : v \in V \}$ is a basis of~$L$.
\end{proposition}

The set of such vectors $\mb{b}_v$ for all $v\in V$ is called the \emph{fundamental basis} of $L$ with respect to $\mb{a}$.
The \emph{fundamental graph} of an isotropic system $S=(V,L)$ with respect to an Eulerian vector $\mb{a}$ is a graph on~$V$ such that two vertices $v$ and $w$ are adjacent if and only if $\mb{b}_v(w) \neq 0$, where $\{\mb{b}_v:v\in V\}$ is the fundamental basis of $S$ with respect to $\mb{a}$.
Let $\mb{b}$ be the complete vector such that $\mb{b}(v) = \mb{b}_v(v)$.
Then $\mb{b}$ is supplementary to $\mb{a}$.
The following proposition shows that $(G,\mb{a},\mb{b})$ is a graphic presentation of~$S$.

\begin{proposition}[Bouchet~{\cite[(4.5)]{Bouchet1988iso}}]
  \label{prop:ig}
  Let $S$ be an isotropic system.
  \begin{enumerate}[label=\rm(\roman*)]
    \item If $(G,\mb{a},\mb{b})$ is a graphic presentation of $S$, then $\mb{a}$ is an Eulerian vector of $S$.
    \item For an Eulerian vector $\mb{a}$ of $S$, let $G$ be the fundamental graph of~$S$ with respect to $\mb{a}$, let $\{\mb{b}_v: v\in V\}$ be the fundamental basis of $S$ with respect to $\mb{a}$, and let $\mb{b}$ be the complete vector such that $\mb{b}(v) = \mb{b}_v(v)$ for all $v\in V$.
    Then $(G,\mb{a},\mb{b})$ is a graphic presentation of $S$.
    Furthermore, if $(G,\mb{a},\mb{b}')$ is a graphic presentation of $S$, then $\mb{b}' = \mb{b}$.
  \end{enumerate}

\end{proposition}

Bouchet~\cite{Bouchet1989iso} explains a relation between the connectivity function of an isotropic system and the cut-rank of its fundamental graph.

\begin{proposition}[Bouchet~{\cite[Theorem~6]{Bouchet1989iso}}]
  \label{prop:igconn}
  Let $G$ be a fundamental graph of an isotropic system $S$.
  Then $c_S(X) = \rho_G(X)$ for every subset $X$ of the vertex set of $G$.
\end{proposition}

\begin{corollary}[Bouchet~{\cite[Theorem~11]{Bouchet1989iso}}]
  \label{cor:igconn2}
  Let $G$ be a fundamental graph of an isotropic system $S$
  with at least four vertices.
  Then $S$ is $3$-connected if and only if $G$ is prime.
\end{corollary}

\begin{lemma}[Bouchet~{\cite[Theorem~23]{Bouchet1989iso}}]\label{lem:c4}
  No isotropic system on $4$ vertices is $3$-connected.
\end{lemma}

An isotropic system is \emph{cyclic} if it has a cycle graph of length at least $5$ as a fundamental graph.

\begin{lemma}[Bouchet~{\cite[Theorem~23]{Bouchet1989iso}}]\label{lem:c5}
  An isotropic system on $5$ vertices is $3$-connected if and only if it is cyclic.
\end{lemma}

For a vertex $v$ of $C_n$ with $n \geq 5$, neither $C_n \setminus v$ nor $C_n / v$ is prime, and therefore we deduce the following.

\begin{lemma}[Allys~{\cite[Lemma~4.2]{Allys1994iso}}]\label{lem:cyclic}
  If $S$ is a cyclic isotropic system with at least $5$ vertices, then $S$ is prime and every vertex is essential.
\end{lemma}

We dedicate the remainder of this subsection to explaining the relation between minors of an isotropic system and vertex-minors of its fundamental graph.

\begin{lemma}[Bouchet~{\cite[(9.4)]{Bouchet1987iso}}]
  \label{lem:eurloc}
  Let $\mb{a}$ be an Eulerian vector of an isotropic system $S = (V,L)$, and let $v$ be a vertex of $S$.
  Let $\mb{a}'$ and $\mb{a}''$ be two complete vectors such that $\mb{a}[V-\{v\}] = \mb{a}'[V-\{v\}] = \mb{a}''[V-\{v\}]$ and $\{\mb{a}(v), \mb{a}'(v), \mb{a}''(v)\} = K - \{0\}$.
  Then exactly one of $\mb{a}'$ and $\mb{a}''$ is an Eulerian vector of $S$.
\end{lemma}

We write $\mb{a} * v$ to denote such an Eulerian vector $\mb{a}'$ or $\mb{a}''$ in Lemma~\ref{lem:eurloc}.

\begin{lemma}[Bouchet~{\cite[(7.1)]{Bouchet1988iso}}]
  \label{lem:eurloc2}
  Let $\mb{a}$ and $\mb{b}$ be Eulerian vectors of an isotropic system $S$.
  Then there is a sequence of vertices $v_1, v_2, \dots, v_k$ such that $\mb{b} = \mb{a}*v_1*v_2*\dots*v_k$.
\end{lemma}

\begin{proposition}[Bouchet~{\cite[(7.6) and~(8.3)]{Bouchet1988iso}}]
  \label{prop:ig2}
  Let $(G,\mb{a},\mb{b})$ be a graphic presentation of $S$.
  For a vertex $u$ and an edge $vw$ of $G$,
  \[
    (G,\mb{a},\mb{b})*u := (G*u, \, \mb{a}+\mb{b}[\{u\}], \, \mb{a}[N_G(u)]+\mb{b})
  \]
  and
  \[
    (G,\mb{a},\mb{b}) \wedge vw := (G \wedge vw, \, \mb{a}[V-\{v,w\}]+\mb{b}[\{v,w\}], \, \mb{a}[\{v,w\}]+\mb{b}[V-\{v,w\}])
  \]
  are graphic presentations of $S$.
\end{proposition}

Therefore, fundamental graphs of an isotropic system are locally equivalent.
We say that two graphic presentations $(G,\mb{a},\mb{b})$ and $(H,\mb{c},\mb{d})$ of an isotropic system are \emph{locally equivalent} if $(H,\mb{c},\mb{d}) = (G,\mb{a},\mb{b}) * v_1 \dots * v_m$ for some vertices $v_1,\dots,v_m$.
They are \emph{pivot-equivalent} if $(H,\mb{c},\mb{d}) = (G,\mb{a},\mb{b}) \wedge e_1 \dots \wedge e_m$ for some edges $e_1,\dots,e_m$.

\begin{proposition}[Bouchet~{\cite[(9.1)]{Bouchet1988iso}}; see Oum~{\cite[Proposition 3.7]{Oum2008wqo}}]
  \label{prop:igmin}
  Let $(G,\mb{a},\mb{b})$ be a graphic presentation of $S = (V,L)$.
  Then one of the following is a graphic presentation of $S|_x^v$.
  \begin{enumerate}[label=\rm(\roman*)]
    \item $(G\setminus v, \, p_{V-\{v\}}(\mb{a}), \, p_{V-\{v\}}(\mb{b}))$ if either $x = \mb{a}(v)$ or $v$ is an isolated vertex,
    \item $(G \wedge vw \setminus v, \, p_{V-\{v\}}(\mb{a}[V-\{v,w\}] + \mb{b}[\{v,w\}]), \, p_{V-\{v\}}(\mb{a}[\{v,w\}] + \mb{b}[V-\{v,w\}]))$ if $x = \mb{b}(v)$ and $w$ is a neighbor of $v$, and
    \item $(G*v\setminus v, \, p_{V-\{v\}}(\mb{a}), \, p_{V-\{v\}}(\mb{a}[N_G(v)] + \mb{b}))$ otherwise.
  \end{enumerate}
\end{proposition}

\begin{corollary}\label{cor:ess1}
  Let $G$ be a fundamental graph of an isotropic system $S$ with at least five vertices.
  A vertex of $S$ is non-essential in $S$ if and only if it is non-essential in $G$.
  \qed
\end{corollary}

\section{Triangles in 3-connected isotropic systems}\label{sec:tri}

Let $S = (V,L)$ be an isotropic system.
A \emph{triangle} in $S$ is a vector in $L$ such that the size of its support is $3$.
Let $H(S)$ be the $3$-uniform hypergraph on $V$ whose edge set is the set of supports of triangles in~$S$.
First, we present several lemmas of Allys~\cite{Allys1994iso} which show the existence of triangles whose supports contain some essential vertices in a $3$-connected isotropic system.

\begin{lemma}[Allys~{\cite[Lemma~3.3]{Allys1994iso}}]\label{lem:tri0}
  Let $S = (V,L)$ be a $3$-connected isotropic system with $|V| \geq 4$.
  If $\mbt$ and~$\mbt'$ are triangles in $S$, then one of the following holds.
  \begin{enumerate}[label=\rm(\roman*)]
    \item $\supp(\mbt)$ and $\supp(\mbt')$ are disjoint.
    \item $\supp(\mbt) \cap \supp(\mbt') = \{v\}$ and $\mbt(v) = \mbt'(v)$ for some $v\in V$.
    \item $\supp(\mbt) \cap \supp(\mbt') = \{v,w\}$, $\mbt(v) \neq \mbt'(v)$, and $\mbt(w) \neq \mbt'(w)$ for some $v,w \in V$.
    \item $\mbt = \mbt'$.
  \end{enumerate}
\end{lemma}

By Lemma~\ref{lem:tri0}, in a $3$-connected isotropic system $S$ with at least $4$ vertices, triangles have distinct supports, and thus there is a bijection from the set of triangles of $S$ to the set of edges of $H(S)$.
Now we investigate what vertices of a tight path in $H(S)$ are essential or non-essential in $S$.
Recall that a vertex $v$ of $S$ is essential if at least two of $S|^v_\alpha$, $S|^v_\beta$, $S|^v_\gamma$ are $3$-connected.

\begin{lemma}\label{lem:not3conn}
  Let $S$ be an isotropic system with at least $5$ vertices, and let $\mbt$ be a triangle in $S$.
  Then for each $v$ in the support of $\mbt$, $S|^v_{\mbt(v)}$ is not $3$-connected.
\end{lemma}
\begin{proof}
  A minor $S|^v_{\mbt(v)}$ has a nonzero vector $p_{V-\{v\}}(\mbt)$ whose support has size $2$.
  By Lemma~\ref{lem:pt}, $S|^v_{\mbt(v)}$ is not $3$-connected.
\end{proof}

\begin{lemma}\label{lem:intess}
  Let $S$ be a $3$-connected isotropic system with at least $5$ vertices.
  Every internal vertex of a tight path in $H(S)$ is essential in $S$.
\end{lemma}
\begin{proof}
  It is enough to prove that if $\mbt$ and $\mbt'$ are triangles in $S$ such that $\supp(\mbt) \cap \supp(\mbt') = \{v,w\}$, then $v$ is essential in $S$.
  By Lemma~\ref{lem:not3conn}, neither $S|^v_{\mbt(v)}$ nor $S|^v_{\mbt'(v)}$ is $3$-connected.
  By Lemma~\ref{lem:tri0}, $\mbt(v) \neq \mbt'(v)$ and, therefore, $v$ is essential in $S$.
\end{proof}  

\begin{lemma}\label{lem:codeg3}
  Let $S$ be a $3$-connected isotropic system.
  For distinct vertices $u$ and $v$ of $H(S)$, there are at most three edges of $H(S)$ incident with both $u$ and $v$.
\end{lemma}
\begin{proof}
  Suppose that there are four distinct triangles $\mbt_1$, $\mbt_2$, $\mbt_3$, and $\mbt_4$ in $S$ whose supports contain both~$u$ and~$v$.
  Then $\mbt_k(u) \in K - \{0\} = \{\alpha,\beta,\gamma\}$ for each $1\leq k\leq 4$.
  So $\mbt_i(u) = \mbt_j(u)$ for some distinct~$i$ and~$j$.
  By Lemma~\ref{lem:tri0}, $\mbt_i = \mbt_j$, which is a contradiction.
\end{proof}

\begin{lemma}\label{lem:codeg2}
  Let $S$ be a $3$-connected isotropic system.
  Let $\{u,v,w_1\}$ and $\{u,v,w_2\}$ be distinct edges in $H(S)$.
  If $e$ is an edge incident with $u$ in $H(S)$, then $e$ is incident with at least one of $v$, $w_1$, and~$w_2$.
\end{lemma}
\begin{proof}
  Let $\mbt_1$ and $\mbt_2$ be triangles in $S$ whose supports are $\{u,v,w_1\}$ and $\{u,v,w_2\}$, respectively.
  Let $\mbt$ be a triangle whose support is $e$.
  By Lemma~\ref{lem:tri0}, $\mbt_1(u) \neq \mbt_2(u)$.
  Without loss of generality, $\mbt(u) \neq \mbt_1(u)$.
  By Lemma~\ref{lem:tri0}, $|\supp(\mbt) \cap \supp(\mbt_1)| = 2$, so $e$ is incident with $v$ or $w_1$.
\end{proof}

\begin{lemma}[Allys~{\cite[Lemma~3.5]{Allys1994iso}}]\label{lem:tri1}
  Let $S = (V,L)$ be a $3$-connected isotropic system with at least $4$ vertices.
  For a vertex $v \in V$ and two distinct $x, y \in K - \{0\}$, if neither $S|_x^v$ nor $S|_y^v$ is $3$-connected, then $S$ has a triangle $\mbt$ such that $\mbt(v) \in \{x,y\}$.
\end{lemma}

Lemma~\ref{lem:tri1} implies that for a $3$-connected isotropic system $S$ with at least $4$ vertices, if a vertex is essential, then it is incident with an edge of $H(S)$.

\begin{lemma}[Allys~{\cite[Lemma~3.5]{Allys1994iso}}\footnote{In~\cite{Allys1994iso}, there are two Lemmas~3.5, and this is the second Lemma~3.5.}]\label{lem:tri2}
  Let $\mbt$ be a triangle in a $3$-connected isotropic system $S$ with at least $4$ vertices, where $\supp(\mbt) = \{u,v,w\}$.
  For $x \in K - \{0, \mbt(u)\}$ and $y \in K - \{0, \mbt(w)\}$, if neither $S|_{x}^u$ nor $S|_{y}^w$ is $3$-connected, then there are triangles $\mbt_1$ and $\mbt_2$ (possibly $\mbt_1 = \mbt_2$) such that $\mbt_1(u) = x$, $\mbt_2(w) = y$, and $\mbt_1(v) = \mbt_2(v)$.
\end{lemma}

\begin{lemma}\label{lem:extpath1}
  For a $3$-connected isotropic system $S$ with at least $4$ vertices, if an edge $e$ of $H(S)$ is incident with at least two essential vertices in $S$, then $H(S)$ has an edge $e'$ such that $|e \cap e'|=2$.
\end{lemma}
\begin{proof}
  Let $u, w \in e$ be distinct essential vertices in $S$ and let $\mbt$ be a triangle in $S$ such that $\supp(\mbt) = e$.
  Because $u$ and $w$ are essential, we have $x \in K - \{0,\mbt(u)\}$ and $y \in K - \{0,\mbt(w)\}$ such that neither $S|^u_x$ nor $S|^w_y$ is $3$-connected.
  By Lemma~\ref{lem:tri2}, $S$ has a triangle $\mbt_1$ such that $\mbt_1(u) = x$.
  Then an edge $e' := \supp(\mbt_1)$ of $H(S)$ satisfies that $|e \cap e'| = 2$ by Lemma~\ref{lem:tri0}.
\end{proof}

Lemma~\ref{lem:extpath1} provides a sufficient condition for extending a tight path of length $1$.
Now we aim to prove that two ends of a maximal tight path in $H(S)$ are non-essential unless $S$ is cyclic.

\begin{figure}
  \centering
  \begin{tikzpicture}
    \foreach \x [count=\p] in {0,...,4} {
      \draw[fill=gray, opacity=0.5]
      (\x*360/5+18:1.7)
      .. controls (\x*360/5+72:2.8) and (\x*360/5+72+36:2.8) ..
      (\x*360/5+144+18:1.7)
      .. controls (\x*360/5+72+54:1.9) and (\x*360/5+72+18:2.4) ..
      (\x*360/5+72+18:1.7)
      .. controls (\x*360/5+72+18:2.4) and (\x*360/5+72-18:1.9) ..
      (\x*360/5+18:1.7) -- cycle;
    };
    \foreach \x [count=\p] in {0,...,4} {
     \draw[fill=gray, opacity=0.5]
      (\x*360/5+18:1.7)
      --
      (\x*360/5+72+18:1.7)
      --
      (\x*360/5+216+18:1.7)
      .. controls (\x*360/5+89:1.35) ..
      (\x*360/5+18:1.7) -- cycle;
    };
    \foreach \x [count=\p] in {0,...,4} {
      \node[shape=circle,fill=black,scale=0.6] (\p) at (\x*360/5 + 18:1.7) {};
    };
    \node () at (-90:2.5) {$C_5$};
  \end{tikzpicture} 
  \hspace{0.5cm}
  \begin{tikzpicture}
    \foreach \x [count=\p] in {0,...,5} {
      \draw[fill=gray, opacity=0.5]
      (\x*360/6:1.7)
      .. controls (\x*360/6+60-15:2.8) and (\x*360/6+60+15:2.8) ..
      (\x*360/6+120:1.7)
      .. controls (\x*360/6+60+30:1.9) and (\x*360/6+60+5:2.6) ..
      (\x*360/6+60:1.7)
      .. controls (\x*360/6+60-5:2.6) and (\x*360/6+60-30:1.9) ..
      (\x*360/6:1.7) -- cycle;
    };
    \foreach \x [count=\p] in {0,1} {
     \draw[fill=gray, opacity=0.5]
     (\x*360/6:1.7)
     .. controls (\x*360/6+60:0.2) ..
     (\x*360/6+120:1.7)
     .. controls (\x*360/6+180:0.2) ..
     (\x*360/6+240:1.7)
     .. controls (\x*360/6+300:0.2) ..
     (\x*360/6:1.7) -- cycle;
    };
    \foreach \x [count=\p] in {0,...,5} {
      \node[shape=circle,fill=black,scale=0.6] (\p) at (\x*360/6:1.7) {};
    };
    \node () at (-90:2.5) {$C_6$};
  \end{tikzpicture}    
  \caption{Illustrations of $H(S)$ for an isotropic system $S$ whose fundamental graph is $C_5$ or $C_6$.}
  \label{fig:c5c6}
\end{figure}
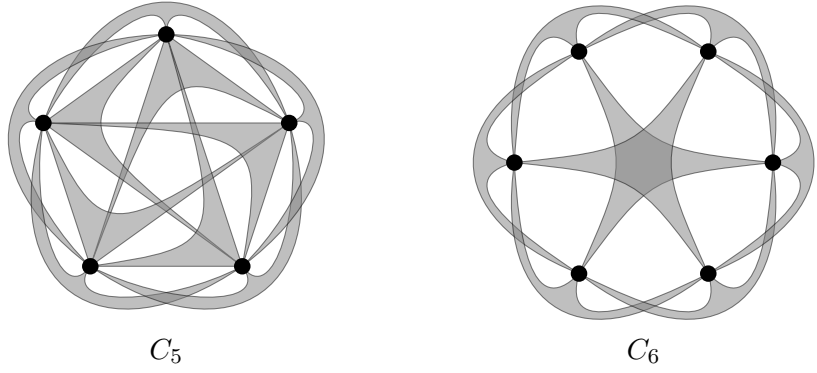

In the next two lemmas, we show that under some assumptions, 
no internal vertex of a tight path of length at least $3$ in $H(S)$
is incident with edges not on the path.
By Figure~\ref{fig:c5c6}, it is necessary to require that 
$S$ does not have $C_5$ or $C_6$ as a fundamental graph.
Our proof of the following lemma is motivated by the proof of Claim~3 in Allys~\cite[Theorem~4.3]{Allys1994iso} proving a weaker statement.

\begin{lemma}\label{lem:v2}
  Let $S = (V,L)$ be a $3$-connected isotropic system and let $P = v_0v_1v_2v_3v_4$ be a tight path in $H(S)$.
  If neither $C_5$ nor $C_6$ is a fundamental graph of $S$, then $P$ contains every edge of $H(S)$ incident with $v_2$. 
\end{lemma}
\begin{proof}
  By Lemma~\ref{lem:c5}, we may assume that $\abs{V}\ge 6$.
  For each $1\leq i\leq 3$, let $\mbt_i$ be a triangle whose support is $\{v_{i-1},v_i,v_{i+1}\}$.
  By Lemma~\ref{lem:tri0} applied to $\mbt_{1}$ and $\mbt_{3}$, we have $\mbt_{1}(v_2) = \mbt_{3}(v_2)$.
  By Lemma~\ref{lem:tri0} applied to $\mbt_{1}$ and $\mbt_{2}$, we have $\mbt_{1}(v_1) \neq \mbt_{2}(v_1)$
  and $\mbt_{1}(v_2) \neq \mbt_{2}(v_2)$.
  Also by applying Lemma~\ref{lem:tri0} to $\mbt_2$ and $\mbt_3$, 
  we have 
  $\mbt_{2}(v_2) \neq \mbt_{3}(v_2)$
  and $\mbt_{2}(v_3) \neq \mbt_{3}(v_3)$.

  Suppose for contradiction that $S$ has a triangle $\mbt$ such that $v_2\in \supp(\mbt)$ and $\mbt \neq \mbt_i$ for all $i \in \{1,2,3\}$.

  We first claim that $\mbt(v_2) = \mbt_2(v_2)$.
  Suppose that   $\mbt(v_2) \neq  \mbt_2(v_2)$.
  By Lemma~\ref{lem:tri0},
  $\abs{\supp(\mbt)\cap \supp(\mbt_2)}=2$.
  Therefore there is a unique $j\in \{1,3\}$ such that $v_j\in \supp(\mbt)\cap \supp(\mbt_2)$ and furthermore $\mbt_2(v_j)\neq \mbt(v_j)$.
  By reversing the path if necessary, we may assume that $j=3$ and $v_1\notin \supp(\mbt)$.
  Since $v_2,v_3\in \supp(t)\cap \supp(\mbt_3)$, 
  by applying Lemma~\ref{lem:tri0}, 
  we deduce that 
  $\mbt(v_3)\neq \mbt_3(v_3)$ and $\mbt(v_2)\neq \mbt_3(v_2)$.
  Since $\mbt_3(v_2)\neq \mbt_2(v_2)$, 
  we deduce that $\mbt(v_2)=\mbt_2(v_2)+\mbt_3(v_2)
  =\mbt_2(v_2)+\mbt_1(v_2)\neq \mbt_1(v_2)$.
  By Lemma~\ref{lem:tri0}, $\abs{\supp(\mbt)\cap \supp(t_1)}=2$ and therefore $v_0\in \supp(\mbt)$. Thus we deduce that $\supp(\mbt)=\{v_0,v_2,v_3\}$.
  As $\mbt_1(v_3)=0$ and $\mbt_2(v_3)\neq \mbt(v_3)$, 
  we deduce that $\mbt_1$, $\mbt_2$, $\mbt$ are linearly independent, so $c(\{v_0,v_1,v_2,v_3\}) = 4- \dim(L|_{\subseteq \{v_0,v_1,v_2,v_3\}}) \leq 1$, where $c$ is the connectivity function of $S$.
  This contracts to the assumption that 
  $S$ is $3$-connected.
  Therefore, $\mbt(v_2) = \mbt_2(v_2)$.

  By Lemma~\ref{lem:tri0} for $\mbt$ and $\mbt_2$, the support of $\mbt$ contains neither $v_1$ nor $v_3$.
  Since $\mbt(v_2)\neq \mbt_1(v_2)$, 
  by Lemma~\ref{lem:tri0},  $\supp(\mbt)$ contains $v_0$
  and $\mbt(v_0)\neq \mbt_1(v_0)$.
  Similarly, as $\mbt(v_2)\neq\mbt_3(v_2)$, the support of $\mbt$ contains $v_4$ and $\mbt(v_4)\neq \mbt_3(v_4)$.
  Hence $\supp(\mbt) = \{v_0,v_2,v_4\}$ 
  and $\mbt_1$, $\mbt_2$, $\mbt_3$, $\mbt$ are linearly independent, so $c(\{v_0,v_1,v_2,v_3,v_4\}) \leq 5 - \dim(L|_{\subseteq \{v_0,v_1,v_2,v_3,v_4\}}) = 1$.
  Since~$S$ is $3$-connected, we have $|V - \{v_0,v_1,v_2,v_3,v_4\}| \leq 1$ and therefore $|V| = 6$.
  Let $v_5$ denote the vertex of~$V$ other than $v_0, \dots, v_4$.
  Because of $\mbt_1$, we have $\dim(L|_{\subseteq\{v_0,v_1,v_2\}}) \geq 1$.
  By Proposition~\ref{prop:conn0}\ref{item:c02}, $c(\{v_3,v_4,v_5\}) = c(\{v_0,v_1,v_2\})$
  and thus $\dim(L|_{\subseteq\{v_3,v_4,v_5\}}) =  \dim(L|_{\subseteq\{v_0,v_1,v_2\}}) \geq 1$.
  Then $S$ has a nonzero vector $\mbt_4$ such that $\supp(\mbt_4) \subseteq \{v_3,v_4,v_5\}$.
  By Lemma~\ref{lem:pt}, $\mbt_4$ is a triangle whose support is $\{v_3,v_4,v_5\}$.
  Similarly, because of $\mbt_2$ and $\mbt_3$, there are triangles $\mbt_5$ and $\mbt_0$ whose supports are $\{v_4,v_5,v_0\}$ and $\{v_5,v_0,v_1\}$, respectively.
  Let $v_{-1} := v_5$, $v_6 := v_0$, $\mbt_{6} := \mbt_0$, and $\mbt_7 := \mbt_1$.
  Let $C_6$ be the cycle graph on $\{v_0,v_1,\dots,v_5\}$ in this order and let $\mb{a}$ and $\mb{b}$ be vectors in $K^V$ such that $\mb{a}(v_i) = \mbt_{i+1}(v_i)$ and $\mb{b}(v_i) = \mbt_{i}(v_i)$ for each $0\leq i\leq 5$.
  Then, for every $0 \le i \le 5$,
  \begin{itemize}
    \item $\mb{a}$ and $\mb{b}$ are supplementary by Lemma~\ref{lem:tri0},
    \item $\mbt_i(v_{i+1}) = \mbt_{i+2}(v_{i+1}) = \mb{a}(v_{i+1})$ by Lemma~\ref{lem:tri0},
    \item $\mbt_i(v_{i-1}) = \mb{a}(v_{i-1})$ by the definition of $\mb{a}$,
    \item $\mbt_i(v_i) = \mb{b}(v_i)$ by the definition of $\mb{b}$, and
    \item $\mbt_i(w) = 0$ for all $w \in V \setminus \supp(\mbt_i)$.
  \end{itemize}
  Then $\mbt_i = \mb{a}[N_{C_6}(v_i)] + \mb{b}[\{v_i\}]$ for each~$i$.
  Therefore $(C_6,\mb{a},\mb{b})$ is a graphic presentation of $S$ and so $C_6$ is a fundamental graph of $S$, contradicting the assumption.
\end{proof}

\begin{lemma}\label{lem:v1}
  Let $S$ be a $3$-connected isotropic system.
  Let $P = v_0v_1v_2v_3v_4v_5$ be a tight path in $H(S)$.
  If %
  $v_0$ is non-essential in $S$, then $P$ contains every edge of $H(S)$ incident with $v_1$.
\end{lemma}
\begin{proof}
  Suppose that there is an edge $e$ of $H(S)$ incident with $v_1$ and not in $P$.
  By Lemma~\ref{lem:codeg2}, $e$ is incident with at least one of $v_0$, $v_2$, and $v_3$.
  Since $v_0$ is non-essential in $S$, 
  $C_6$ is not a fundamental graph of $S$ 
  and therefore neither $v_2$ nor $v_3$ is incident with $e$ by Lemma~\ref{lem:v2}. Thus $e$ is incident with both $v_0$ and $v_1$.
  Let $v\neq v_0,v_1$ be a vertex incident with $e$. 
  Then $vv_0v_1v_2$ is a tight path and by Lemma~\ref{lem:intess}, $v_0$ is essential in $S$, contradicting the assumption. 
\end{proof}

By Lemmas~\ref{lem:v2} and~\ref{lem:v1}, we deduce the following.
\begin{proposition}\label{prop:v12}
  For a $3$-connected isotropic system $S$, if $H(S)$ has a tight path $P = v_0v_1\cdots v_{k+1}$ of length $k\geq 4$ such that $v_0$ and $v_{k+1}$ are non-essential in $S$, then
  $P$ contains every edge of $H(S)$ incident with at least one of $v_1,v_2,\dots,v_k$. \qed
\end{proposition}

\begin{lemma}\label{lem:lx}
  Let $S = (V,L)$ be a $3$-connected isotropic system.
  If $X$ is a subset of $V$ such that $\min\{|X|, |V-X|\} \geq 2$ and $\dim(L|_{\subseteq X}) \geq |X|-2$, then $\dim(L|_X) = |X|+2$.
\end{lemma}
\begin{proof}
  Recall that by Lemma~\ref{lem:ortho}, $\dim(L|_{\subseteq X}) + \dim(L|_{X}) = \dim(K^X) = 2|X|$ and $c(X) = |X| - \dim(L|_{\subseteq X}) = \dim(L|_X) - |X|$.
  Hence $\dim(L|_{X}) = 2|X| - \dim(L|_{\subseteq X}) \leq |X|+2$.
  Since $S$ is $3$-connected and $\min\{|X|, |V-X|\} \geq 2$, we have $c(X) \geq 2$ and so $\dim(L|_{X}) = |X| + c(X) \geq |X|+2$.
  Therefore, $\dim(L|_{X}) = |X|+2$.
\end{proof}

\begin{lemma}\label{lem:cyc}
  Let $S = (V,L)$ be a $3$-connected isotropic system.
  Let $v_0v_1\cdots v_{k+1}$ be a tight path of length $k\geq 3$ in $H(S)$.
  If $\{v_k,v_{k+1},v_0\}$ is an edge of $H(S)$, then $C_{k+2}$ is a fundamental graph of $S$.
\end{lemma}
\begin{proof}
  For each $1\leq i\leq k$, let $\mbt_i$ be a triangle in $S$ whose support is $\{v_{i-1},v_i,v_{i+1}\}$.
  Let $\mbt_{k+1}$ be a triangle in $S$ whose support is $\{v_k,v_{k+1},v_0\}$.
  By Lemma~\ref{lem:tri0}, we deduce that $\mbt_{k+1}(v_0) = \mbt_{1}(v_0)$ and $\mbt_{i-1}(v_i) = \mbt_{i+1}(v_i)$ for every $2 \leq i \leq k$.
  Also by Lemma~\ref{lem:tri0}, we have that $\mbt_i(v_i) \neq \mbt_{i+1}(v_i)$ and $\mbt_i(v_{i+1}) \neq \mbt_{i+1}(v_{i+1})$ for all $1\leq i\leq k$.

  Let $X = \{v_2,v_3,\dots,v_k\}$.
  We claim that $p_X(\mbt_1), p_X(\mbt_2), \dots, p_X(\mbt_{k+1})$ are linearly independent.
  Suppose that $\sum_{i=1}^{k+1}c_i p_X(\mbt_i)=0$ for some $c_1,c_2,\ldots,c_{k+1}\in \GF(2)$.
  For $1<j< k+1$,
  $\sum_{i=1}^{k+1} c_i \mbt_i(v_{j}) = c_{j-1} \mbt_{j-1}(v_j)+c_j \mbt_j (v_{j}) +c_{j+1} \mbt_{j+1}(v_j) 
  =c_j \mbt_j (v_{j}) +(c_{j-1}+c_{j+1}) \mbt_{j+1}(v_j) $
  and thus $c_j=0$ 
  because $\mbt_j(v_j)$ and $\mbt_{j+1}(v_j)$ are linearly independent in $K$.
  So $c_1 p_X(\mbt_1)+c_{k+1}p_X(\mbt_{k+1})=0$. Since $0=c_1 \mbt_1(v_2)+c_{k+1}\mbt_{k+1}(v_2)=c_1 \mbt_1(v_2)$, we deduce that $c_1=0$ and so $c_{k+1}=0$. 
  Therefore $p_X(\mbt_1), p_X(\mbt_2), \dots, p_X(\mbt_{k+1})$ are linearly independent.
  This also implies that 
  $\mbt_1, \mbt_2, \dots, \mbt_{k+1}$ are linearly independent.

  Hence $c(\{v_0, v_1, \dots, v_{k+1}\}) = k+2 - \dim(L|_{\subseteq \{v_0,v_1,\dots,v_{k+1}\}}) \leq 1$.
  Since $S$ is $3$-connected, $|V - \{v_0,v_1, \dots, v_{k+1}\}| \leq 1$.
  Hence $|V| = k+3$ or $k+2$.

  We have $\dim(L|_{\subseteq X}) \geq k-3 = |X|-2$ because $\mbt_3, \mbt_4, \dots, \mbt_{k-1}$ are linearly independent.
  By Lemma~\ref{lem:lx}, $\dim(L|_X) = |X|+2 = k+1$.
  Therefore  
  $p_X(\mbt_1), p_X(\mbt_2), \dots, p_X(\mbt_{k+1})$  form a basis of~$L|_X$.

  Suppose $|V| = k+3$.
  Let $w$ be the vertex of $V$ other than $v_0,v_1,\dots,v_{k+1}$.
  Let $\mb{a}$ and $\mb{b}$ be vectors in $L$ such that $\{\mbt_1, \dots, \mbt_{k+1}, \mb{a}, \mb{b}\}$ is a basis of $L$.
  Since $\{ p_X(\mbt_1), p_X(\mbt_2), \dots, p_X(\mbt_{k+1}) \}$ is a basis of $L|_X$, we may assume that $p_X(\mb{a}) = 0$ and $p_X(\mb{b}) = 0$.
  Hence the supports of $\mb{a}$ and~$\mb{b}$ are subsets of $V-X = \{v_0,v_1,v_{k+1},w\}$.
  Then $\supp(\mb{a}) \cap \supp(\mbt_2) \subseteq \{v_1\}$, so $0 = \inn{\mb{a}}{\mbt_2} = \innk{\mb{a}(v_1)}{\mbt_2(v_1)}$.
  It implies that $\mb{a}(v_1) \in \{0, \mbt_2(v_1)\}$ and similarly $\mb{b}(v_1) \in \{0, \mbt_2(v_1)\}$.
  Then one of $\mb{a}$, $\mb{b}$, and $\mb{a}+\mb{b}$, say $\mb{c}$, satisfies $\mb{c}(v_1) = 0$.
  Then $\supp(\mb{c}) \subseteq \{v_0,v_{k+1},w\}$.
  Since $S$ is $3$-connected and $\mb{c}$ is a nonzero vector, $\supp(\mb{c}) = \{v_0,v_{k+1},w\}$ by Lemma~\ref{lem:pt}.
  Then $\supp(\mb{c}) \cap \supp(\mbt_{k+1}) = \{v_0,v_{k+1}\}$.
  Since $v_0 \in \supp(\mbt_1)$, by Lemma~\ref{lem:codeg2}, the support of $\mbt_1$ contains $v_{k}$, $v_{k+1}$, or $w$, contradicting that $\supp(\mbt_1) = \{v_0,v_1,v_2\}$.
  Therefore, $|V| = k+2$.
  
  Let $\mbt_0$ be a vector in $L$ such that $\{\mbt_0,\mbt_1,\dots,\mbt_{k+1}\}$ is a basis of $S$.
  Since $\{ p_X(\mbt_1), \dots, p_X(\mbt_{k+1}) \}$ is a basis of $L|_X$, 
  we may assume that $p_X(\mbt_0)=0$
  and therefore 
  the support of $\mbt_0$ is a subset of $V-X = \{v_{k+1},v_0,v_1\}$.
  Since $S$ is $3$-connected and $\mbt_0$ is nonzero, 
  by Lemma~\ref{lem:pt}, $\supp(\mbt_0) = \{v_{k+1},v_0,v_1\}$.
  Let $\mb{a}, \mb{b}\in K^V$  be vectors 
  such that 
  $\mb{a}(v_i) = \mbt_{i+1}(v_i)$
  and 
  $\mb{b}(v_i)=\mbt_i(v_i)$
  for all $0\leq i\leq k+1$, 
  where $\mbt_{k+2}:=\mbt_0$.
  Let $C_{k+2}$ be the cycle graph on $\{v_0,v_1,\ldots,v_{k+1}\}$ in this order. 
  Let $v_{-1}:=v_{k+1}$, $v_{k+2}:=v_0$, and $\mbt_{k+3}:=\mbt_1$. Then, for all $i\in\{0,1,2,\ldots,k+1\}$,
  \begin{itemize}
    \item $\innk{\mb a(v_i)}{\mb b(v_i)}=\innk{\mbt_{i+1}(v_i)}{\mbt_i(v_i)}=1$  by Lemma~\ref{lem:tri0},
  \item   $\mbt_i(v_{i+1})=\mbt_{i+2}(v_{i+1})=\mb a(v_{i+1})$   by Lemma~\ref{lem:tri0},
  \item     $\mbt_i(v_{i-1})=\mb a(v_{i-1})$ by the definition of $\mb a$,
  \item $\mbt_i (v_i)= \mb b(v_i)$ by the definition of $\mb b$, and 
  \item $\mbt_i(u)=0$ for all $u\in V$ with $u\neq v_{i-1},v_i,v_{i+1}$ because $\supp(\mbt_i) = \{v_{i-1},v_i,v_{i+1}\}$.
  \end{itemize}
  Thus $\mbt_i=\mb a[N_{C_{k+2}}(v_i)]+\mb b[\{v_i\}]$ for all $i$ and therefore 
  $(C_{k+2},\mb a,\mb b)$ is a graphic presentation of $S$. So $C_{k+2}$ is a fundamental graph of $S$.
\end{proof}

The following lemma provides a sufficient condition to extend a tight path.

\begin{lemma}\label{lem:extpath}
  Let $S = (V,L)$ be a $3$-connected isotropic system.
  Let $v_0v_1\dots v_{k+1}$ be a tight path of length $k \geq 2$ in $H(S)$.
  If $S$ is not cyclic and $v_{k+1}$ is essential in $S$, then $H(S)$ has a vertex $v_{k+2}$ such that
  \begin{enumerate}[label=\rm(\roman*)]
      \item $v_0v_1v_2v_3v_4$ or $v_0v_2v_1v_3v_4$ is a tight path in $H(S)$ if $k=2$, and
      \item $v_0v_1\dots v_{k+1}v_{k+2}$ is a tight path in $H(S)$ if $k \geq 3$.
  \end{enumerate}
\end{lemma}
\begin{proof}
  Observe that $|V| \geq 6$ 
  by Lemmas~\ref{lem:c4} and~\ref{lem:c5}
  because $S$ is $3$-connected and not cyclic and $\abs{V}\ge k+2\ge 4$.

  Let $\mb{a}$ and $\mb{b}$ be triangles in $S$ whose supports are $\{v_{k-2},v_{k-1},v_{k}\}$ and $\{v_{k-1},v_{k},v_{k+1}\}$, respectively.
  By Lemma~\ref{lem:tri0}, $\mb{a}(v_{k-1}) \neq \mb{b}(v_{k-1})$ and $\mb{a}(v_{k}) \neq \mb{b}(v_{k})$.

  Since $v_{k+1}$ is essential, there is $x \in K-\{0,\mb{b}(v_{k+1})\}$ such that $S|_x^{v_{k+1}}$ is not $3$-connected.
  By Lemma~\ref{lem:pt}, $S|^{v_k}_{\mb{a}(v_k)}$ is not $3$-connected because the support of $p_{V-\{v_k\}}(\mb{a})$ has size $2$.
  Applying Lemma~\ref{lem:tri2} for $\mb{b}$, $S|^{v_k}_{\mb{a}(v_k)}$, and $S|^{v_{k+1}}_{x}$, we obtain a triangle $\mb{c}$ in~$S$ such that
  \begin{equation*}
    \mb{c}(v_{k+1}) = x \neq \mb{b}(v_{k+1}).
  \end{equation*}
  By Lemma~\ref{lem:tri0} for $\mb{b}$ and $\mb{c}$, the support of $\mb{c}$ contains exactly one of $v_{k-1}$ and $v_{k}$.

  As $S$ is $3$-connected, 
  $c(\{v_{k-2},v_{k-1},v_{k},v_{k+1}\}) \geq 2$ and therefore $\dim(L|_{\subseteq \{v_{k-2},v_{k-1},v_{k},v_{k+1}\}}) \leq 2$.
  Observe that $\mb{a}$, $\mb{b}$, and $\mb{c}$ are linearly independent. Since the supports of $\mb{a}$ and $\mb{b}$ are subsets of $\{v_{k-2},v_{k-1},v_{k},v_{k+1}\}$, the support of $\mb{c}$ is not a subset of $\{v_{k-2},v_{k-1},v_{k},v_{k+1}\}$
  because otherwise $L|_{\subseteq \{v_{k-2},v_{k-1},v_{k},v_{k+1}\}}$ contains three linearly independent vectors.
  Thus $v_{k-2} \not\in \supp(\mb{c})$,
  because $\abs{\supp(\mb{c})\cap \{v_{k-1},v_k,v_{k+1}\}}=2$.
  If $k=2$, then $\mb{c} = \{v_2,v_3,v_4\}$ or $\{v_1,v_3,v_4\}$ for some $v_4 \in V \setminus \{v_0,v_1,v_2,v_3\}$ and therefore $v_0v_1v_2v_3v_4$ or $v_0v_2v_1v_3v_4$ is a tight path in $H(S)$.
  Hence we may assume that $k\geq 3$.

  Since neither $C_5$ nor $C_6$ is a fundamental graph of $S$, 
  by Lemma~\ref{lem:v2}, 
  a tight path $v_{k-3}v_{k-2}v_{k-1}v_{k}v_{k+1}$ 
  contains every edge of $H(S)$ incident with $v_{k-1}$
  and so $v_{k-1}\notin\supp(\mb c)$.
  This implies that $v_k,v_{k+1}\in \supp(\mb c)$.
  For each $i\in\{0,\dots,k-3\}$, as a cycle $C_{k-i+2}$ is not a fundamental graph of $S$, by Lemma~\ref{lem:cyc} applied to a tight path $v_iv_{i+1} \cdots v_{k} v_{k+1}$, a set $\{v_i,v_k,v_{k+1}\}$ is not an edge of $H(S)$.
  Hence none of $v_0$, $v_1$, $\ldots$, $v_{k-3}$ is in the support of $\mb{c}$.
  Recall that $v_{k-2} \notin \supp(\mb{c})$.
  Therefore $\supp(\mb{c}) = \{v_k,v_{k+1},v_{k+2}\}$ for some $v_{k+2} \in V \setminus \{v_0,\dots,v_{k+1}\}$ and $v_0v_1v_2\cdots v_{k+1}v_{k+2}$ is a tight path in $S$.
\end{proof}

For a $3$-connected isotropic system $S$ with at least five vertices, Lemma~\ref{lem:intess} states that if $v_0v_1\dots v_{k+1}$ is a tight path of length $k\geq 2$ in $H(S)$, then $v_1,v_2,\dots,v_k$ are essential in $S$.
The following proposition provides a feature of ends of a maximal tight path when $S$ is not cyclic.
Recall that a tight path has exactly two ends if its length is at least $2$.

\begin{proposition}\label{prop:maxpath}
  Let $S$ be a $3$-connected isotropic system with at least $5$ vertices.
  If $S$ is not cyclic, then at least two ends of a maximal tight path in $H(S)$ are non-essential in $S$.
\end{proposition}
\begin{proof}
  Let $v_0v_1\dots v_{k+1}$ be a maximal tight path in $H(S)$.
  If $k\geq 2$, then $v_0$ and $v_{k+1}$ are non-essential in $S$ by Lemma~\ref{lem:extpath}.
  If $k=1$, then at least two of $v_0, v_1, v_2$ are non-essential in $S$ by Lemma~\ref{lem:extpath1}.
\end{proof}

By Proposition~\ref{prop:igconn} and Lemma~\ref{lem:pt}, for a prime graph $G$ with at least five vertices and $X \subseteq V(G)$, a subset $X$ has size $3$ and $|\rho_G(X)| = 2$ if and only if $S$ has a triangle $\mbt$ with $\supp(\mbt) = X$, where $S$ is an isotropic system having $G$ as its fundamental graph.
Therefore, Proposition~\ref{prop:maxpath} and Theorem~\ref{thm:maxpathG} are equivalent.

It is straightforward to prove Theorem~\ref{thm:noness2} from Proposition~\ref{prop:maxpath}.

\begin{proof}[Proof of Theorem~\ref{thm:noness2}]
  Let $G = (V,E)$ be a prime graph with at least four vertices which is not locally equivalent to a cycle graph.
  As no graph on four vertices is prime, $|V| \ge 5$.
  Let $\mb{a}$ and $\mb{b}$ be supplementary vectors in $K^V$.
  Let $S$ be an isotropic system having a graphic presentation $(G,\mb{a},\mb{b})$.
  Then $S$ is not cyclic because all fundamental graphs of $S$ are locally equivalent.
  By Corollary~\ref{cor:igconn2}, $S$ is $3$-connected.
  By Corollary~\ref{cor:ess1}, $v\in V$ is non-essential in $G$ if and only if it is non-essential in $S$.
  Therefore, it suffices to show that $S$ has at least two non-essential vertices.
  
  We may assume that $S$ has an essential vertex.
  Then by Lemma~\ref{lem:tri1}, $S$ has a triangle and so $H(S)$ has a maximal tight path.
  By Proposition~\ref{prop:maxpath}, at least two ends of the maximal tight path are non-essential, and therefore $S$ has at least two non-essential vertices.
\end{proof}

In the remainder, we describe a structure of a $3$-connected isotropic system in terms of triangles, which will be a major ingredient to prove Theorem~\ref{thm:noness3}.
We first define $5$ types of partial hypergraphs.
Let $H = (V,E)$ be a $3$-uniform hypergraph and $N$ be a subset of $V$.
\begin{figure}
  \centering  
    \begin{tikzpicture}
      \node (dum) at (2,-0.8) {};
      \node (v0) at (0,0) {};
      \node (v7) at (4,1.25) {};
      \node (v1) at (-0.1,1.2) {};
      \node (v2) at (0.7,0.8) {};
      \node (v3) at (0.9,1.6) {};
      \node (v4) at (1.7,1.2) {};
      \node (v5) at (2.3,1.8) {};
      \node (v6) at (2.7,1.1) {};

      \begin{scope}
        \draw[fill=gray!70, opacity=0.5] ($(v0)+(-0.8,0)$) 
        to[out=90,in=180] ($(dum) + (-0.3,0.5)$) 
        to[out=0,in=180] ($(v7) + (0.3,0.8)$)
        to[out=0,in=0] ($(dum) + (0.5,-0.8)$)
        to[out=180,in=270] ($(v0)+(-0.8,0)$);
        \draw[rounded corners=15,fill=red,opacity=0.5] (-0.2,-0.6)--(-0.4,1.65)--(1.2,0.9)--cycle;
        \draw[rounded corners=15,fill=orange,opacity=0.5] (-0.65,1.2)--(0.9,0.3)--(1.2,2)--cycle;
        \draw[rounded corners=15,fill=yellow,opacity=0.5] (0.3,0.4)--(0.7,2)--(2.2,1.15)--cycle;
        \draw[rounded corners=15,fill=green,opacity=0.5] (0.35,1.75)--(1.7,0.7)--(2.8,2.1)--cycle;
        \draw[rounded corners=15,fill=blue,opacity=0.5] (1.2,1)--(2.4,2.3)--(3.1,0.75)--cycle;
        \draw[rounded corners=15,fill=violet,opacity=0.5] (2,2.2)--(2.45,0.7)--(4.55,1.2)--cycle;
      \end{scope}

      \foreach \v in {0,1,2,...,7} {
          \fill (v\v) circle (0.1);
      };

      \node at (2.7,-0.7) {$N$};
      \node at (2.2,-2) {$N$-ear};
    \end{tikzpicture}
    \begin{tikzpicture}
      \node (v1) at (0,0.5) {};
      \node (v2) at (0.7,-0.4) {};
      \node (v3) at (1.4,0.2) {};

      \begin{scope}
        \draw[fill=gray!70, opacity=0.5] ($(v1)+(-0.8,0)$) 
        to[out=0,in=180] ($(v2) + (-0.3,1.8)$) 
        to[out=0,in=180] ($(v3) + (0.3,0.8)$)
        to[out=0,in=0] ($(v2) + (0.5,-0.8)$)
        to[out=180,in=180] ($(v1)+(-0.8,0)$);
        \draw[rounded corners=15,fill=red,opacity=0.5] (-0.5,0.8)--(0.65,-0.85)--(1.9,0.4)--cycle;
      \end{scope}

      \foreach \v in {1,2,3} {
          \fill (v\v) circle (0.1);
      };

      \node at (1.5,-0.7) {$N$};
      \node at (0.6,-1.6) {$N$-triangle};
    \end{tikzpicture}

    \begin{tikzpicture}
      \node (v1) at (0,0) {};
      \node (v2) at (0.4,-0.1) {};
      \node (v3) at (2.5,0) {};
      \node (v4) at (3,0.3) {};
      \node (v5) at (1.3,1.2) {};

      \begin{scope}
        \draw[fill=gray!70, opacity=0.5] ($(v1)+(-0.8,0)$) 
        to[out=90,in=90] ($(v2) + (-0.3,0.5)$) 
        to[out=270,in=90] ($(v4) + (0.5,0.3)$)
        to[out=270,in=90] ($(v3) + (0.7,-0.2)$)
        to[out=270,in=0] ($(v2) + (0.5,-0.8)$)
        to[out=180,in=270] ($(v1)+(-0.8,0)$);
        \draw[rounded corners=25,fill=red,opacity=0.5] (-0.85,-0.3)--(0.8,-0.6)--(1.9,2)--cycle;
        \draw[rounded corners=25,fill=blue,opacity=0.5] (2.2,-0.6)--(3.7,0.4)--(0.6,1.7)--cycle;
      \end{scope}

      \foreach \v in {1,2,...,5} {
          \fill (v\v) circle (0.1);
      };

      \node at (1.3,-0.5) {$N$};
      \node at (1.3,-1.35) {$N$-windmill};
    \end{tikzpicture}
    \begin{tikzpicture}
      \node (v1) at (0,0) {};
      \node (v2) at (2,-0.4) {};
      \node (v3) at (2.8,0.8) {};
      \node (v4) at (1.3,1.2) {};
      \node (v5) at (1.3,1.8) {};

      \begin{scope}
        \draw[fill=gray!70, opacity=0.5] ($(v1)+(-0.8,0)$) 
        to[out=90,in=180] ($(v2) + (-0.3,0.5)$) 
        to[out=0,in=180] ($(v3) + (0.3,0.8)$)
        to[out=0,in=0] ($(v2) + (0.5,-0.8)$)
        to[out=180,in=270] ($(v1)+(-0.8,0)$);
        \draw[rounded corners=22,fill=red,opacity=0.5] (-0.7,-0.7)--(1.7,0.9)--(1.45,2.6)--cycle;
        \draw[rounded corners=22,fill=blue,opacity=0.5] (2.45,-1.4)--(0.8,1.15)--(1.3,2.6)--cycle;
        \draw[rounded corners=22,fill=green,opacity=0.5] (3.5,0.45)--(0.75,0.9)--(1.2,2.5)--cycle;
      \end{scope}

      \foreach \v in {1,2,...,5} {
          \fill (v\v) circle (0.1);
      };

      \node at (1,-0.5) {$N$};
      \node at (1.5,-1.6) {$N$-tripod};
    \end{tikzpicture}
    \begin{tikzpicture}
      \node (v1) at (0,0) {};
      \node (v2) at (2.2,-0.8) {};
      \node (v3) at (3.3,0.6) {};
      \node (v4) at (1.9,1) {};
      \node (v5) at (1.4,2) {};
      \node (v6) at (1,0.9) {};

      \begin{scope}
        \draw[fill=gray!70, opacity=0.5] ($(v1)+(-0.8,0)$) 
        to[out=90,in=180] ($(v2) + (-0.3,0.5)$) 
        to[out=0,in=180] ($(v3) + (0.3,0.8)$)
        to[out=0,in=0] ($(v2) + (0.5,-0.8)$)
        to[out=180,in=270] ($(v1)+(-0.8,0)$);
        \draw[rounded corners=22,fill=red,opacity=0.5] (-0.7,-0.7)--(1.6,0.9)--(1.55,2.65)--cycle;
        \draw[rounded corners=22,fill=blue,opacity=0.5] (2.65,-1.8)--(0.4,1.15)--(2.25,1.44)--cycle;
        \draw[rounded corners=22,fill=green,opacity=0.5] (3.95,0.3)--(1.4,0.7)--(1.1,2.55)--cycle;
        \draw[rounded corners=17,fill=yellow,opacity=0.5] (2.45,0.75)--(1.3,2.6)--(0.55,0.55)--cycle;
      \end{scope}

      \node at (1,-0.8) {$N$};
      \node at (1.75,-1.9) {$N$-table};

      \foreach \v in {1,2,...,6} {
          \fill (v\v) circle (0.1);
      };
    \end{tikzpicture}    
  \caption{$5$ types of partial hypergraphs.}
  \label{fig:mtp}
\end{figure}
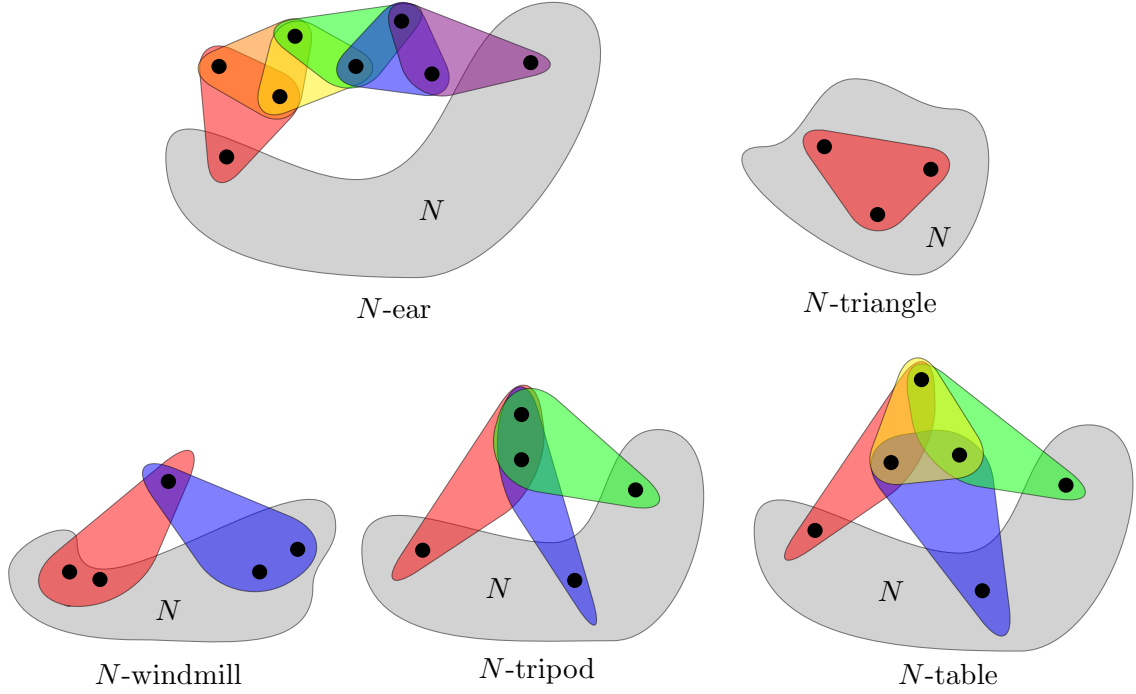
\begin{itemize}
  \item An \emph{$N$-ear} in $H$ is a tight path $P$ of length at least $2$ such that two ends are in $N$, all internal vertices are in $V-N$, and no edge in $E(H) - E(P)$ is incident with an internal vertex of $P$.
  \item An \emph{$N$-triangle} in $H$ is a partial hypergraph $(V',E')$ of $H$ without isolated vertices 
  such that %
  \begin{enumerate}[label=\rm($\triangle$\arabic*)]
    \item   $V'\subseteq N$, 
    \item   $\abs{E'}=1$,
    \item\label{item:delta3} $\abs{V'\cap e}\neq 2$ for all edges $e$ of $H$.
  \end{enumerate}

  \item An \emph{$N$-windmill} in~$H$ 
  is a partial hypergraph $(V',E')$ of $H$ without isolated vertices
  for which there is a vertex $v\notin N$ such that 
  \begin{enumerate}[label=\rm(W\arabic*)]
    \item\label{item:W1}   $V'-\{v\}\subseteq N$, 
    \item\label{item:W2}   $E'$ is the set of all edges of $H$ incident with $v$, 
    \item\label{item:W3}   $E'\neq \emptyset$,
    \item\label{item:W4}  $\abs{e\cap e'}\neq 2$ for all edges $e$ of $H$ and all edges $e'$ in $E'$.
  \end{enumerate}

  \item An \emph{$N$-tripod} in $H$ is a partial hypergraph $(V',E')$ of $H$ without isolated vertices for which there are two distinct vertices $v,w\notin N$
  such that 
  \begin{enumerate}[label=\rm(Y\arabic*)]
    \item\label{item:Y1} $V'-\{v,w\}\subseteq N$,
    \item\label{item:Y2} $E'$ is the set of all edges of $H$ incident with both $v$ and $w$,
    \item\label{item:Y3} $\abs{E'}=3$, 
    \item\label{item:Y4} no edge of $H$ is incident with exactly one of $v$ and $w$.
  \end{enumerate}

  \item An \emph{$N$-table} in $H$ is a partial hypergraph $(V',E')$ of $H$ without isolated vertices for which there are three distinct vertices $u,v,w\notin N$ 
  such that 
  \begin{enumerate}[label=\rm(T\arabic*)]
    \item\label{item:T1} $V'-\{u,v,w\}\subseteq N$,
    \item\label{item:T2} $E'$ is the set of all edges of $H$ incident with at least two of $u$, $v$, and $w$, 
    \item\label{item:T3} $\abs{E'}=4$,
    \item\label{item:T4} $e\cap \{u,v,w\}\neq e'\cap \{u,v,w\}$ and $e-\{u,v,w\}\neq e'-\{u,v,w\}$ for distinct edges $e$, $e'$ in $E'$,
    \item\label{item:T5} no edge of $H$ is incident with exactly one of $u$, $v$, and $w$.
  \end{enumerate}

\end{itemize}
See Figure~\ref{fig:mtp} for illustrations for these $5$ types of hypergraphs.
It is easy to observe the following two lemmas from the definition.

\begin{lemma}\label{lem:disjoint}
  Let $H$ be a $3$-uniform hypergraph and $N \subseteq V(H)$.
  If each of $H_1$ and $H_2$ is an $N$-ear, an $N$-triangle, an $N$-windmill, an $N$-tripod, or an $N$-table in $H$, then 
  $H_1 = H_2$ or 
  $(V(H_1)-N)\cap (V(H_2)-N)=\emptyset$.
  \qed
\end{lemma}
\begin{lemma}\label{lem:paths-in-gadgets}
  Let $H$ be a $3$-uniform hypergraph and $N\subseteq V(H)$.
  Let $H'$ be an $N$-ear,  
  an $N$-triangle, an $N$-windmill, an $N$-tripod, or an $N$-table of $H$.
  \begin{enumerate}[label=\rm(\roman*)]
    \item\label{item:pathsgadget1} If $P$ is a tight path of $H$ such that $E(P)\cap E(H')\neq\emptyset$, then $E(P)\subseteq E(H')$.
    \item\label{item:pathsgadget2} If $P$ is a maximal tight path of $H$ contained in $H'$, then 
      $V(P)-N=V(H')-N$.\qed
  \end{enumerate}
\end{lemma}

\begin{theorem}\label{thm:fan}
  Let $S$ be a $3$-connected isotropic system with at least $5$ vertices and let $N$ be the set of non-essential vertices in $S$.
  If $N \neq \emptyset$,
  then the set of edge sets of all $N$-ears, $N$-triangles, $N$-windmills, $N$-tripods, and $N$-tables in $H(S)$
  is a partition of 
  the edge set of $H(S)$.
\end{theorem}

\begin{proof}
  As $N \neq \emptyset$, by Lemma~\ref{lem:cyclic}, $S$ is not cyclic.
  Note that for distinct $x$ and $y$ in $N$, there is at most one edge of $H(S)$ containing both $x$ and $y$ by Lemma~\ref{lem:intess}.
  It suffices to show that each edge $e$ of $H(S)$ is contained in an $N$-ear, an $N$-triangle, an $N$-windmill, an $N$-tripod, or an $N$-table,
  because if $e$ is contained in two of such partial hypergraphs $H_1$ and $H_2$, then either both $H_1$ and $H_2$ are $N$-triangles, meaning that $H_1=H_2$, 
  or $e$ is incident with a vertex not in $N$, implying that $H_1$ and $H_2$ share a vertex not in $N$, thus $H_1=H_2$ by Lemma~\ref{lem:disjoint}.

  If all three vertices incident with $e$ are non-essential in $S$, then a  hypergraph $(e,\{e\})$ is an $N$-triangle in $H(S)$
  because \ref{item:delta3} holds by Lemma~\ref{lem:intess}.
  Therefore, we may assume that $e$ is incident with an essential vertex in $S$.
  Let $P$ be a maximal tight path in $H(S)$ containing $e$.
  We denote $P$ by a sequence $av_1v_2 \cdots v_kb$ of distinct vertices, where $k$ is the length of $P$.
  If $k\geq 2$, then the internal vertices $v_1,\dots,v_k$ of $P$ are essential in $S$ by Lemma~\ref{lem:intess}.  If $k=1$, then by relabelling, we may assume that $v_1$ is essential in $S$.
  By Proposition~\ref{prop:maxpath}, two ends $a$ and $b$ of $P$ are non-essential in $S$.

  Suppose that $P$ includes every edge of $H(S)$ incident with some of $v_1,v_2,\dots,v_k$.
  If $k=1$, then $P$ is an $N$-windmill because 
  \ref{item:W4} holds by the assumption that $P$ is maximal.
  If $k>1$, then $P$ is an $N$-ear.
  Therefore, we may assume that $H(S)$ has an edge $f \not \in E(P)$ incident with some of $v_1,\dots,v_k$.
  If the length of $P$ is more than $3$, then no such $f$ exists by %
  Proposition~\ref{prop:v12}.
  Hence $k\leq 3$.

  \textbf{Case I.} $k=3$.

  By Lemma~\ref{lem:v2}, $f$ is not incident with $v_2$.
  Without loss of generality, we may assume that $f$ is incident with $v_1$.
  Since $v_1 \in f$ and $v_2 \not\in f$, by Lemma~\ref{lem:codeg2} applied to edges $\{a,v_1,v_2\}$ and $\{v_1,v_2,v_3\}$,
  we deduce that 
  $f$ is incident with $a$ or $v_3$.
  By Lemma~\ref{lem:intess}, $f$ is not incident with $a$ because otherwise %
  $f$ and $\{a,v_1,v_2\}$
  form a tight path, implying that $a$ is essential in~$S$. 
  Therefore, $f = \{v_1,v_3,c\}$ for some $c \in V(H(S)) - V(P)$.

  Now we show that $H':=(V(P)\cup \{c\}, E(P) \cup \{f\})$ is an $N$-table in $H(S)$.
  Note that $V(H') = \{a,b,c,v_1,v_2,v_3\}$ and $E(H') = \{\{a,v_1,v_2\}, \{b,v_2,v_3\}, \{c,v_1,v_3\}, \{v_1,v_2,v_3\}\}$.
  Thus \ref{item:T3} holds and \ref{item:T4} holds for $\{u,v,w\} = \{v_1,v_2,v_3\}$.
  By Lemma~\ref{lem:v2} applied to tight paths $av_1v_2v_3b$, $av_2v_1v_3c$, and $bv_2v_3v_1c$, we deduce that $H'$ has all edges of $H$ incident with $v_1$, $v_2$, or $v_3$.
  This implies not only  \ref{item:T2} and \ref{item:T5}, 
  but also $c\in N$  
  by Proposition~\ref{prop:maxpath}
  because 
  $av_2v_1v_3c$ is a maximal tight path. It follows that \ref{item:T1} holds.

  \textbf{Case II.} $k=2$.

  We prove that $H' := (V(P) \cup f, E(P) \cup \{f\})$ is an $N$-tripod in $H(S)$.
  Trivially \ref{item:Y3} holds.
  To see \ref{item:Y4} with $v:=v_1$ and $w:=v_2$,
  suppose that there is an edge $g$ of $H(S)$ incident with exactly one of $v_1$ and $v_2$, say $v_1$ by symmetry.
  By Lemma~\ref{lem:codeg2} applied to $\{v_1,v_2,a\}$ and $\{v_1,v_2,b\}$, we deduce that $g$ is incident with $a$ or $b$.
  However, 
  by Lemma~\ref{lem:intess}, 
  $g\cap \{a,v_1,v_2\}\neq\{a,v_1\}$
  and $g\cap \{v_1,v_2,b\}\neq \{v_1,b\}$,
  contradicting our previous conclusion.
  This proves \ref{item:Y4}.
  This also implies that $f = \{v_1,v_2,c\}$ for some $c \in V(H(S)) - V(P)$.

  By Lemma~\ref{lem:codeg3}, $H'$ satisfies \ref{item:Y2}.
  By \ref{item:Y4}, 
  $av_1v_2c$ is a maximal tight path in $H$.
  By Proposition~\ref{prop:maxpath}, $c \in N$ and so $V(H') - \{v_1,v_2\} = \{a,b,c\} \subseteq N$, implying \ref{item:Y1}.

  \textbf{Case III.} $k=1$.
  
  Let $E'$ be the set of edges incident with $v_1$ and let $V'$ be the set of all vertices incident with an edge in $E'$.
  We show that $H' := (V',E')$ is an $N$-windmill.
  By definition, \ref{item:W2} holds and $e=\{a,v_1,b\}\in E'$, implying \ref{item:W3}.

  To see~\ref{item:W1} with $v:=v_1$, 
  suppose that 
  there is an edge $g$ 
  incident with both $v_1$ and a vertex not in $N\cup\{v_1\}$. 
  Let $Q$ be a maximal tight path containing~$g$.
  By Proposition~\ref{prop:maxpath}, the length of $Q$ is 
  at least $2$
  and $v_1$ is not an end of $Q$.
  Then $Q$ has two edges $e_1$ and $e_2$ such that both are incident with $v_1$ and $|e_1\cap e_2| = 2$.
  By Lemma~\ref{lem:codeg2} applied to $e_1$ and $e_2$, 
  $|\{a,v_1,b\} \cap e_i| \geq 2$ for some $i \in \{1,2\}$.
  We may assume that $|\{a,v_1,b\} \cap e_1| \geq 2$.
  Since $P$ is a maximal tight path, 
  we deduce that $|\{a,v_1,b\} \cap e_1| \neq 2$
  and therefore $\{a,v_1,b\} = e_1$.
  Then $P$ is a proper subpath of $Q$, contradicting the assumption that $P$ is a maximal tight path. Thus, \ref{item:W1} holds.

  If there are edges $g \in E(H(S))$ and $g' \in E'$ such that $|g \cap g'| = 2$, then by Lemma~\ref{lem:intess}, $g'$ is incident with at least two essential vertices in $S$, contradicting \ref{item:W1}. 
  Therefore, \ref{item:W4} holds.
\end{proof}

\begin{corollary}\label{cor:partition}
  Let $S$ be a $3$-connected isotropic system with at least $5$ vertices and let $N$ be the set of non-essential vertices in $S$.
  Let $P_1$ and $P_2$ be maximal tight paths of $H(S)$.
  If $N \neq \emptyset$,
  then $V(P_1) - N$ and $V(P_2) - N$ are equal or disjoint.
\end{corollary}
\begin{proof}
  Suppose that there is a vertex $v\in (V(P_1)-N)\cap (V(P_2)-N)$.
  For $i\in \{1,2\}$, let $e_i$ be an edge of $P_i$ incident with $v$.
  By Theorem~\ref{thm:fan}, 
  for each $i\in\{1,2\}$, there is a partial hypergraph $H_i$ of $H$ 
  such that $H_i$ includes $e_i$
  and $H_i$ is an $N$-ear, an $N$-windmill, an $N$-tripod, 
  or an $N$-table. 
  We remark that since $v\notin N$, $H_i$ is not an $N$-triangle.
  By Lemma~\ref{lem:paths-in-gadgets}\ref{item:pathsgadget1}, $E(P_i) \subseteq E(H_i)$.
  By Lemma~\ref{lem:disjoint}, 
  $H_1=H_2$.
  Since $P_1$ and $P_2$ are maximal tight paths of $H$ 
  contained in $H_1=H_2$, 
  by Lemma~\ref{lem:paths-in-gadgets}\ref{item:pathsgadget2}, 
  $V(P_1)-N=V(P_2)-N$.
\end{proof}

\section{Prime graphs with at most two non-essential vertices}\label{sec:main}
We prove a part of Theorem~\ref{thm:noness3}, that is, if a prime graph with at least $5$ vertices has at most $2$ non-essential vertices, then it is locally equivalent to $\theta(\ell_1,\dots,\ell_m)$ for some~$m$ and~$\ell_i$.
We first prove some lemmas in order to obtain a specific Eulerian vector of an isotropic system.

For two vectors $\mb{a}\in K^A$ and $\mb{b}\in K^B$ with disjoint sets $A$ and $B$, let $\mb{a} \oplus \mb{b}$ be a vector in $K^{A\cup B}$ such that
\begin{equation*}
  (\mb{a}\oplus\mb{b})(v) =
  \begin{cases}
    \mb{a}(v) & \text{if $v\in A$}, \\
    \mb{b}(v) & \text{otherwise}.
  \end{cases}
\end{equation*}
If $|B| = 1$ and $\mb{b}(w) = x$ for $w\in B$,
then we simply write $\mb{a} \oplus \mb{b}$ as $\mb{a} \oplus x$.
Let $0_A$ be the zero vector in~$K^A$.

\begin{lemma}\label{lem:eur1}
  Let $S = (V,L)$ be an isotropic system, $v$ be a vertex in $V$, and $x \in K - \{0\}$.
  If $\mb{a} \in K^{V-\{v\}}$ is an Eulerian vector of $S|^v_x$ and $0_{V-\{v\}} \oplus x \not \in L$, then $\mb{a} \oplus x$ is an Eulerian vector of $S$
\end{lemma}
\begin{proof}
  Suppose that $\mb{a} \oplus x$ is not an Eulerian vector of $S$.
  We have a nonempty subset $X$ of $V$ such that $(\mb{a}\oplus x)[X] \in L$.
  Then $\mb{a}[X -\{v\}] \in L|^v_x$.
  Since $\mb{a}$ is an Eulerian vector of $S|^v_x$, we deduce $X-\{v\} = \emptyset$.
  Hence $X = \{v\}$ and therefore $0_{V-\{v\}} \oplus x = (\mb{a} \oplus x)[X] \in L$, which is a contradiction.
\end{proof}

\begin{lemma}\label{lem:eur2}
  Let $W = \{w_1, w_2, \dots, w_k\}$ be a set of vertices in an isotropic system $S$, and let $x_i \in K-\{0\}$ for $1\leq i \leq k$.
  Let $\mb{a}_1, \mb{a}_2, \dots, \mb{a}_k$ be vectors in $L$ such that $\innk{\mb{a}_i(w_i)}{x_i} = 1$ for every $1\leq i\leq k$, and $\innk{\mb{a}_i(w_j)}{x_j} = 0$ for all $1\leq j < i \leq k$.
  Then $S$ has an Eulerian vector $\mb{c}$ such that $\mb{c}(w_i) = x_i$ for every $1\leq i \leq k$.
\end{lemma}
\begin{proof}
  We proceed by induction on $k\geq 0$.
  For $k=0$, by Lemma~\ref{lem:eur}, $S$ has an Eulerian vector.
  Now we assume that $k\geq 1$.
  Observe that $p_{V-\{w_1\}}(\mb{a}_i) \in L|^{w_1}_{x_1}$ for all $2 \leq i \leq k$, and $\innk{p_{V-\{w_1\}}(\mb{a}_i)(w_i)}{x_i} = \innk{\mb{a}_i(w_i)}{x_i} = 1$ for all $2 \leq i \leq k$, and $\innk{p_{V-\{w_1\}}(\mb{a}_i)(w_j)}{x_j} = \innk{\mb{a}_i(w_j)}{x_j} = 0$ for all $2 \leq j < i \leq k$.
  By the induction hypothesis, $S|^{w_1}_{x_1}$ has an Eulerian vector $\mb{c}'$ such that $\mb{c}' (w_i) = x_i$ for every $2 \leq i \leq k$.
  A vector $0_{V-\{w_1\}} \oplus x_1$ is not in $L$, since $\inn{\mb{a}_1}{0_{V-\{w_1\}} \oplus x_1} = \innk{\mb{a}_1(w_1)}{x_1} = 1$.
  Thus, by Lemma~\ref{lem:eur1}, the proof is completed.
\end{proof}

\begin{lemma}\label{lem:vec}
  Let $S = (V,L)$ be a $3$-connected isotropic system with at least $5$ vertices.
  Then for distinct $u,v \in V$ and nonzero $x,y \in K$, there is a vector $\mb{a}$ in $L$ such that $\innk{\mb{a}(u)}{x} = 1$ and $\innk{\mb{a}(v)}{y} = 0$.
\end{lemma}
\begin{proof}
  Since $S$ is $3$-connected, 
  by Lemma~\ref{lem:pt}, $0_{V-\{v\}} \oplus x \notin L=L^\perp$ and therefore $L$ has a vector $\mb{a}_1$ such that $\innk{\mb{a}_1(u)}{x} = 1$.

  Let $c\in K^V$ be a vector such that $c(u)=x$, $c(v)=y$, and $c(w)=0$ for all $w\in V-\{u,v\}$. 
  Again by Lemma~\ref{lem:pt},
  $c\notin L=L^\perp$, and therefore $L$ has a vector $\mb{a}_2$ such that $\innk{\mb{a}_2(u)}{x} \neq \innk{\mb{a}_2(v)}{y}$.

  We may assume that $\innk{\mb{a}_1(v)}{y} = 1$, since otherwise we finish the proof by taking $\mb{a} = \mb{a}_1$.
  We may assume that $\innk{\mb{a}_2(u)}{x} = 0$ and $\innk{\mb{a}_2(v)}{y} = 1$, since otherwise we finish the proof by taking $\mb{a} = \mb{a}_2$.
  Then $\mb{a} = \mb{a}_1 + \mb{a}_2$ satisfies the desired condition.
\end{proof}

\begin{proposition}\label{prop:gn2}
  Let $G$ be a prime graph with at least $5$ vertices.
  If $G$ has at most $2$ non-essential vertices, then $G$ is locally equivalent to a graph isomorphic to $\theta(\ell_1,\dots,\ell_m)$ for some $m$ and $\ell_i$.
\end{proposition}

\begin{proof}
  Since a cycle graph of length $k$ is isomorphic to $\theta(1,k-1)$, we may assume that $G$ is not locally equivalent to a cycle graph.
  Then by Theorem~\ref{thm:noness2}, $G$ has exactly two non-essential vertices $u$ and $v$.

  Let $V = V(G)$ and let $S$ be an isotropic system having $G$ as a fundamental graph.
  Then $S$ is not cyclic, since all fundamental graphs of $S$ are locally equivalent.
  Furthermore, $S$ is $3$-connected and has exactly two non-essential vertices by Corollaries~\ref{cor:igconn2} and~\ref{cor:ess1}.

  By Lemma~\ref{lem:tri1}, every essential vertex in $S$ is in a tight path of $H(S)$.
  By Proposition~\ref{prop:maxpath}, every maximal tight path has $u$ and $v$ as its ends.
  By Corollary~\ref{cor:partition}, there are maximal tight paths $P_1,P_2,\dots,P_m$ of $H(S)$ such that $V(P_1)-\{u,v\}, V(P_2)-\{u,v\}, \dots, V(P_m)-\{u,v\}$ partition the set of essential vertices in $S$.
  Let us denote $P_i$ as a sequence $u v_{i,1}\dots v_{i,\ell(i)} v$.
  We may assume that $\ell(1) \leq \ell(2) \leq \dots \leq \ell(m)$.
  Let $\mbt_{i,j}$ be a triangle in~$S$ whose support is $\{v_{i,j-1},v_{i,j},v_{i,j+1}\}$ for each $1\leq i\leq m$ and $1\leq j\leq \ell(i)$, where $v_{i,0} := u$ and $v_{i,\ell(i)+1} := v$ for each $i$.

  Since all fundamental graphs of $S$ are locally equivalent, it is enough to show that $S$ has a fundamental graph isomorphic to $\theta(\ell(1)+1,\ell(2)+1,\dots,\ell(m)+1)$ or $\theta(1,\ell(1)+1,\ell(2)+1,\dots,\ell(m)+1)$.

  Since $u$, $v$, and $v_{i,j}$ for all $1\leq i\leq m$ and $1\leq j\leq \ell(i)$ are distinct, by symmetry of the nonzero elements in $K$, we can assume that $\mbt_{1,1}(u) = \mbt_{1,\ell(1)}(v) = \alpha$ and $\mbt_{i,j}(v_{i,j}) = \beta$ for all $1\leq i\leq m$ and $1\leq j\leq \ell(i)$.
  If $\ell(1) = \ell(2) = 1$, then $v_{1,1} u v v_{2,1}$ is a tight path properly containing $P_1$ and $P_2$, which contradicts that $P_1$ and $P_2$ are maximal tight paths.
  Thus, $2\leq \ell(2) \leq \dots \leq \ell(m)$.
  Then by Lemma~\ref{lem:tri0}, $\mbt_{i,1}(u) = \mbt_{1,1}(u) = \alpha$ and $\mbt_{i,\ell(i)}(v) = \mbt_{1,\ell(1)}(v) = \alpha$ for all $2\leq i\leq m$.
  Applying Lemma~\ref{lem:tri0}, for all $1\leq i\leq m$ and $2\leq j\leq \ell(i)-1$, we have $\mbt_{i,j-1}(v_{i,j}) = \mbt_{i,j+1}(v_{i,j}) \neq \mbt_{i,j}(v_{i,j}) = \beta$.
  By symmetry in $K-\{0\}$, we can assume that $\mbt_{i,j-1}(v_{i,j}) = \mbt_{i,j+1}(v_{i,j}) = \alpha$ for all $1\leq i\leq m$ and $2\leq j\leq \ell(i)-1$.
  Similarly, for all $1\leq i\leq m$ with $\ell(i) \geq 2$, we have $\mbt_{i,2}(v_{i,1}) \neq \mbt_{i,1}(v_{i,1}) = \beta$ and $\mbt_{i,\ell(i)-1}(v_{i,\ell(i)}) \neq \mbt_{i,\ell(i)}(v_{i,\ell(i)}) = \beta$, and thus we can assume that $\mbt_{i,2}(v_{i,1}) = \mbt_{i,\ell(i)-1}(v_{i,\ell(i)}) = \alpha$.
  In short, we assumed that
  \begin{equation*}
    \mbt_{i,j}(v_{i,j-1}) = \mbt_{i,j}(v_{i,j+1}) = \alpha
    \text{ and }
    \mbt_{i,j}(v_{i,j}) = \beta
  \end{equation*}
  for all $1\leq i\leq m$ and $1\leq j\leq \ell(i)$.

  By Lemma~\ref{lem:vec}, there exist vectors $\mb{a}$ and $\mb{b}$ in $L$ such that $\innk{\mb{a}(u)}{\alpha} = 1$, $\innk{\mb{a}(v)}{\alpha} = 0$, $\innk{\mb{b}(u)}{ \alpha} = 0$, and $\innk{\mb{b}(v)}{\alpha} = 1$.
  Let us denote $T := \{ \mbt_{i,j} : 1\leq i\leq m \text{ and } 1\leq j\leq \ell(i) \}$.
  Since $\innk{\mbt_{i,j}(u)}{\alpha} = \innk{\mbt_{i,j}(v)}{\alpha} = 0$ for all $1\leq i\leq m$ and $1\leq j\leq \ell(m)$, $\{\mb{a},\mb{b}\} \cup T$ is linearly independent.
  Then $\{\mb{a},\mb{b}\} \cup T$ is a basis of $L$ because $|T| = \sum_{i=1}^m \ell(i) = |V|-2$,

  Let $w_1, w_2, \dots, w_n$ be all vertices of $V$ such that $w_1 = u$ and $w_2 = v$.
  Let $\mb{a}_1 = \mb{a}$ and $\mb{a}_2 = \mb{b}$.
  For $k\geq 3$, let $\mb{a}_k = \mbt_{i,j}$ if $w_k = v_{i,j}$.
  Then $\innk{\mb{a}_k(w_k)}{\alpha} = 1$ for all $k$, and $\innk{\mb{a}_k(w_s)}{\alpha} = 0$ for $1 \leq s < k \leq n$.
  By applying Lemma~\ref{lem:eur2} for $w_1,\dots,w_n$ and $\mb{a}_1, \dots, \mb{a}_n$, we obtain an Eulerian vector~$\mb{c}$ of~$S$ such that $\mb{c}(w_k) = \alpha$ for all $1 \leq k \leq n$.
  Let $G'$ be the fundamental graph of $S$ with respect to~$\mb{c}$.
  Then $v_{i,j}$ is only adjacent to $v_{i,j-1}$ and $v_{i,j+1}$ in $G'$ for each $1\leq i\leq m$ and $1\leq j\leq \ell(i)$ because $\mbt_{i,j}$ is a vector in the fundamental basis of $S$ with respect to~$\mb{c}$.
  Therefore, $G'$ is isomorphic to $\theta(\ell(1)+1,\dots,\ell(m)+1)$ or $\theta(1,\ell(1)+1,\dots,\ell(m)+1)$ depending on the adjacency between $u$ and~$v$.
\end{proof}

\section{Graphs consisting of internally-disjoint paths}\label{sec:theta}

To complete the proof of Theorem~\ref{thm:noness3}, we investigate the condition that $\theta(\ell_1,\dots,\ell_m)$ is prime and has at most $2$ non-essential vertices.

The following lemma provides three ways to extend a prime graph.
For a graph $G$ and its induced subgraph $H$, a sequence $v_0,v_1,\dots,v_{\ell}$  of distinct vertices of $G$ is a \emph{handle} of $H$ if $\ell \geq 3$, $\{v_0,\dots,v_{\ell}\} \cap V(H) = \{v_0,v_{\ell}\}$, and
$v_i$ is only adjacent to $v_{i-1}$ and $v_{i+1}$ in $G[V(H) \cup \{v_1,\dots,v_{\ell-1}\}]$ for every $1\leq i \leq \ell-1$.
We say that $G[V(H) \cup \{v_1,\dots,v_{\ell-1}\}]$ is obtained from $H$ by \emph{adding a handle} of length~$\ell$.

\begin{lemma}[Geelen~\cite{Geelen1996thesis}]
  \label{lem:pext}
  Let $G$ be a graph with at least $5$ vertices.
  \begin{enumerate}[label=\rm(\alph*)]
    \item\label{item:p1} 
    If $G$ has a vertex $v$ of degree at least $2$ such that $G \setminus v$ is prime and $v$ has no twin in $G$, then $G$ is prime.
    \item\label{item:p2}
    If $G$ is obtained from its prime induced subgraph with at least $4$ vertices by adding a handle, then $G$ is prime.
    \item\label{item:p3} If $G$ has an edge $e$ such that both ends of $e$ have degree $2$ and $G/e$ is prime, then $G$ is prime.
  \end{enumerate}
\end{lemma}
\begin{proof}
  Both~\ref{item:p1} and~\ref{item:p2} were proved by Geelen in Lemma~5.3 and Proposition~5.5, respectively, of~\cite{Geelen1996thesis}.
  For~\ref{item:p3}, $G/e$ is isomorphic to $G*v \setminus v$, where $e=vw$.
  Since $G*v\setminus v$ is prime with at least $4$ vertices, it is easy to check that none of the two neighbors of $v$ is a twin of $v$ in $G*v$, and the neighbor of $w$ other than $v$ is not a twin of $v$ in $G*v$.
  Hence $v$ has no twin in $G*v$.
  By~\ref{item:p1}, $G*v$ is prime and therefore $G$ is prime.
\end{proof}

\begin{proposition}\label{prop:theta1}
  Let $m\geq 2$ and $\ell_1,\dots,\ell_m$ be positive integers, and $G = \theta(\ell_1,\dots,\ell_m)$ be a graph with at least $5$ vertices.
  Then $G$ is prime if and only if $|\{i:\ell_i = 2\}| \leq 1$.
\end{proposition}

\begin{proof}
  If $|\{i: \ell_i = 2\}| \geq 2$, then $G$ has twins and thus it is not prime.
  Now, let us prove the backward direction.
  By Lemma~\ref{lem:pext}\ref{item:p3}, it suffices to show that $\theta(1,2,3,\dots,3)$, $\theta(2,3,\dots,3)$, $\theta(1,3,\dots,3)$, and $\theta(3,\dots,3)$ are prime.
  By Lemma~\ref{lem:pext}\ref{item:p2}, it is enough to show that $\theta(1,2,3)$, $\theta(2,3)$, $\theta(1,3,3)$, and $\theta(3,3)$ are prime.
  Since $\theta(2,3)$ and $\theta(3,3)$ are cycles of length $5$ and $6$, respectively, they are prime.
  For the unique common neighbor $v$ of two degree-$3$ vertices in $\theta(1,2,3)$, 
  the graph $\theta(1,2,3) * v$ is isomorphic to $C_5$ and therefore $\theta(1,2,3)$ is prime.
  For an edge $e$ in $\theta(1,3,3)$ whose both ends have degree $2$, 
  the graph $\theta(1,3,3) \wedge e$ is isomorphic to $C_6$ and so $\theta(1,3,3)$ is prime.
\end{proof}

The following lemma is useful for finding pivotal vertices in a graph.
Remember that all pivotal vertices are essential.

\begin{lemma}\label{lem:essv}
  If $v$ and $w$ are adjacent vertices of degree $2$ in a graph $G$ with at least $5$ vertices, then both $v$ and $w$ are pivotal in $G$.
\end{lemma}
\begin{proof}
  By Lemma~\ref{lem:prime}, 
  neither $G \setminus v$ nor $G \wedge vw \setminus v$ is prime and therefore 
  $v$ is pivotal in $G$. Similarly, $w$ is pivotal in $G$.
\end{proof}

To find non-essential vertices of $\theta(\ell_1,\dots,\ell_m)$, we will use the next two lemmas.
A graph is \emph{outerplanar} if it has a planar embedding such that every vertex lies on the boundary of the outer face.

\begin{lemma}\label{lem:2outer}
  An outerplanar graph with at least $5$ vertices is prime if and only if it is $2$-connected.
\end{lemma}
\begin{proof}
  By Lemma~\ref{lem:prime}, it is enough to prove the backward direction. 
  Let $G$ be a $2$-connected outerplanar graph with at least $5$ vertices.
  We fix an embedding of~$G$ into the plane such that the boundary of the outer face contains every vertex of $G$.
  Since $G$ is $2$-connected, there is a cycle $C$ in $G$ corresponding to the boundary of the outer face.
  Let $v_1,v_2,\dots,v_n$ be the vertices of $C$ in the clockwise order, starting at a vertex $v_1$.
  Suppose that $G$ has a split $(X,Y)$.
  We may assume that $\abs{X}\ge 3$ by swapping $X$ and $Y$ if necessary.
  By rotational symmetry, we may assume that $v_{n-1} \in X$ and $v_n \in Y$.

  We claim that if $v_i \in Y$ for some $i\in \{1,2,\dots,n-3\}$, then $v_{i+1} \in Y$.
  Suppose that $v_{i+1}\in X$.
  Since $(X,Y)$ is a split, $v_i$ is adjacent to $v_{n-1}$, and $v_{i+1}$ is adjacent to $v_{n}$, contradicting the assumption that $G$ is outerplanar.
  This proves the claim.

  By the claim, $Y-\{v_{n-1},v_n\}=Y-\{v_n\}=\{v_i:j\le i\le n-2\}$ for some $j\in \{1,2,\ldots,n-1\}$.
  Since $\abs{X}\ge 3$ and $\abs{Y}\ge 2$, we deduce that $3\le j\le n-2$.
  Since $(X,Y)$ is a split, $v_1,v_{j-1}\in X$, $v_n,v_{j}\in Y$, and $v_1v_n,v_{j-1}v_j\in E(G)$, we deduce that 
  $v_1v_j, v_{j-1}v_n\in E(G)$, contradicting the assumption that $G$ is outerplanar.
\end{proof}

\begin{lemma}\label{lem:primepart}
  Let $G$ be a graph and $V_1,V_2,V_3$ be disjoint subsets of $V(G)$ such that $|V_i|\geq 2$ for all~$i$.
  If $G-V_i$ is prime for each~$i$, then $G$ is prime.
\end{lemma}
\begin{proof}
  Suppose that $G$ is not prime.
  Then $G$ has a split $(X,Y)$.
  We call vertices in $X$ \emph{red} and vertices in $Y$ \emph{blue}.

  If neither $V_i$ nor $V_j$ is monochromatic for some distinct $i,j$, then $G-V_k$ has a split for $k\neq i,j$, contradicting that $G-V_k$ is prime.
  Thus, there is at most one non-monochromatic $V_i$.
  We may assume that $V_1$ and $V_{2}$ are monochromatic.
  If $V_1$ and $V_2$ have different colors, then $G-V_3$ has a split, which is a contradiction.
  Therefore, we may assume that $V_1$ and $V_2$ are red.
  Then all blue vertices belong to $V(G) - (V_1 \cup V_2)$ and therefore $G-V_1$ has a split, which is a contradiction.
\end{proof}

Now we prove two lemmas presenting non-essential vertices of $\theta(\ell_1,\dots,\ell_m)$.

\begin{lemma}\label{lem:tp2}
  Let $G = \theta(\ell_1,\dots,\ell_m)$ such that $m \geq 3$, $\ell_1 \leq \dots \leq \ell_m$, $\ell_1 \neq 2$, $\ell_2\geq 3$, and if $m=3$, then $\ell_1 \geq 3$ or $\ell_2 \geq 4$.
  For a vertex $x$ of $G$, the following are equivalent: (i) $x$ has degree larger than $2$, (ii) $x$ is non-pivotal, and (iii) $x$ is non-essential.
  In particular, $G$ has exactly two non-essential vertices.
\end{lemma}
\begin{proof}
  Trivially, (iii) implies (ii).
  Observe that $|V(G)| = 2+\sum_{i=1}^m (\ell_i-1) \geq 6$.
  By Lemma~\ref{lem:essv}, (ii) implies (i).
  Therefore, it suffices to show that (i) implies (iii).
  Let $u$ and $v$ be the two distinct vertices of degree $m$ in $G$, which are the only vertices of degree larger than $2$.
  We prove that $u$ and $v$ are non-essential in $G$.
  By symmetry, it is enough to show that $u$ is non-essential.

  We claim that both $G*u\setminus u$ and $G\wedge uw \setminus u$ are prime.
  This claim implies that $u$ is non-essential in~$G$.

  We proceed by induction on $m \geq 3$.
  Suppose that $\ell_1 = 1$ and $\ell_2 = 3$.
  Then $m\geq 4$.
  For the middle edge $e$ of a path of length $3$ between $u$ and $v$ in $G$, 
  a graph $H:=G \wedge e$ is isomorphic to $\theta(3,\ell_3,\dots,\ell_m)$.
  By the inductive hypothesis, both $H*u\setminus u$ and $H \wedge uw \setminus u$ are prime.
  Hence both
  $G*u\setminus u = (H \wedge e) * u \setminus u = (H * u \setminus u) \wedge e$
  and $G \wedge uw \setminus u = (H \wedge e) \wedge uw \setminus u = (H \wedge uw \setminus u) \wedge e$ are prime. 
  Therefore
  we can assume that $\ell_1 \geq 3$ or $\ell_2 \geq 4$.

  Let $P_1,\dots,P_m$ be internally-disjoint paths between $u$ and $v$ of lengths $\ell_1,\dots,\ell_m$, respectively, in~$G$ and let $V_i := V(P_i) - \{u,v\}$ for each $i$.
  Then $|V_i| = \ell_i - 1 \geq 2$ for all $i \geq 2$.
  Let~$w$ be the neighbor of $u$ in $P_1$.

  We first consider the case that $m=3$
  and $\ell_1 \geq 3$; see Figure~\ref{fig:tp2i}.
  Then $|V_1| = \ell_1 -1 \geq 2$.
  For each $i \in \{1,2,3\}$,
  $(G*u \setminus u) - V_i$ is a cycle of length $|V_j| + |V_k| + 1 \geq 5$, where $\{i,j,k\} = \{1,2,3\}$, and thus it is prime.
  Therefore, $G*u \setminus u$ is prime by Lemma~\ref{lem:primepart}.
  Note that $G\wedge uw \setminus u \setminus w$ is isomorphic to $\theta(\ell_1 -2, \ell_2, \ell_3)$ and thus it is prime by Proposition~\ref{prop:theta1}.
  Since $w$ has degree $2$ and has no twin in $G\wedge uw \setminus u$, by Lemma~\ref{lem:pext}\ref{item:p1}, $G\wedge uw \setminus u$ is prime.

  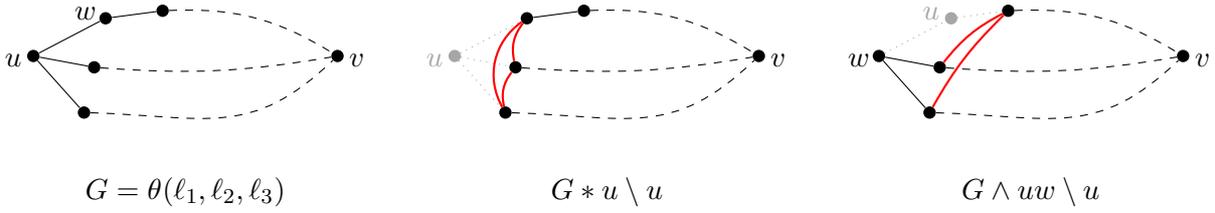
\begin{figure}
    \centering
    \begin{tikzpicture}
      \node[shape=circle,fill=black, scale=0.45] (u) at (0,0) {};
      \node[shape=circle,fill=black, scale=0.45] (v) at (4,0) {};

      \node[shape=circle,fill=black, scale=0.45] (x11) at (0.95,0.5) {};
      \node[shape=circle,fill=black, scale=0.45] (x12) at (1.7,0.6) {};

      \node[shape=circle,fill=black, scale=0.45] (x21) at (0.8,-0.15) {};

      \node[shape=circle,fill=black, scale=0.45] (x31) at (2/3,-0.75) {};

      \draw (u)+(-0.265,-0.05) node {$u$};
      \draw (v)+(0.255,-0.05) node {$v$};
      \draw (x11)+(-0.265,0.07) node {$w$};

      \draw (u) -- (x11);
      \draw (u) -- (x21);
      \draw (u) -- (x31);

      \draw (x11) -- (x12);
      \draw[dashed] (x12) to [bend left = 17] (v);
      \draw[dashed] (x21) to [bend right = 7] (v);
      \draw[dashed] (x31) .. controls (2.3,-0.8) and (2.85,-1.1) .. (v);
  
      \draw (2,-1.8) node {$G = \theta(\ell_1,\ell_2,\ell_3)$};
    \end{tikzpicture}
    \hspace{0.3cm}
    \begin{tikzpicture}
      \node[shape=circle,fill=black, scale=0.45, opacity=0.35] (u) at (0,0) {};
      \node[shape=circle,fill=black, scale=0.45] (v) at (4,0) {};

      \node[shape=circle,fill=black, scale=0.45] (x11) at (0.95,0.5) {};
      \node[shape=circle,fill=black, scale=0.45] (x12) at (1.7,0.6) {};

      \node[shape=circle,fill=black, scale=0.45] (x21) at (0.8,-0.15) {};

      \node[shape=circle,fill=black, scale=0.45] (x31) at (2/3,-0.75) {};

      \draw[opacity=0.35] (u)+(-0.265,-0.05) node {$u$};
      \draw (v)+(0.255,-0.05) node {$v$};

      \draw[opacity=0.35,dotted] (u) -- (x11);
      \draw[opacity=0.35,dotted] (u) -- (x21);
      \draw[opacity=0.35,dotted] (u) -- (x31);

      \draw (x11) -- (x12);
      \draw[dashed] (x12) to [bend left = 17] (v);
      \draw[dashed] (x21) to [bend right = 7] (v);
      \draw[dashed] (x31) .. controls (2.3,-0.8) and (2.85,-1.1) .. (v);

      \draw[color=red,thick] (x11) to [bend right=25] (x21);
      \draw[color=red,thick] (x21) to [bend right=25] (x31);
      \draw[color=red,thick] (x11) to [bend right=42] (x31);
  
      \draw (2,-1.8) node {$G*u\setminus u$};
    \end{tikzpicture}
    \hspace{0.3cm}
    \begin{tikzpicture}
      \node[shape=circle,fill=black, scale=0.45] (u) at (0,0) {};
      \node[shape=circle,fill=black, scale=0.45] (v) at (4,0) {};

      \node[shape=circle,fill=black, scale=0.45, opacity=0.35] (x11) at (0.95,0.5) {};
      \node[shape=circle,fill=black, scale=0.45] (x12) at (1.7,0.6) {};

      \node[shape=circle,fill=black, scale=0.45] (x21) at (0.8,-0.15) {};

      \node[shape=circle,fill=black, scale=0.45] (x31) at (2/3,-0.75) {};

      \draw (u)+(-0.265,-0.05) node {$w$};
      \draw (v)+(0.255,-0.05) node {$v$};
      \draw[opacity=0.35] (x11)+(-0.265,0.07) node {$u$};

      \draw[opacity=0.35,dotted] (u) -- (x11);
      \draw (u) -- (x21);
      \draw (u) -- (x31);

      \draw[opacity=0.35,dotted] (x11) -- (x12);
      \draw[dashed] (x12) to [bend left = 17] (v);
      \draw[dashed] (x21) to [bend right = 7] (v);
      \draw[dashed] (x31) .. controls (2.3,-0.8) and (2.85,-1.1) .. (v);

      \draw[color=red,thick] (x12) to [bend right=12] (x21);
      \draw[color=red,thick] (x12) to [bend right=8] (x31);
  
      \draw (2,-1.8) node {$G \wedge uw \setminus u$};
    \end{tikzpicture}
    \caption{$\theta(\ell_1,\ell_2,\ell_3)$ with $3 \le \ell_1 \le \ell_2 \le \ell_3$ and its vertex-minors.}
    \label{fig:tp2i}
  \end{figure}

  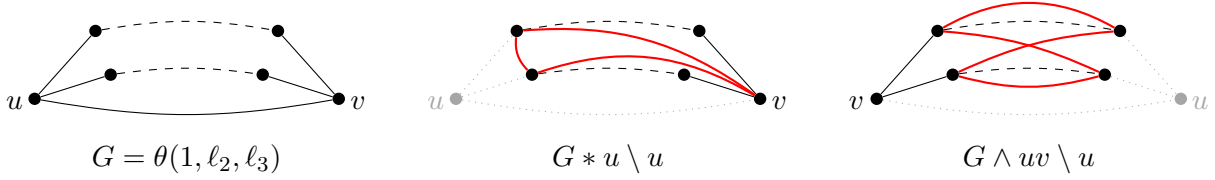
\begin{figure}
    \centering
    \begin{tikzpicture}
      \node[shape=circle,fill=black, scale=0.45] (u) at (0,0) {};
      \node[shape=circle,fill=black, scale=0.45] (v) at (4,0) {};

      \node[shape=circle,fill=black, scale=0.45] (x11) at (0.8,0.9) {};
      \node[shape=circle,fill=black, scale=0.45] (x14) at (3.2,0.9) {};

      \node[shape=circle,fill=black, scale=0.45] (x21) at (1,0.32) {};
      \node[shape=circle,fill=black, scale=0.45] (x23) at (3,0.32) {};

      \draw (u)+(-0.265,-0.05) node {$u$};
      \draw (v)+(0.255,-0.05) node {$v$};

      \draw (u) to [bend right=10] (v);

      \draw (u) -- (x11);
      \draw (u) -- (x21);
      
      \draw (x14) -- (v);
      \draw (x23) -- (v);

      \draw[dashed] (x11) to [bend left = 10] (x14);
      \draw[dashed] (x21) to [bend left = 8] (x23);
  
      \draw (2,-0.8) node {$G = \theta(1,\ell_2,\ell_3)$};
    \end{tikzpicture}
    \hspace{0.3cm}
    \begin{tikzpicture}
      \node[shape=circle,fill=black, scale=0.45, opacity=0.35] (u) at (0,0) {};
      \node[shape=circle,fill=black, scale=0.45] (v) at (4,0) {};

      \node[shape=circle,fill=black, scale=0.45] (x11) at (0.8,0.9) {};
      \node[shape=circle,fill=black, scale=0.45] (x14) at (3.2,0.9) {};

      \node[shape=circle,fill=black, scale=0.45] (x21) at (1,0.32) {};
      \node[shape=circle,fill=black, scale=0.45] (x23) at (3,0.32) {};

      \draw[opacity=0.35] (u)+(-0.265,-0.05) node {$u$};
      \draw (v)+(0.255,-0.05) node {$v$};

      \draw[opacity=0.35,dotted] (u) to [bend right=10] (v);

      \draw[opacity=0.35,dotted] (u) -- (x11);
      \draw[opacity=0.35,dotted] (u) -- (x21);
      
      \draw (x14) -- (v);
      \draw (x23) -- (v);

      \draw[dashed] (x11) to [bend left = 10] (x14);
      \draw[dashed] (x21) to [bend left = 8] (x23);

      \draw[color=red,thick] (x11) to [bend right=25] (x21);

      \draw[color=red,thick] (x11) to [bend left=20] (v);
      \draw[color=red,thick] (x21) to [bend left=25] (v);

      \draw (2,-0.8) node {$G*u \setminus u$};
    \end{tikzpicture}
    \hspace{0.3cm}
    \begin{tikzpicture}
      \node[shape=circle,fill=black, scale=0.45] (u) at (0,0) {};
      \node[shape=circle,fill=black, scale=0.45, opacity=0.35] (v) at (4,0) {};

      \node[shape=circle,fill=black, scale=0.45] (x11) at (0.8,0.9) {};
      \node[shape=circle,fill=black, scale=0.45] (x14) at (3.2,0.9) {};

      \node[shape=circle,fill=black, scale=0.45] (x21) at (1,0.32) {};
      \node[shape=circle,fill=black, scale=0.45] (x23) at (3,0.32) {};

      \draw (u)+(-0.265,-0.05) node {$v$};
      \draw[opacity=0.35] (v)+(0.255,-0.05) node {$u$};

      \draw[opacity=0.35,dotted] (u) to [bend right=10] (v);

      \draw (u) -- (x11);
      \draw (u) -- (x21);
      
      \draw[opacity=0.35,dotted] (x14) -- (v);
      \draw[opacity=0.35,dotted] (x23) -- (v);

      \draw[dashed] (x11) to [bend left = 10] (x14);
      \draw[dashed] (x21) to [bend left = 8] (x23);

      \draw[color=red,thick] (x11) to [bend left=30] (x14);
      \draw[color=red,thick] (x21) to [bend right=15] (x23);

      \draw[color=red,thick] (x11) to [bend left=10] (x23);
      \draw[color=red,thick] (x21) to [bend left=10] (x14);

      \draw (2,-0.8) node {$G\wedge uv \setminus u$};
    \end{tikzpicture}
    \caption{$\theta(1,\ell_2,\ell_3)$ with $4 \le \ell_2 \le \ell_3$ and its vertex-minors.}
    \label{fig:tp2ii}
  \end{figure}
  
  Next we consider the case that $m=3$ and $\ell_1 = 1$; see Figure~\ref{fig:tp2ii}.
  Then $\ell_2 \geq 4$ and $w=v$ because $v$ is the neighbor of $u$ in $P_1$.
  Both $G*u\setminus u$ and $G \wedge uv \setminus u \setminus v$ are $2$-connected outerplanar and thus they are prime by Lemma~\ref{lem:2outer}.
  Since $v$ has degree $2$ and has no twin in $G \wedge uv \setminus u$, by Lemma~\ref{lem:pext}\ref{item:p1}, $G \wedge uv \setminus u$ is prime.

  Now it remains to consider the case that $m \geq 4$.
  For each $2 \le i \le 4$, as $G-V_i$ is isomorphic to $\theta(\ell_1,\dots,\ell_{i-1},\ell_{i+1},\dots,\ell_m)$, 
  we deduce that both $(G*u \setminus u) -V_i = (G-V_i)*u \setminus u$ and $(G\wedge uw \setminus u)-V_i = (G-V_i)\wedge uw \setminus u$ are prime by the inductive hypothesis.
  Thus by Lemma~\ref{lem:primepart}, both $G*u \setminus u$ and $G \wedge uw \setminus u$ are prime.
\end{proof}

\begin{lemma}\label{lem:tp3}
  Let $G$ be $\theta(1,2,\ell_1,\dots,\ell_m)$ or $\theta(2,\ell_1,\dots,\ell_m)$ with $m\geq 2$ and $\min\{ \ell_1, \dots, \ell_m \} \geq 3$.
  For a vertex $x$ of $G$, the following are equivalent: (i) $x$ has degree larger than $2$ or has no neighbor of degree $2$, (ii) $x$ is non-pivotal, and (iii) $x$ is non-essential.
  In particular, $G$ has exactly three non-essential vertices.
\end{lemma}
\begin{proof}
  Note that $|V(G)| = 3 + \sum_{i=1}^m (\ell_i - 1) \geq 7$.
  By definition, (iii) implies (ii) and by Lemma~\ref{lem:essv}, (ii) implies (i).
  Thus, it suffices to show that (i) implies (iii).
  Let $u$ and $v$ be the two vertices of degree at least $3$ in $G$ and let $w$ be the common neighbor of $u$ and~$v$.
  Note that $w$ is the unique vertex that has degree $2$ and has no neighbor of degree $2$.
  We claim that $u$, $v$, and $w$ are non-essential.

  Since $\theta(1,2,\ell_1,\dots,\ell_m)$ and $\theta(2,\ell_1,\dots,\ell_m)$ are locally equivalent by applying a local complementation at $w$, we may assume that $G = \theta(2,\ell_1,\dots,\ell_m)$.
  Observe that $G \setminus w$ is isomorphic to $\theta(\ell_1, \dots, \ell_m)$, and $G * w \setminus w$ is isomorphic to $\theta(1, \ell_1, \dots, \ell_m)$.
  By Proposition~\ref{prop:theta1}, $G \setminus w$ and $G * w \setminus w$ are prime and therefore $w$ is non-essential in~$G$.

  Since a cycle graph of length at least $5$ has no non-essential vertex, $G$ is not locally equivalent to a cycle graph.
  By Theorem~\ref{thm:noness2}, $G$ has at least two non-essential vertices.
  Therefore, $u$ or $v$ is non-essential in~$G$, and by symmetry, both $u$ and $v$ are non-essential in $G$.
\end{proof}

A graph $\theta(\ell_1,\ell_2)$ is a cycle of length $\ell_1+\ell_2$, which is prime and has no non-essential vertex if $\ell_1+\ell_2 \geq 5$.
By Proposition~\ref{prop:theta1}, for positive integers $m \geq 3$ and $\ell_1 \leq \dots \leq \ell_m$ with $\ell_2 \geq 2$, a graph $\theta(\ell_1,\dots,\ell_m)$ is prime if and only if either (i) $\ell_2 \geq 3$ or (ii) $\ell_1 = 1$, $\ell_2 = 2$, and $\ell_3 \geq 3$.
In the next proposition, we 
determine the number of non-essential vertices in 
$\theta(\ell_1,\dots,\ell_m)$
when it is prime and $m \ge 3$.

\begin{proposition}\label{prop:theta2}
  Let $m$ and $\ell_1,\dots,\ell_m$ be positive integers with $m\geq 3$, $\ell_1 \leq \ell_2 \leq \dots \leq \ell_m$, and $\ell_3 \geq 3$ such that either $\ell_2 \geq 3$ or $(\ell_1, \ell_2) = (1,2)$.
  Let $G = \theta(\ell_1,\dots,\ell_m)$.
  \begin{enumerate}[label=\rm(\arabic*)]
    \item\label{item:th0} If $\ell_1=1$, $m=3$, and $\ell_2 \leq 3$, then $G$ is locally equivalent to a cycle of length $\ell_2+\ell_3$ and has no non-essential vertex.
    \item\label{item:th21} If $\ell_1 = 1$, $m=3$, and $\ell_2 \geq 4$, then $G$ has exactly $2$ non-essential vertices.
    \item\label{item:th31} If $\ell_1=1$, $m \geq 4$, and $\ell_2=2$, then $G$ has exactly $3$ non-essential vertices.
    \item\label{item:th22} If $\ell_1 = 1$, $m \geq 4$, and $\ell_2 \geq 3$, then $G$ has exactly $2$ non-essential vertices.
    \item\label{item:th32} If $\ell_1 = 2$, then $G$ has exactly $3$ non-essential vertices.
    \item\label{item:th23} If $\ell_1 \geq 3$, then $G$ has exactly $2$ non-essential vertices.
  \end{enumerate}
\end{proposition}

\begin{proof}
  It is obvious to check \ref{item:th0}.
  Lemma~\ref{lem:tp2} implies~\ref{item:th21},~\ref{item:th22}, and~\ref{item:th23}.
  Lemma~\ref{lem:tp3} implies~\ref{item:th31} and~\ref{item:th32}.
\end{proof}

By the preceding proposition, every graph in $\Theta$ has at most $2$ non-essential vertices.
We now prove Theorem~\ref{thm:noness3} from Propositions~\ref{prop:gn2} and~\ref{prop:theta2}.

\begin{proof}[Proof of Theorem~\ref{thm:noness3}]
  Let $G$ be a prime graph with at least four vertices.
  Then $|V(G)| \ge 5$ because no graph on four vertices is prime.
  Suppose that $G$ is locally equivalent to a graph $H \in \Theta$ 
  consisting of $m$ internally-disjoint paths 
  between two fixed distinct vertices having no common neighbors, where $m\geq 2$.
  If $m = 2$, then $H$ is a cycle and thus $G$ has no non-essential vertex.
  If $m\geq 3$, then by Proposition~\ref{prop:theta2}, $G$ has at most $2$ non-essential vertices.

  Now, we prove the backward direction.
  Suppose that $G$ has at most $2$ non-essential vertices.
  By Proposition~\ref{prop:gn2}, $G$ is locally equivalent to a graph isomorphic to $\theta(\ell_1,\dots,\ell_m)$ for some $m$ and $\ell_i$.
  Since $G$ is prime, $m\geq 2$.
  We may assume that $G$ is not locally equivalent to a cycle graph because a cycle graph of length $k$ is isomorphic to $\theta(1,k-1)$.
  Thus $m \geq 3$.
  By Proposition~\ref{prop:theta2}, we conclude that $\ell_i \neq 2$ for all~$i$.
\end{proof}

\section{Pivot-minors and non-pivotal vertices}\label{sec:nonpiv}

We prove Theorems~\ref{thm:nonpiv2} and~\ref{thm:nonpiv3}, which are analogues of Theorems~\ref{thm:noness2} and~\ref{thm:noness3} for pivot-minors.
We also prove Corollaries~\ref{cor:bippiv} and~\ref{cor:bippiv3}.
We present useful results first.

\begin{figure}
  \centering
  \begin{tikzpicture}
    \begin{scope}[xshift=-5cm]
    \node[shape=circle,fill=black, scale=0.35] (1) at (360/5*0+18:1.2) {};
    \node[shape=circle,fill=black, scale=0.35] (2) at (360/5*1+18:1.2) {};
    \node[shape=circle,fill=black, scale=0.35] (3) at (360/5*2+18:1.2) {};
    \node[shape=circle,fill=black, scale=0.35] (4) at (360/5*3+18:1.2) {};
    \node[shape=circle,fill=black, scale=0.35] (5) at (360/5*4+18:1.2) {};
    \draw (1) -- (2) -- (3) -- (4) -- (5) -- (1);
    \node () at (0,-1.45) {$C_5$};
    \end{scope}

    \draw[<-,dashed] (-3.5,0) -- (-1.5,0) node[midway,above] {\scriptsize{pivoting $e$}};

    \node[shape=circle,fill=black, scale=0.35] (1) at (360/5*0+18:1.2) {};
    \node[shape=circle,fill=black, scale=0.35] (2) at (360/5*1+18:1.2) {};
    \node[shape=circle,fill=black, scale=0.35] (3) at (360/5*2+18:1.2) {};
    \node[shape=circle,fill=black, scale=0.35] (4) at (360/5*3+18:1.2) {};
    \node[shape=circle,fill=black, scale=0.35] (5) at (360/5*4+18:1.2) {};
    \draw (1) -- (2) -- (3) -- (4) -- (5);
    \draw (2) -- (4);
    \draw (1) -- (5) node[midway,left] {$e$};
    \node () at (0,-1.45) {$G_1$};
  
    \draw[<-,dashed] (1.5,0) -- (3.5,0) node[midway,above] {\scriptsize{pivoting $e'$}};

    \begin{scope}[xshift=5cm]
    \node[shape=circle,fill=black, scale=0.35] (1) at (360/5*0+18:1.2) {};
    \node[shape=circle,fill=black, scale=0.35] (2) at (360/5*1+18:1.2) {};
    \node[shape=circle,fill=black, scale=0.35] (3) at (360/5*2+18:1.2) {};
    \node[shape=circle,fill=black, scale=0.35] (4) at (360/5*3+18:1.2) {};
    \node[shape=circle,fill=black, scale=0.35] (5) at (360/5*4+18:1.2) {};
    \draw (1) -- (2) -- (3);
    \draw (4) -- (5) -- (1);
    \draw (2) -- (4);
    \draw (2) -- (5);
    \draw (3) -- (4) node[midway,left] {$e'$};
    \node () at (0,-1.45) {$G_2$};
    \end{scope}
  \end{tikzpicture}
  \caption{All graphs locally equivalent to $C_5$ up to isomorphism.}
  \label{fig:locC5}
\end{figure}
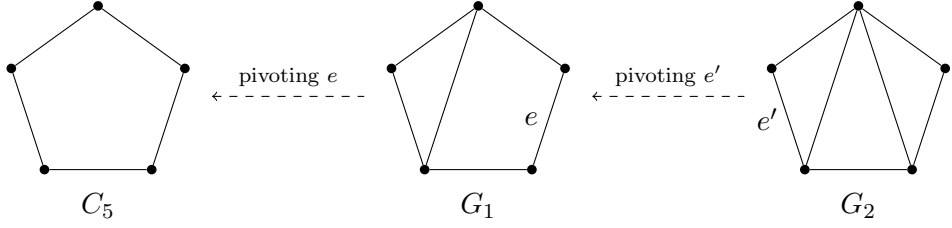

\begin{lemma}
  \label{lem:locC5}
  A graph is locally equivalent to a cycle of length $5$ if and only if it is pivot-equivalent to a cycle of length $5$.
\end{lemma}
\begin{proof}
  The backward direction is trivial.
  The forward direction is easily seen by Figure~\ref{fig:locC5} which depicts all graphs locally equivalent to a cycle of length $5$ up to isomorphism.
  In Figure~\ref{fig:locC5}, $G_1 \wedge e$ is isomorphic to $C_5$, and $G_2 \wedge e'$ is isomorphic to $G_1$.
  Therefore $C_5$, $G_1$, and $G_2$ are pivot-equivalent.
\end{proof}

\begin{theorem}
  \label{thm:plp}
  Let $S = (V,L)$ be an isotropic system, and let $(G,\mb{a},\mb{b})$ and $(H,\mb{c},\mb{d})$ be graphic presentations of $S$.
  Then there are nonnegative integers $m$, $k$, $\ell$, vertices $v_1$, $v_2$, $\ldots$, $v_m$, and edges $e_1$, $e_2$, $\ldots$, $e_k$, $e_1'$, $e_2'$, $\ldots$, $e_\ell'$ such that the following hold.
  \begin{enumerate}[label=\rm(\roman*)]
    \item\label{item:plp1} For $1\le i\le k$, $e_i$ is an edge of $G\wedge e_1\wedge \cdots \wedge e_{i-1}$.
    \item\label{item:plp2} $\{v_1,v_2,\ldots,v_m\}$ is an independent set of vertices  in $G \wedge e_1 \wedge \dots \wedge e_k$.
    \item\label{item:plp3} For $1\le j\le \ell$, $e_j'$ is an edge of $G\wedge e_1\wedge \cdots \wedge e_{k}*v_1*\cdots *v_m \wedge e_1'\wedge \cdots \wedge e_{j-1}'$.
    \item\label{item:plp4} $(H,\mb{c},\mb{d}) = (G,\mb{a},\mb{b}) \wedge e_1\wedge \dots \wedge e_k * v_1* \dots * v_m \wedge e_{1}'\wedge \dots \wedge e_\ell'$
  \end{enumerate}
\end{theorem}

Theorem~\ref{thm:plp} is a slight strengthening of Fon-Der-Flaass~\cite{FonDerFlaass1989} in Russian; see also~{\cite[Theorem~3.4]{FonDerFlaass1996}}.
His theorem states that for two locally equivalent graphs $G$ and $H$, there are vertices $v_1,\dots,v_m$ and edges $e_1,\dots,e_k,e_1',\dots,e_{\ell}'$ satisfying \ref{item:plp1}--\ref{item:plp3} and the following replacement of~\ref{item:plp4}: 
\begin{enumerate}
  \item[\textrm{(iv$'$)}] $H = G \wedge e_1 \wedge \cdots \wedge e_k * v_1 * \cdots * v_m \wedge e_1 \wedge \cdots \wedge e_{\ell}'$.
\end{enumerate}
Our proof of Theorem~\ref{thm:plp} uses the \emph{divergence}. %
This technique was introduced by Fon-Der-Flaass~\cite{FonDerFlaass1988} for graphs, and used for isotropic systems by Bouchet~\cite{Bouchet1991algo}.

\begin{proof}%
  Let $(G',\mb{a}',\mb{b}')$ and $(H',\mb{c}',\mb{d}')$ be graphic presentations pivot-equivalent to $(G,\mb{a},\mb{b})$ and $(H,\mb{c},\mb{d})$, respectively.
  For $v\in V$, let
  \begin{align*}
    d_{G',H'}(v) :=
    \begin{cases}
      0   & \text{if $\mb{c}'(v) = \mb{a}'(v)$,} \\
      1   & \text{if $\mb{c}'(v) = \mb{b}'(v)$,} \\
      2   & \text{otherwise.} \\
    \end{cases}
  \end{align*}
  and let $D(G',H') := \sum_{v\in V} d_{G',H'}(v)$.
  For each $i\in\{0,1,2\}$, let $A_i$ be the set of vertices $v\in V$ such that $d_{G',H'}(v) = i$.

  We take $G'$ and $H'$ minimizing $D(G',H')$.
  To complete the proof, it suffices to show that $(H',\mb{c}',\mb{d}') = (G',\mb{a}',\mb{b}') * v_1 \dots * v_m$ for some independent set $\{v_1,\dots,v_m\}$ in $G'$.
  We present four claims step by step.

  \medskip\noindent
  \textbf{Claim I.} 
  No vertex in $A_1$ is adjacent to vertices in $A_1\cup A_2$ in $G'$.
  Suppose that there is $vw \in E(G')$ for some $v\in A_1$ and $w \in A_1 \cup A_2$.
  By Proposition~\ref{prop:ig2},  
  \begin{align*}
    d_{G'\wedge vw,H'}(v) &= 0,\\
    d_{G'\wedge vw, H'}(w) &=
    \begin{cases}
      0   & \text{if $w \in A_1$,} \\
      2   & \text{otherwise,} \\
    \end{cases}
    \\
    d_{G'\wedge vw, H'}(x) &= d_{G', H'}(x) &\text{for all }x \in V - \{v,w\}.
  \end{align*}
  Therefore, $D(G'\wedge vw,H') \leq D(G',H') - 1$, which contradicts our choice of $G'$ and $H'$. This proves Claim I.

  \medskip\noindent 
  \textbf{Claim II.} $A_1 = \emptyset$.
  Suppose that $A_1$ has a vertex $v$.
  By Claim I, $N_{G'}(v) \subseteq A_0$.
  Since $(G',\mb{a}',\mb{b}')$ is a graphic presentation of $S = (V,L)$, $L$ has a vector $\mb{b}'_v$ such that
  \begin{align*}
    \mb{b}'_v(w) =
    \begin{cases}
      \mb{b}'(v)  & \text{if $w=v$,} \\
      \mb{a}'(w)  & \text{if $w$ is a neighbor of $v$ in $G'$,} \\
      0           & \text{otherwise.} \\
    \end{cases}
  \end{align*}
  Then $\mb{c}'[\{v\} \cup N_{G'}(v)] = \mb{b}'_v \in L$, violating that $\mb{c}'$ is an Eulerian vector of $S$. This proves Claim II.

  \medskip 
  Note that for each $v\in A_2$, $\mb{d}'(v)$ is either $\mb{a}'(v)$ or $\mb{b}'(v)$ because $\mb{c}'(v) = \mb{a}'(v) + \mb{b}'(v)$.
  Let $B$ be the set of vertices $v\in A_2$ such that $\mb{d}'(v) = \mb{a}'(v)$. 

  \medskip\noindent 
  \textbf{Claim III.} The set $A_2$ is independent in $H'$.
  Suppose that there is $vw \in E(H')$ for some $v,w \in A_2$.
  By Proposition~\ref{prop:ig2}, for each $x\in \{v,w\}$,
  \begin{align*}
    d_{G',H'\wedge vw}(x) =
    \begin{cases}
      0   & \text{if $x \in B$,} \\
      1   & \text{otherwise,} \\
    \end{cases}
  \end{align*}
  and for each $y \in V - \{v,w\}$, $d_{G', H'\wedge vw}(y) = d_{G', H'}(y)$.
  Therefore, $D(G',H'\wedge vw) \leq D(G',H') - 2$, a~contradiction.
  This proves Claim III.

  \medskip\noindent
  \textbf{Claim IV.} $B = \emptyset$.
  Suppose that $B$ has a vertex $v$.
  By Claims II and III, $N_{H'}(v) \subseteq A_0$.
  Since $(H',\mb{c}',\mb{d}')$ is a graphic presentation of $S$, there is a vector $\mb{d}'_v$ in $L$ such that
  \begin{align*}
    \mb{d}'_v(w) =
    \begin{cases}
      \mb{d}'(v)  & \text{if $w=v$,} \\
      \mb{c}'(w)  & \text{if $w$ is a neighbor of $v$ in $H'$,} \\
      0           & \text{otherwise.} \\
    \end{cases}
  \end{align*}
  Then $\mb{a}'[\{v\} \cup N_{H'}(v)] = \mb{d}'_v \in L$, which contradicts that $\mb{a}'$ is an Eulerian vector of $S$.
  This proves Claim IV.

  \medskip 
  In conclusion, $V - A_0$ is independent in $H'$ and for each $v \in V - A_0$, $\mb{c}'(v) = \mb{a}'(v) + \mb{b}'(v)$ and $\mb{d}'(v) =  \mb{b}'(v)$.
  By Propositions~\ref{prop:ig2} and~\ref{prop:ig}(ii), $(G',\mb{a}',\mb{b}') = (H',\mb{c}',\mb{d}') * v_m \dots * v_1$ where $\{v_1,\dots,v_m\} = V - A_0$.
  Equivalently, $(H',\mb{c}',\mb{d}') = (G',\mb{a}',\mb{b}') * v_1 \dots * v_m$.
  It is easy to check that $\{v_1,\dots,v_m\}$ is independent in $G' = H' * v_m \dots * v_1$.
\end{proof}

\begin{corollary}
  \label{cor:gamma}
  Let $G$ be a prime graph with at least $5$ vertices and let $S$ be an isotropic system associated with a graphic presentation $(G,\mb{a},\mb{b})$ such that $\{\mb{a}(v), \mb{b}(v)\} = \{\alpha, \beta\}$ for all $v\in V(G)$.
  For each graphic presentation of $(H,\mb{c},\mb{d})$ of $S$, either $\{v\in V(G): \mb{c}(v) = \gamma \text{ or } \mb{d}(v) = \gamma \}$ is empty or has at least $3$ vertices.
\end{corollary}
\begin{proof}
  Denote $W:=\{v\in V(G): \mb{c}(v) = \gamma \text{ or } \mb{d}(v) = \gamma \}$.
  By Proposition~\ref{prop:ig2} and Theorem~\ref{thm:plp}, we may assume that $(H,\mb{c},\mb{d}) = (G,\mb{a},\mb{b})  * v_1 * \dots * v_m$ for some distinct and pairwisely non-adjacent vertices $v_1,\dots,v_m$ in $G$.
  By Proposition~\ref{prop:ig2}, $\mb{c}(v_i) = \gamma$ for all $i\in[m]$.
  Thus, we may assume that $m \leq 2$.
  Suppose that $m = 1$.
  Then for each $w \in N_G(v_1)$, $\mb{d}(w) = \gamma$ by Proposition~\ref{prop:ig2} and therefore $W = \{v_1\} \cup N_G(v_1)$.
  Since $G$ is prime and $|V(G)| \geq 5$, the minimum degree of $G$ is $2$ or more, which implies that $|W| \geq 3$.
  Hence we may assume that $m = 2$.
  Because $G$ has no twins, $|N_{G}(v_1) \triangle N_{G}(v_2)| \geq 1$.
  By Proposition~\ref{prop:ig2}, $W = \{v_1,v_2\} \cup (N_{G}(v_1) \triangle N_{G}(v_2))$ and therefore $|W| \geq 3$.
\end{proof}

The following corollary was conjectured by Bouchet~\cite{Bouchet1990} and proved by Fon-Der-Flaass~\cite{FonDerFlaass1989,FonDerFlaass1996}.
It will be used to prove Corollary~\ref{cor:bippiv}.
We remark that this corollary is also a consequence of the following facts on binary matroids and their relation to bipartite graphs. Seymour~\cite{Seymour1988} showed that if two connected binary matroids on the same ground set have identical connectivity functions, then they are equal up to duality. 
It is well known that a binary matroid is uniquely determined by its fundamental graph, which is bipartite.
Oum~\cite{Oum2005} observed that 
the connectivity function of a binary matroid is precisely the cut-rank function of its fundamental graph and locally equivalent graphs have the same cut-rank functions. 

\begin{corollary}[Fon-Der-Flaass~\cite{FonDerFlaass1989,FonDerFlaass1996}]\label{cor:bipLocPiv}
  If two bipartite graphs are locally equivalent, then they are pivot-equivalent.
\end{corollary}

For the convenience of readers, we include a proof by Fon-Der-Flaass using Theorem~\ref{thm:plp}.

\begin{proof}
  Let $G$ and $H$ be locally equivalent bipartite graphs.
  By Theorem~\ref{thm:plp}, $H = G \wedge e_1\wedge \dots \wedge e_k * v_1* \dots * v_m \wedge e_{1}'\wedge \dots \wedge e_\ell'$
  for some edges $e_1$, $\ldots$, $e_k$, $e_1'$, $\ldots$, $e_\ell'$
  and an independent set $\{v_1,\dots,v_m\}$ of $G\wedge e_1\wedge \dots \wedge e_k$.
  Let $G':=G\wedge e_1\wedge \dots \wedge e_k$ and $H':=H\wedge e_\ell'\wedge\cdots \wedge e_1'$. 
  Observe that both $G'$ and $H'$ are bipartite and $H'=G'*v_1*\cdots *v_m$.
  Now it is easy to see that $G'=H'$.
\end{proof}

\begin{lemma}[Oum~{\cite[Proposition~10.1]{Oum2008wqo}}]
  \label{lem:abminor}
  Let $(G_1,\mb{a}_1,\mb{b}_2)$ and $(G_2,\mb{a}_2,\mb{b}_2)$ be graphic presentations of an isotropic system.
  If $\{\mb{a}_1(v),\mb{b}_1(v)\} = \{\mb{a}_2(v),\mb{b}_2(v)\}$ for each vertex $v$, then $G_1$ and $G_2$ are pivot-equivalent.
\end{lemma}

Now we prove Theorem~\ref{thm:nonpiv2} using Theorem~\ref{thm:noness2} together with preceding results.

\begin{proof}[Proof of Theorem~\ref{thm:nonpiv2}]
  Let $G$ be a prime graph with at least four vertices which has at most $1$ non-pivotal vertex.
  Denote $V:= V(G)$.
  Then $|V| \ge 5$ since no graph on four vertices is prime.
  If $|V|=5$, then $G$ is locally equivalent to a cycle and therefore it is pivot-equivalent to a cycle by Lemma~\ref{lem:locC5}.
  Thus, we may assume that $G$ has at least $6$ vertices.
  
  Since every non-essential vertex is non-pivotal, $G$ has at most $1$ non-essential vertex.
  By Theorem~\ref{thm:noness2}, $G$ is locally equivalent to a cycle~$H$.
  It is enough to show that $G$ is pivot-equivalent to~$H$.

  Let $S$ be an isotropic system associated with a graphic presentation $(G,\mb{a},\mb{b})$
  where $\mb{a}, \mb{b} \in K^V$ are supplementary vectors 
  such that $\mb{a}(v) = \alpha$ and $\mb{b}(v) = \beta$ for each $v\in V$.
  Since $H$ is locally equivalent to $G$, 
  there exist supplementary vectors $\mb{c}$ and $\mb{d}$
  such that $(H,\mb{c},\mb{d})$ is a graphic presentation of $S$.
  A vertex $v$ is non-pivotal in $G$ if and only if $S|^v_\alpha$ or $S|^v_\beta$ is $3$-connected by Proposition~\ref{prop:igmin}.

  For every vertex $v$, a vertex-minor $H*v\setminus v$ is prime and 
  therefore $S|^v_{\mb{c}(v)+\mb{d}(v)}$ is $3$-connected by Proposition~\ref{prop:igmin}.
  So if $\mb{c}(v)+\mb{d}(v)\in\{\alpha,\beta\}$, then $v$ is non-pivotal in $G$.
  Since $G$ has at most $1$ non-pivotal vertex,
  by Corollary~\ref{cor:gamma}, 
  $\mb{c}(v)+\mb{d}(v)\notin\{\alpha,\beta\}$ for all $v\in V(G)$
  and therefore
  $\{\mb{c}(v),\mb{d}(v)\}=\{\alpha,\beta\}$. 
  By Lemma~\ref{lem:abminor}, $G$ and $H$ are pivot-equivalent.
\end{proof}

\begin{proof}[Proof of Theorem~\ref{thm:nonpiv3}]
  By Lemma~\ref{lem:tp2}, every graph in $\Theta$ has at most $2$ non-pivotal vertices.
  Thus, it suffices to show that if a prime graph $G$ with at least $5$ vertices has at most $2$ non-pivotal vertices, then $G$ is pivot-equivalent to a graph in $\Theta$.
  Denote $V:=V(G)$.
  Then $|V| \ge 5$ because no graph on four vertices is prime.
  If $\abs{V}=5$, then $G$ is locally equivalent to a cycle and therefore it is pivot-equivalent to a cycle by Lemma~\ref{lem:locC5}.
  Thus, we may assume that $G$ has at least $6$ vertices.

  Because $G$ has at most $2$ non-essential vertices, by Theorem~\ref{thm:noness3}, $G$ is locally equivalent to a graph $H \in \Theta$ consisting of internally-disjoint paths $P_1,\dots,P_m$ between two fixed vertices $x$ and $y$ such that $m\geq 2$ and no $P_i$ has length $2$.

  Let $S$ be an isotropic system associated with a graphic presentation $(G,\mb{a},\mb{b})$ where $\mb{a}, \mb{b} \in K^V$ are supplementary vectors such that $\mb{a}(v) = \alpha$ and $\mb{b}(v) = \beta$ for each $v\in V$.
  Since $G$ and $H$ are locally equivalent, 
  there are supplementary vectors $\mb{c}$ and $\mb{d}$ such that $(H,\mb{c},\mb{d})$ is a graphic presentation of $S$.
  By Proposition~\ref{prop:igmin}, a vertex $v$ is non-pivotal in $G$ if and only if $S|^v_\alpha$ or $S|^v_\beta$ is $3$-connected.

  If $v$ is a vertex of degree~$2$ in $H$, then 
  $H*v\setminus v$ is prime 
  by Proposition~\ref{prop:theta1}.
  Therefore, $S|^v_{\mb{c}(v)+\mb{d}(v)}$ is $3$-connected 
  by Proposition~\ref{prop:igmin}.
  So if $\mb{c}(v)+\mb{d}(v)\in\{\alpha,\beta\}$, then $v$ is non-pivotal in $G$.

  If $H$ is a cycle, then 
  since $G$ has at most $2$ non-pivotal vertices, 
  by Corollary~\ref{cor:gamma}, 
  $\mb{c}(v)+\mb{d}(v)\notin\{\alpha,\beta\}$ for all $v\in V$
  and therefore
  $\{\mb{c}(v),\mb{d}(v)\}=\{\alpha,\beta\}$.

  If $H$ is not a cycle, then 
  by Lemma~\ref{lem:tp2}, $H$ has exactly $2$ non-essential vertices that are $x$ and~$y$.
  Then $x$ and $y$ are non-essential in $G$ and thus they are non-pivotal in $G$. Therefore no vertex of degree~$2$ in $H$ is non-pivotal in $G$ and so 
  $\{\mb{c}(v),\mb{d}(v)\}=\{\alpha,\beta\}$ for every vertex $v$ of degree~$2$ in $H$.
  By Corollary~\ref{cor:gamma}, 
  $\{\mb{c}(x),\mb{d}(x)\} = \{\mb{c}(y),\mb{d}(y)\}= \{\alpha,\beta\} $.
  
  In both cases, it implies that $G$ is pivot-equivalent to $H\in\Theta$
  by Lemma~\ref{lem:abminor}.
\end{proof}

Here is an easy observation on bipartite graphs.
\begin{lemma}\label{lem:pivot and bipartite}
  If $G$ is bipartite and $uv$ is an edge of $G$, 
  then $G\wedge uv$ is bipartite.
  \qed
\end{lemma}

Finally, we are ready to prove Corollaries~\ref{cor:bippiv} and~\ref{cor:bippiv3}.

\begin{proof}[Proof of Corollary~\ref{cor:bippiv} using Theorem~\ref{thm:nonpiv2}]
  Let $G$ be a prime bipartite graph such that $|V(G)| \ge 4$ and $G$ has fewer than two non-pivotal vertices.
  Then by Theorem~\ref{thm:nonpiv2}, $G$ is pivot-equivalent to a cycle~$C$.
  By Lemma~\ref{lem:pivot and bipartite}, $C$ is bipartite and so it is an even cycle.
\end{proof}

\begin{proof}[Proof of Corollary~\ref{cor:bippiv3} using Theorem~\ref{thm:nonpiv3}]
  Let $G$ be a prime bipartite graph such that $|V(G)| \ge 4$ and $G$ has fewer than three non-pivotal vertices.
  Then by Theorem~\ref{thm:nonpiv3} and Lemma~\ref{lem:pivot and bipartite}, $G$ is pivot-equivalent to some bipartite graph in $\Theta$.
\end{proof}

In the remainder of this section, we will show that Corollary~\ref{cor:bippiv} can be deduced directly from Theorem~\ref{thm:noness2} 
without using Theorem~\ref{thm:nonpiv2}.
For that, we will need several properties of bipartite graphs.

\begin{lemma}[Allys~{\cite[Lemma~5.2]{Allys1994iso}}]\label{lem:oddcyc}
  No bipartite graph is locally equivalent to an odd cycle of length at least five.
\end{lemma}

Allys proved the preceding lemma by using isotropic systems.
Lemma~\ref{lem:oddcyc} is also implied by the following theorem of Fon-Der-Flaass~\cite{FonDerFlaass1988}, which was published earlier than~\cite{Allys1994iso}.
We provide another proof of Lemma~\ref{lem:oddcyc} in Appendix~\ref{app:oddcyc}.

\begin{theorem}[Fon-Der-Flaass~{\cite[Theorem~5.1]{FonDerFlaass1988}}]\label{thm:Hamiltonian}
  Every graph locally equivalent to a cycle of length at least five is Hamiltonian.
\end{theorem}

Allys stated that Lemma~\ref{lem:oddcyc} can be used to show that Theorem~4.3 of~{\cite{Allys1994iso}} implies the wheel and whirl theorem of Tutte~{\cite[8.2]{Tutte1966cm}} for binary matroids.
However, we believe that the next lemma is necessary to complete this implication.

\begin{lemma}\label{lem:evencyc}
  A bipartite graph with at least five vertices is locally equivalent to a cycle if and only if it is pivot-equivalent to an even cycle.
\end{lemma}

\begin{proof}
  Since the backward direction is trivial, we now prove the forward direction.
  Let $G$ be a bipartite graph locally equivalent to a cycle $C$ of length at least $5$.
  By Lemma~\ref{lem:oddcyc}, $C$ must be an even cycle and by Corollary~\ref{cor:bipLocPiv}, $G$ and $C$ are pivot-equivalent.
\end{proof}

\begin{proof}[Proof of Corollary~\ref{cor:bippiv} using Theorem~\ref{thm:noness2}]
  It is straightforward from Theorem~\ref{thm:noness2} and Lemma~\ref{lem:evencyc}.
\end{proof}

\section{Proofs of Applications}\label{sec:applications}

We present proofs of Corollaries~\ref{cor:rooted prime reduction1},~\ref{cor:rooted prime reduction2}, and~\ref{cor:rooted prime reduction0}.

\begin{proof}[Proof of Corollary~\ref{cor:rooted prime reduction1}]
  It suffices to show that $G$ has a vertex $v \in V(G) - \{x,y\}$ such that $G \setminus v$, $G*v \setminus v$, or $G / v$ is prime.
  By Theorem~\ref{thm:noness3}, we may assume that $G$ is locally equivalent to a graph $G'$ in $\Theta$.
  As $|V(G')| = |V(G)| \ge 6$ and $G'$ has at most two vertices of degree larger than $2$, there is a vertex $v \in V(G') - \{x,y\}$ such that $\deg_{G'}(v) = 2$.
  By Proposition~\ref{prop:theta1}, $(G' * v) \setminus v$ is prime and so $G \setminus v$, $G*v \setminus v$, or $G / v$ is prime.
\end{proof}

\begin{proof}[Proof of Corollary~\ref{cor:rooted prime reduction2}]
  It suffices to show that $G$ has a non-pivotal vertex $v$ in $V(G) - \{x,y\}$ unless~\ref{item:rpr2i} or~\ref{item:rpr2ii} holds.
  By Theorem~\ref{thm:nonpiv3}, we may assume that $G$ is pivot-equivalent to a graph $G'$ in $\Theta$.
  By~\ref{item:rpr2i}, we may assume that $G'$ is a graph consisting of at least three internally-disjoint paths between two vertices $a$ and $b$ having no common neighbor.
  Note that $\theta(1,3,\ell)$ with $\ell \ge 3$ is pivot-equivalent to a cycle of length $\ell+3$.
  By~\ref{item:rpr2ii}, $\{a,b\} \ne \{x,y\}$ and we take $v \in \{a,b\} - \{x,y\}$.
  Then by Lemma~\ref{lem:tp2}, $v$~is non-pivotal in $G'$ and so it is non-pivotal in $G$.
\end{proof}

Now let us show Corollary~\ref{cor:rooted prime reduction0}.

\begin{lemma}\label{lem:cycle only one deletion}
  Let $G$ be a graph locally equivalent to the cycle graph $C_n = (\{v_1,\dots,v_n\}, \{v_iv_{i+1} : 1\le i\le n\})$, where $n\ge 6$ and $v_{n+1} := v_1$.
  If $G / v_i$ is prime for every $1\le i\le n-1$,
  then $n$ is odd and $G = H_n$.
\end{lemma}
\begin{proof}
  Let $S$ be an isotropic system associated with a graphic presentation $(C_n, \mb{a}, \mb{b})$, where $\mb{a}$ and $\mb{b}$ are supplementary vectors in $K^{V(C_n)}$ such that $\mb{a}(v_i) = \alpha$ and $\mb{b}(v_i) = \beta$ for all $i$.
  Since $G$ is locally equivalent to $C_n$, by Proposition~\ref{prop:ig2}, there is a pair of supplementary vectors $\mb{c}$ and $\mb{d}$ in $K^{V(C_n)}$ such that $(G,\mb{c},\mb{d})$ is a graphic presentation of $S$.
  As $C_n * v_i \setminus v_i$ is prime for each $1\le i\le n$, by Proposition~\ref{prop:igmin} and Corollary~\ref{cor:igconn2}, $S|^{v_i}_{\gamma}$ is $3$-connected.
  Note that neither $C_n \setminus v_i$ nor $C_n / v_i$ is prime, implying that neither $S|^{v_i}_\alpha$ nor $S|^{v_i}_\beta$ is $3$-connected.
  From the assumption that $G / v_i$ is prime for all $1\le i\le n-1$, again by Proposition~\ref{prop:igmin} and Corollary~\ref{cor:igconn2}, we deduce that 
  $\mb{d}(v_i) = \gamma$ for each $1 \le i \le n-1$.

  Let $\{\mbt_i: 1\le i\le n\}$ and $\{\mb{u}_i: 1\le i\le n\}$ be the fundamental bases of $L$ with respect to Eulerian vectors $\mb{a}$ and $\mb{c}$, respectively, such that $\mbt_i(v_i) = \mb{b}(v_i) = \beta$ and $\mb{u}_i(v_i) = \mb{d}(v_i)$ for all $i$.
  Then for each $i \in \{1,\ldots,n\}$, we have $\mbt_i(v_{i-1}) = \mb{a}(v_{i-1}) = \alpha$, $\mbt_i(v_{i+1}) = \mb{a}(v_{i+1}) = \alpha$, and $\supp(\mbt_i) = \{v_{i-1},v_i,v_{i+1}\}$ where $v_{0} := v_n$.
  For all $i,j \in \{2,\ldots,n\}$, we have $\mb{u}_j(v_i) = \mb{d}(v_i) = \gamma$ if $i=j$ and $\mb{u}_j(v_i) = \mb{c}(v_i) \neq \gamma$ otherwise.

  \medskip\noindent 
  \textbf{Claim I.}
  For each $1 \le i\le n-1$, $\mb{u}_i(v_n) \in \{\beta,\gamma\}$ and one of the following holds:
  \begin{itemize}
    \item $\mb{u}_i(v_{i+1}) \in \{0,\alpha\}$ and $\mb{u}_i(v_j) = \beta$ for all $1 \le j \le i-1$.
    \item $\mb{u}_i(v_{i-1}) \in \{0,\alpha\}$ and $\mb{u}_i(v_j) = \beta$ for all $i+1 \le j \le n-1$.
  \end{itemize} 

  Let us fix $i \in \{1,\ldots,n-1\}$.
  As $0 = \inn{\mb{u}_i}{\mbt_i} = \innk{\mb{u}_i(v_{i-1})}{\alpha} + \innk{\mb{u}_i(v_{i})}{\beta} + \innk{\mb{u}_i(v_{i+1})}{\alpha} = 1+ \innk{\mb{u}_i(v_{i-1})}{\alpha} + \innk{\mb{u}_i(v_{i+1})}{\alpha}$, we have that for some $s \in \{1,-1\}$, $\mb{u}_{i}(v_{i-s}) \in \{\beta,\gamma\}$ and $\mb{u}_{i}(v_{i+s}) \in \{0,\alpha\}$.
  By reversing the labels of vertices if necessary, we may assume that $s=-1$.
  Recall that $\mb{u}_i(v_j) \neq \mb{d}(v_j) = \gamma$ for all $j \in \{1,\ldots,n\} - \{i,n\}$.
  Hence $\mb{u}_i(v_{i+1}) = \beta$ unless $i+1 = n$.

  Now let us show that
  \[
    \mb{u}_i(v_j) =
    \begin{cases}
      \beta   
      & \text{if $i+1 \le j \le n-1$}, \\
      \beta \text{ or } \gamma
      & \text{if $j=n$}.
    \end{cases}
  \]
  We proceed by induction on $j$.
  We already show this for $j = i+1$ and
  therefore we may assume that $j \ge i+2$.
  By the induction hypothesis, $\mb{u}_i(v_{j-1}) = \beta$.
  Again by the induction hypothesis, $\mb{u}_i(v_{j-2}) = \beta$ if $j-2 > i$.
  Note that $\mb{u}_i(v_{i}) = \gamma$ and therefore $\innk{\mb{u}_i(v_{j-2})}{\alpha} = 1$.
  Hence $0=\inn{\mb{u}_i}{\mbt_{j-1}}= \innk{\mb{u}_i(v_{j-2})}{\alpha} + \innk{\mb{u}_i(v_{j-1})}{\beta} + \innk{\mb{u}_i(v_{j})}{\alpha} = 1+ \innk{\mb{u}_i(v_{j})}{\alpha}$ and so we deduce that $\mb{u}_i(v_{j}) \in \{\beta,\gamma\}$.
  If $j < n$, then $\mb{u}_i(v_{j}) = \beta$ because $\mb{u}_i(v_j) \ne \mb{d}(v_j) = \gamma$.
  Therefore Claim~I is proved.

  \medskip\noindent 
  \textbf{Claim II.}
  For all $1 \le i \le n-1$, $\mb{c}(v_i) = \beta$.

  Let us fix $3 \le i \le n-3$.
  By Claim~I, either $\mb{u}_i(v_{i+1}) = \mb{u}_i(v_{i+2}) = \cdots = \mb{u}_i(v_{n-1}) = \beta$ or $\mb{u}_i(v_{i-1}) = \mb{u}_i(v_{i-2}) = \cdots = \mb{u}_i(v_{1}) = \beta$.
  By reversing the labels of vertices if necessary, we may assume that $\mb{u}_i(v_{i+1}) = \mb{u}_i(v_{i+2}) = \cdots = \mb{u}_i(v_{n-1}) = \beta$.
  So $\mb{c}(v_{i+1}) = \mb{c}(v_{i+2}) = \cdots = \mb{c}(v_{n-1}) = \beta$.
  Then by Claim~I, we may assume that $\mb{u}_{i+1}(v_{i+2}) = \beta$, because otherwise $\mb{c}(v_j) = \mb{u}_{i+1}(v_j) = \beta$ for all $1 \le j \le i$.
  As $\mb{u}_{i+1}(v_{i+2}) = \beta \ne 0$, 
  by Proposition~\ref{prop:basis},
  $\mb{u}_{i+2}(v_{i+1}) \ne 0$.
  Then $\mb{u}_{i+2}(v_{i+1}) = \mb{c}(v_{i+1}) = \beta$.
  By Claim~I, for all $1 \le j \le i+1$, $\mb{u}_{i+1}(v_j) = \beta$ and thus $\mb{c}(v_j) = \beta$.
  Hence Claim~II is proved.

  \medskip
  
  Let $\mbt := \sum_{i=1}^n \mbt_i$.
  Then $\mbt(v_i) = \beta$ for all $i$.
  By Claim~II, $\mb{u}_n(v_j) \in \{0,\mb{c}(v_j)\} = \{0,\beta\}$ for all $1\le j\le n-1$.
  Then as $0 = \inn{\mb{u}_n}{\mbt} = \innk{\mb{u}_n(v_n)}{\beta}$, we deduce that $\mb{u}_n(v_n) \in \{0,\beta\}$.
  Then $\mb{u}_n(v_n) = \beta$ 
  because $\mb{u}_n(v_n) = \mb{d}(v_n) \ne 0$.
  Thus,
  \[
    \mb{d}(v_i) =
    \begin{cases}
      \beta  & \text{if $i=n$}, \\
      \gamma   & \text{otherwise}.
    \end{cases}
  \]
  Then $\mb{c}(v_n)$ is $\alpha$ or $\gamma$ because $\mb{c}(v_n) \ne \mb{d}(v_n) = \beta$.
  By Claim~I, $\mb{c}(v_n) = \mb{u}_1(v_n) \in \{\beta,\gamma\}$ and therefore by Claim~II,
  \[
    \mb{c}(v_i) =
    \begin{cases}
      \gamma  & \text{if $i=n$}, \\
      \beta   & \text{otherwise}.
    \end{cases}
  \]
  Hence by Claim~I, $\mb{u}_{i}(v_n) = \mb{c}(v_n) = \gamma$ for every $1\le i \le n-1$.

  \medskip\noindent 
  \textbf{Claim III.}
  If $1\le i\le n-2$ and $\mb{u}_{i}(v_{i+1}) = 0$, then for each $j$ with $i+1 \le j \le n-1$,
  \[
    \mb{u}_i(v_j) =
    \begin{cases}
      0     & \text{if $j \not\equiv i \pmod{2}$}, \\
      \beta & \text{otherwise}.
    \end{cases}
  \]

  We may assume that $i+2 \le n-1$.
  By Claim~II, $\mb{u}_i(v_{i+2}) \in \{0, \beta\}$.
  Since $\mb{u}_i(v_{i}) = \gamma$ and $\mb{u}_i(v_{i+1}) = 0$, we have $0 = \inn{\mb{u}_i}{\mbt_{i+1}} = \innk{\mb{u}_i(v_{i})}{\alpha} + \innk{\mb{u}_i(v_{i+1})}{\beta} + \innk{\mb{u}_i(v_{i+2})}{\alpha} =  1 + \innk{\mb{u}_i(v_{i+2})}{\alpha}$ and so $\mb{u}_i(v_{i+2}) = \beta$.
  Therefore we may assume that $j \ge i+3$.
  By Claim~II, $\mb{u}_i(v_{j-2})$, $\mb{u}_i(v_{j-1})$, and $\mb{u}_i(v_{j})$ are in $\{0,\beta\}$.
  Then $0 = \inn{\mb{u}_i}{\mbt_{j-1}} = \innk{\mb{u}_i(v_{j-2})}{\alpha} + \innk{\mb{u}_i(v_{j})}{\alpha}$ and so
  $\mb{u}_i(v_{j-2}) = \mb{u}_i(v_{j})$.
  This proves Claim~III.

  \medskip

  By Claim~I, 
  $\mb{u}_1(v_{n}) = \mb{c}(v_n) = \gamma$ and $\mb{u}_1(v_2) = 0$.
  Then $0 = \inn{\mb{u}_1}{\mbt_n} = \innk{\mb{u}_1(v_1)}{\alpha} + \innk{\mb{u}_1(v_n)}{\beta} + \innk{\mb{u}_1(v_{n-1})}{\alpha} = \innk{\mb{u}_1(v_{n-1})}{\alpha}$.
  Hence $\mb{u}_1(v_{n-1}) = 0$ because $\mb{c}(v_{n-1}) = \beta$.
  Then by Claim~III, $n$ is odd.
  
  By Claims~I and~III,
  \[
    \mb{u}_1(v_{j}) =
    \begin{cases}
      0         & \text{if $j \in \{2,4,6,\ldots,n-1\}$}, \\
      \gamma    & \text{if $j\in \{1,n\}$}, \\
      \beta     & \text{otherwise}.
    \end{cases}
  \]
  Then
  by Proposition~\ref{prop:basis}, $\mb{u}_j(v_1) = 0$ for $j \in \{2,4,6,\ldots,n-1\}$.
  Then again by Claims~I and~III, for each $i \in \{2,4,6,\ldots,n-1\}$,
  \[
    \mb{u}_i(v_j) =
    \begin{cases}
      0       & \text{if $j \in \{1,3,5,\ldots,i-1\}$}, \\
      \gamma  & \text{if $j\in\{i,n\}$}, \\
      \beta   & \text{otherwise}. \\
    \end{cases}
  \]
  By symmetry, we also conclude that for each $i\in \{1,3,5,\ldots,n-2\}$,
  \[
    \mb{u}_i(v_j) =
    \begin{cases}
      0       & \text{if $j \in \{n-1,n-3,n-5,\ldots,i+1\}$}, \\
      \gamma  & \text{if $j\in\{i,n\}$}, \\
      \beta   & \text{otherwise}.
    \end{cases}
  \]
  Note that for $1\le i < j \le n$, two vertices $v_i$ and $v_j$ are adjacent in $G$ if and only if $\mb{u}_i(v_j) \ne 0$.
  Therefore, $G = H_n$.
\end{proof}

\begin{proof}[Proof of Corollary~\ref{cor:rooted prime reduction0}]
  If $G$ is not locally equivalent to a cycle, then by Theorem~\ref{thm:noness2}, $G$ has a non-essential vertex $v \in V(G) - \{x\}$ and so $G\setminus v$ or $G * v \setminus v$ is prime.
  Therefore we may assume that $G$ is locally equivalent to a cycle $C$.
  Suppose that for all $v \in V(G) - \{x\}$, neither $G \setminus v$ nor $G * v \setminus v$ is prime.
  Because $C * v \setminus v$ is prime, $G / v$ is prime.
  By Lemma~\ref{lem:cycle only one deletion}, $|V(G)|$ is odd, $G$ is isomorphic to $H_{|V(G)|}$, and $x$ is adjacent to all other vertices, a contradiction.
\end{proof}

\section*{Acknowledgments}
The authors would like to thank anonymous reviewers for their valuable suggestions.


\begin{thebibliography}{10}

  \bibitem{Allys1994iso}
  Lo{\"{i}}c Allys.
  \newblock Minimally {$3$}-connected isotropic systems.
  \newblock {\em Combinatorica}, 14(3):247--262, 1994.
  
  \bibitem{Bouchet1987iso}
  Andr\'{e} Bouchet.
  \newblock Isotropic systems.
  \newblock {\em European J. Combin.}, 8(3):231--244, 1987.
  
  \bibitem{Bouchet1987cir}
  Andr\'{e} Bouchet.
  \newblock Reducing prime graphs and recognizing circle graphs.
  \newblock {\em Combinatorica}, 7(3):243--254, 1987.
  
  \bibitem{Bouchet1988iso}
  Andr\'{e} Bouchet.
  \newblock Graphic presentations of isotropic systems.
  \newblock {\em J. Combin. Theory Ser. B}, 45(1):58--76, 1988.
  
  \bibitem{Bouchet1989iso}
  Andr\'{e} Bouchet.
  \newblock Connectivity of isotropic systems.
  \newblock In {\em Combinatorial {M}athematics: {P}roceedings of the {T}hird {I}nternational {C}onference ({N}ew {Y}ork, 1985)}, volume 555 of {\em Ann. New York Acad. Sci.}, pages 81--93. New York Acad. Sci., New York, 1989.
  
  \bibitem{Bouchet1990}
  Andr\'{e} Bouchet.
  \newblock {$\kappa$}-transformations, local complementations and switching.
  \newblock In {\em Cycles and rays ({M}ontreal, {PQ}, 1987)}, volume 301 of {\em NATO Adv. Sci. Inst. Ser. C: Math. Phys. Sci.}, pages 41--50. Kluwer Acad. Publ., Dordrecht, 1990.
  
  \bibitem{Bouchet1991algo}
  Andr\'{e} Bouchet.
  \newblock An efficient algorithm to recognize locally equivalent graphs.
  \newblock {\em Combinatorica}, 11(4):315--329, 1991.
  
  \bibitem{Bouchet1994circle}
  Andr\'{e} Bouchet.
  \newblock Circle graph obstructions.
  \newblock {\em J. Combin. Theory Ser. B}, 60(1):107--144, 1994.
  
  \bibitem{FonDerFlaass1988}
  Dmitrii~G. Fon-Der-Flaass.
  \newblock On local complementations of graphs.
  \newblock In {\em Combinatorics ({E}ger, 1987)}, volume~52 of {\em Colloq. Math. Soc. J\'{a}nos Bolyai}, pages 257--266. North-Holland, Amsterdam, 1988.
  
  \bibitem{FonDerFlaass1989}
  Dmitrii~G. Fon-Der-Flaass.
  \newblock Distance between locally equivalent graphs.
  \newblock {\em Metody Diskret. Analiz.}, (48):85--94, 106--107, 1989.
  
  \bibitem{FonDerFlaass1996}
  Dmitrii~G. Fon-Der-Flaass.
  \newblock {\em Local Complementations of Simple and Directed Graphs}, pages 15--34.
  \newblock Springer Netherlands, Dordrecht, 1996.
  
  \bibitem{Geelen1996thesis}
  Jim Geelen.
  \newblock {\em Matchings, matroids and unimodular matrices}.
  \newblock ProQuest LLC, Ann Arbor, MI, 1996.
  \newblock Thesis (Ph.D.)--University of Waterloo (Canada).
  
  \bibitem{Geelen2009circle}
  Jim Geelen and Sang-il Oum.
  \newblock Circle graph obstructions under pivoting.
  \newblock {\em J. Graph Theory}, 61(1):1--11, 2009.
  
  \bibitem{Godsil2001}
  Chris Godsil and Gordon Royle.
  \newblock {\em Algebraic graph theory}, volume 207 of {\em Graduate Texts in Mathematics}.
  \newblock Springer-Verlag, New York, 2001.
  
  \bibitem{Lang1987}
  Serge Lang.
  \newblock {\em Linear algebra}.
  \newblock Undergraduate Texts in Mathematics. Springer-Verlag, New York, third edition, 1987.
  
  \bibitem{Oum2005thesis}
  Sang-il Oum.
  \newblock {\em Graphs of bounded rank-width}.
  \newblock PhD thesis, Princeton University, 2005.
  
  \bibitem{Oum2005}
  Sang-il Oum.
  \newblock Rank-width and vertex-minors.
  \newblock {\em J. Combin. Theory Ser. B}, 95(1):79--100, 2005.
  
  \bibitem{Oum2008wqo}
  Sang-il Oum.
  \newblock Rank-width and well-quasi-ordering.
  \newblock {\em SIAM J. Discrete Math.}, 22(2):666--682, 2008.
  
  \bibitem{Oxley2000b}
  James Oxley and Haidong Wu.
  \newblock Matroids and graphs with few non-essential elements.
  \newblock {\em Graphs Combin.}, 16(2):199--229, 2000.
  
  \bibitem{Oxley2000a}
  James Oxley and Haidong Wu.
  \newblock On the structure of 3-connected matroids and graphs.
  \newblock {\em European J. Combin.}, 21(5):667--688, 2000.
  
  \bibitem{Oxley2004}
  James Oxley and Haidong Wu.
  \newblock The 3-connected graphs with exactly three non-essential edges.
  \newblock {\em Graphs Combin.}, 20(2):233--246, 2004.
  
  \bibitem{Seymour1988}
  Paul Seymour.
  \newblock On the connectivity function of a matroid.
  \newblock {\em J. Combin. Theory Ser. B}, 45(1):25--30, 1988.
  
  \bibitem{Tutte1961}
  William~T. Tutte.
  \newblock A theory of {$3$}-connected graphs.
  \newblock {\em Nederl. Akad. Wetensch. Proc. Ser. A 64 = Indag. Math.}, 23:441--455, 1961.
  
  \bibitem{Tutte1966cm}
  William~T. Tutte.
  \newblock Connectivity in matroids.
  \newblock {\em Canadian J. Math.}, 18:1301--1324, 1966.
  
  \end{thebibliography}

%

\appendix

\section{Lemma~\ref{lem:oddcyc}: An alternative proof}\label{app:oddcyc}

To prove Lemma~\ref{lem:oddcyc} without using isotropic systems and Theorem~\ref{thm:Hamiltonian}, we first review double occurrence words, chord diagrams, and circle graphs based on Bouchet~\cite{Bouchet1994circle}.

A \emph{double occurrence word} is a sequence of letters such that each letter appears exactly twice.
The \emph{circle graph} $A(m)$ of a double occurrence word $m$ is a graph such that (i) its vertex set is the set of letters appearing in $m$, and (ii) two vertices $v$ and $w$ are adjacent if and only if they appear alternatively in $m$, that is, $m = (\cdots v \cdots w \cdots v \cdots w \cdots)$ or $(\cdots w \cdots v \cdots w \cdots v \cdots)$.
For instance, $A(abacbc)$ is a path with $3$ vertices.

\begin{figure}
  \centering
  \begin{tikzpicture}
    \foreach \x [count=\p] in {0,...,5} {
      \node[shape=circle,fill=black,scale=0.5] (\p) at (\x*180/3:2) {};
    };
    \draw (0:2.33) node {$a$};
    \draw (60:2.33) node {$b$};
    \draw (120:2.33) node {$a$};
    \draw (180:2.33) node {$c$};
    \draw (240:2.33) node {$c$};
    \draw (300:2.33) node {$b$};
    \draw[thick] (1) arc (0:360:2);
    \draw (1) -- (3);
    \draw (2) -- (6);
    \draw (4) -- (5);
    \node at (0,-2.8) {$D(m)$};
  \end{tikzpicture}
  \hspace{1cm}
  \begin{tikzpicture}
    \foreach \x [count=\p] in {0,...,2} {
      \node[shape=circle,fill=black,scale=0.5] (\p) at (\x*120-30:1.5) {};
    };
    \draw (-30:1.8) node {$b$};
    \draw (90:1.8) node {$a$};
    \draw (210:1.2) node {$c$};
    \draw (1.center) to [bend right=30] (2.center);
    \draw (1.center) to (2.center);
    \draw (1.center) to [bend left=30] (2.center);
    \draw (2.center) to [bend right=20] (3.center);
    \draw (1.center) to [bend left=20] (3.center);
    \draw (3) to [loop,distance=35,in=170,out=250] (3);
    \node at (0,-2.4) {$T(m)$};
  \end{tikzpicture}
  \caption{Drawings of $D(m)$ and $T(m)$, where $m = (ababcc)$.}
\label{fig:wdwtw}
\end{figure}
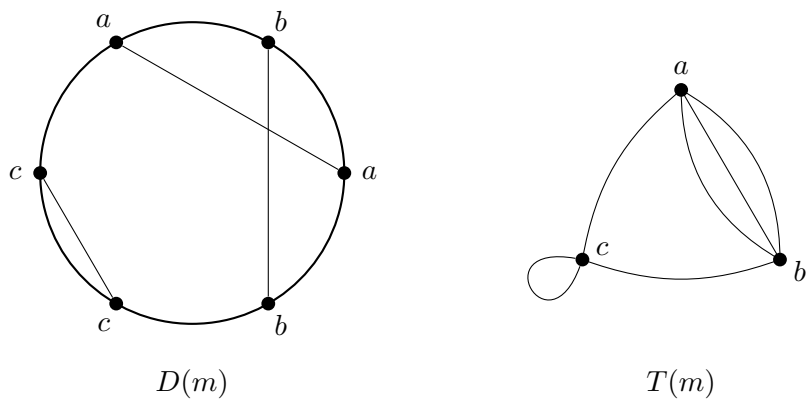

Given a double occurrence word $m$, let us write letters of the word along a circle.
Then by connecting the same letters by a chord, we obtain a chord diagram $D(m)$; see Figure~\ref{fig:wdwtw}.
We can regard every chord diagram as a cubic multigraph whose vertices are ends of chords, and edges are chords and segments of the circle cut by ends of chords.
By contracting all chords of $D(m)$, we obtain a $4$-regular multigraph $T(m)$.

For a double occurrence word $m$ and a letter $v$ in $m$, let $m * v$ be the double occurrence word obtained from $m$ by reversing the sub-word between two $v$'s.
For instance, $(abcdabcd)*b = (abadcbcd)$.
One may observe that $A(m*v) = A(m)*v$ and $T(m*v) = T(m)$.

The \emph{square} $G^2$ of a graph $G$ is a graph on $V(G)$ such that two vertices are adjacent in $G^2$ if and only if their distance is at most two in $G$.

\begin{lemma}[see Godsil and Royle~{\cite[Lemma~17.5.1]{Godsil2001}}]
  \label{lem:planar}
  Let $m$ be a double occurrence word.
  Then $D(m)$ is a planar graph if and only if $A(m)$ is bipartite.
\end{lemma}

\begin{proof}[Proof of Lemma~\ref{lem:oddcyc}]
  Suppose that a bipartite graph $G$ is locally equivalent to an odd cycle $C_k$ with $k\geq 5$.
  The odd cycle $C_k$ is a circle graph, and let $m_k$ be a double occurrence word such that $A(m_k) = C_k$; see Figure~\ref{fig:C7}.
  Since $G$ and $C_k$ are locally equivalent, there is a sequence $a_1, a_2, \dots, a_t$ of vertices such that $G = C_k * a_1 * \dots * a_t = A(m_k) * a_1 * \dots * a_t = A(m_k * a_1 * \dots * a_t)$.
  Let $m := m_k * a_1 * \dots * a_t$.
  By Lemma~\ref{lem:planar}, $D(m)$ is planar.
  Then $T(m)$ is planar because $T(m)$ is a minor of $D(m)$.

  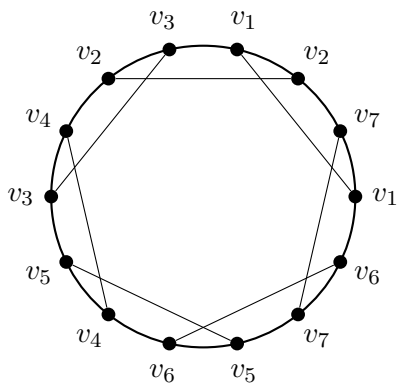
\begin{figure}
    \centering
    \begin{tikzpicture}
      \foreach \x [count=\p] in {0,...,13} {
        \node[shape=circle,fill=black, scale=0.5] (\p) at (\x*180/7:2) {};
      };
      \foreach \x [count=\p] in {0,...,6} {
        \draw (\x*360/7:2.4) node {$v_\p$};
        \draw (3*180/7+\x*360/7:2.4) node {$v_\p$};
      }; 
      \draw[thick] (1) arc (0:360:2);
      \draw (1) -- (4);
      \draw (3) -- (6);
      \draw (5) -- (8);
      \draw (7) -- (10);
      \draw (9) -- (12);
      \draw (11) -- (14);
      \draw (13) -- (2);
    \end{tikzpicture}
    \caption{A drawing of $D(m_7)$.}
  \label{fig:C7}
  \end{figure}

  It is easy to see that $T(m_k)$ is isomorphic to $C_k^2$.
  Note that $C_5^2$ is isomorphic to the complete graph with $5$ vertices, and $C_{\ell+2}^2$ with $\ell\geq 3$ has a minor isomorphic to $C_{\ell}^2$.
  Therefore, $T(m_k)$ is not planar, and this contradicts that $T(m_k) = T(m)$ is planar.
\end{proof}

\section{Binary matroids with few non-essential elements}\label{app:binary matroid}
In this section, we will explain why 
our theorems restricted to bipartite graphs are equivalent to 
known theorems on binary matroids due to Oxley and Wu.
We use the same terminology as used in Oxley and Wu~\cite{Oxley2000a,Oxley2000b}.
An element $e$ of a $3$-connected matroid $M$ is \emph{deletable} if $M\setminus e$ is $3$-connected, and it is \emph{contractible} if $M/e$ is $3$-connected.
We say $e$ is \emph{non-essential} if it is deletable or contractible.

\subsection{Fewer than two non-essential elements}

We first review the result of Oxley and Wu~\cite{Oxley2000a} extending Tutte's wheel and whirl theorem~{\cite[8.3]{Tutte1966cm}}.

\begin{theorem}[Oxley and Wu~{\cite[Corollary~3.5]{Oxley2000a}}]
  Let $M$ be a $3$-connected matroid with at least four elements.
  Then $M$ has fewer than two non-essential elements
  if and only if $M$ is a wheel or a whirl.
\end{theorem}

Since whirls are non-binary, we obtain the following consequence for binary matroids.

\begin{corollary}\label{cor:binmat1}
  Let $M$ be a $3$-connected binary matroid with at least four elements.
  Then $M$ has fewer than two non-essential elements
  if and only if $M$ is a wheel.
\end{corollary}

We will show the equivalence of Corollary~\ref{cor:binmat1} and Corollary~\ref{cor:bippiv}.
For a base $B$ of a matroid $M$, the \emph{fundamental graph of $M$ with respect to $B$} is a graph on $E(M)$ such that 
two vertices $x$ and $y$ are adjacent if and only if 
$B\triangle \{x,y\}$ is a base of $M$.
It is well known that fundamental graphs determine connected binary matroids up to duality.

The \emph{connectivity function} of a matroid $M$ with the rank function~$r$ is defined as a function, denoted by $\lambda_M$, on subsets of $E(M)$ 
such that $\lambda_M(X)=r(X)+r(E(M)-X)-r(E(M))$ for each $X\subseteq E(M)$.
Oum~\cite{Oum2005} observed that if $G$ is a fundamental graph of a binary matroid~$M$, then 
$\lambda_M(X)$ is equal to the cut-rank function of $G$
and therefore $M$ is $3$-connected if and only if $G$ is connected and prime. Since all prime graphs with at least four vertices are connected, we deduce the following easily. 

\begin{lemma}\label{lem:fundamental graphs of binary matroids 2}
  Let $M$ be a binary matroid with at least four elements
  and let $G$ be its fundamental graph.
  Then $M$ is $3$-connected if and only if 
  $G$ is prime.
  \qed
\end{lemma}
Minors of binary matroids correspond to pivot-minors of their fundamental graphs as follows.
\begin{lemma}[Oum~\cite{Oum2005}]\label{lem:fundamental graphs of binary matroids}
  Let $G$ be the fundamental graph of a binary matroid $M$ with respect to a base~$B$.
  Then the following hold.
  \begin{enumerate}[label=\rm(\roman*)]
    \item\label{item:fgbm2} For an edge $uv$ of $G$, $G\wedge uv$ is the fundamental graph of $M$ with respect to $B\triangle \{u,v\}$.
    \item\label{item:fgbm3} For $v\in E(M)-B$, $G \setminus v$ is the fundamental graph of $M\setminus v$ with respect to $B$.
    \item\label{item:fgbm4} For $v\in B$, $G\setminus v$ is the fundamental graph of $M/v$ with respect to $B-\{v\}$.
  \end{enumerate}
\end{lemma}

It is easy to see that $U_{2,4}$ is the unique $3$-connected matroid on four elements
and therefore there is no binary $3$-connected matroid on four elements. So we deduce the following by using the previous two lemmas.
\begin{lemma}\label{lem:fundamental graphs of binary matroids 3}
  Let $M$ be a $3$-connected binary matroid with at least four elements and let $G$ be its fundamental graph.
  Then 
  $e\in E(M)$ is non-essential in~$M$ if and only if it is non-pivotal in~$G$.
  \qed
\end{lemma}

Now we are ready to show how Corollary~\ref{cor:bippiv} implies Corollary~\ref{cor:binmat1}.

\begin{proof}[Proof of Corollary~\ref{cor:binmat1} using Corollary~\ref{cor:bippiv}]
  Suppose that $M$ has fewer than two non-essential elements.
  Let $G$ be a fundamental graph of $M$.
  Then $V(G) = E(M)$ and 
  by Lemmas~\ref{lem:fundamental graphs of binary matroids 2} and \ref{lem:fundamental graphs of binary matroids 3}, $G$ is prime and has fewer than two non-pivotal vertices.
  By Corollary~\ref{cor:bippiv}, $G$ is pivot-equivalent to an even cycle $C$.
  By Lemma~\ref{lem:fundamental graphs of binary matroids}\ref{item:fgbm2}, $C$ is a fundamental graph of $M$
  and so $M$ is a wheel.

  Conversely, a wheel with at least six elements has no non-essential elements. This completes the proof. 
\end{proof}

Finally, we prove that Corollary~\ref{cor:binmat1} implies Corollary~\ref{cor:bippiv}.

\begin{proof}[Proof of Corollary~\ref{cor:bippiv} using Corollary~\ref{cor:binmat1}]
  Suppose that $G$ is a prime bipartite graph
  with at least four vertices 
  having fewer than two non-pivotal vertices.
  Let $M$ be a binary matroid having $G$ as a fundamental graph.
  Then by Lemmas~\ref{lem:fundamental graphs of binary matroids 2} and \ref{lem:fundamental graphs of binary matroids 3}, $M$ is $3$-connected and has fewer than two non-essential elements.
  By Corollary~\ref{cor:binmat1}, $M$ is a wheel and by Lemma~\ref{lem:fundamental graphs of binary matroids}\ref{item:fgbm2}, $G$ is pivot-equivalent to an even cycle.
\end{proof}

\subsection{Exactly two non-essential elements}

We now discuss the relation between Corollary~\ref{cor:bippiv3} and the results of Oxley and Wu~\cite{Oxley2000b}.
Before stating theorems in~\cite{Oxley2000b} and their corollary restricted to binary matroids,
we quickly review some graphs and binary matroids defined in~\cite{Oxley2000b}.

\begin{figure}
  \centering
  \begin{tikzpicture}
    \begin{scope}[xshift=-9cm]
      \node[shape=circle,fill=black, scale=0.5] (v1) at (0,0) {};
      \node[shape=circle,fill=black, scale=0.5] (v2) at (1.5,0) {};
      \node[shape=circle,fill=black, scale=0.5] (v3) at (3.0,0) {};
      \node[shape=circle,fill=black, scale=0.5] (v4) at (4.5,0) {};

      \node at (-0.2,-0.2) {$y_1$};
      \node at (1.25,+0.2) {$x_2$};
      \node at (3.27,-0.2) {$x_1$};
      \node at (4.77,-0.2) {$y_2$};

      \node[shape=circle,fill=black, scale=0.5] (e1) at (2.08,0.733) {};
      \node[shape=circle,fill=black, scale=0.5] (e2) at (3.0,0.94) {};
      \node[shape=circle,fill=black, scale=0.5] (e3) at (3.92,0.733) {};

      \node[shape=circle,fill=black, scale=0.5] (f1) at (0.72,-0.805) {};
      \node[shape=circle,fill=black, scale=0.5] (f2) at (2.28,-0.805) {};
  
      \draw[very thick] (v1) -- (v2) -- (v3) -- (v4);
      \draw (v1) to [bend right = 80] (v3);
      \draw (v2) to [bend left = 80] (v4);
      \draw (v1) to [bend left = 85] (v4);

      \draw[very thick] (v3) -- (e1);
      \draw[very thick] (v3) -- (e2);
      \draw[very thick] (v3) -- (e3);

      \draw[very thick] (v2) -- (f1);
      \draw[very thick] (v2) -- (f2);

    \end{scope}
    
    \begin{scope}[scale=1.3]
      \node[shape=circle,fill=black,scale=0.5] (a1) at (90+180/5:1.2) {};
      \node[shape=circle,fill=black,scale=0.5] (a2) at (90+2*180/5:1.2) {};
      \node[shape=circle,fill=black,scale=0.5] (a3) at (90+3*180/5:1.2) {};
      \node[shape=circle,fill=black,scale=0.5] (a4) at (90+4*180/5:1.2) {};
      
      \node (t1) at (1.2*0.7071,1.2*0.7071) {};
      \node (t2) at (1.2,0) {};
      \node (t3) at (1.2*0.7071,-1.2*0.7071) {};
      
      \node[shape=circle,fill=black,scale=0.5] (b1) at (1.2*0.7071 + 0.2, 1.2*0.7071 + 0.1) {};
      \node[shape=circle,fill=black,scale=0.5] (b2) at (1.2 + 0.2, 0 + 0.2) {};
      \node[shape=circle,fill=black,scale=0.5] (b3) at (1.2*0.7071 + 0.2, -1.2*0.7071 + 0.1) {};
      
      \node[shape=circle,fill=black,scale=0.5] (c1) at (1.2*0.7071 - 0.2, 1.2*0.7071 - 0.1) {};
      \node[shape=circle,fill=black,scale=0.5] (c2) at (1.2 - 0.3, 0 - 0.3) {};
      \node[shape=circle,fill=black,scale=0.5] (c3) at (1.2*0.7071 - 0.3, -1.2*0.7071 - 0.25) {};
      
      \node[shape=circle,fill=black,scale=0.5] (c) at (0,0) {};
      \node[shape=circle,fill=black,scale=0.5] (x) at (0,1.2) {};
      \node[shape=circle,fill=black,scale=0.5] (y) at (0,-1.2) {};

      \draw (-0.098,-0.30) node {$c$};

      \draw[very thick] (x) -- (c) -- (y);
      
      \draw (x) -- (a1) -- (a2) -- (a3) -- (a4) -- (y);
      \draw[very thick] (c) -- (a1);
      \draw[very thick] (c) -- (a2);
      \draw[very thick] (c) -- (a3);
      \draw[very thick] (c) -- (a4);
      
      \draw (x) -- (b1) -- (b2) -- (b3) -- (y);
      \draw[very thick] (c) -- (b1);
      \draw[very thick] (c) -- (b2);
      \draw[very thick] (c) -- (b3);
      
      \draw (x) -- (c1) -- (c2) -- (c3) -- (y);
      \draw[very thick] (c) -- (c1);
      \draw[very thick] (c) -- (c2);
      \draw[very thick] (c) -- (c3);
    \end{scope}
  \end{tikzpicture}
  \caption{A twisted wheel (left) and a multidimensional wheel (right).}
  \label{fig: TMwheel}
\end{figure}
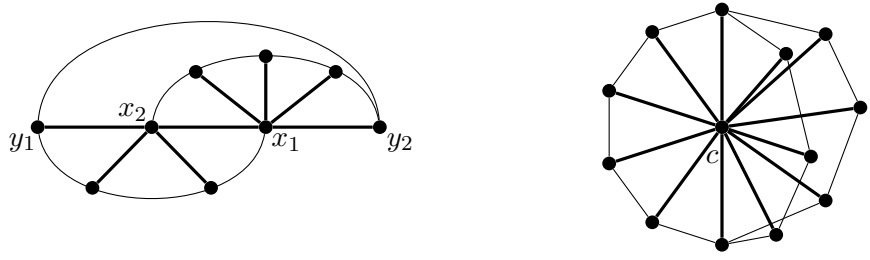

A \emph{twisted wheel} is a graph obtained from the disjoint union of two paths $P_1$ from $x_1$ to $y_1$ and $P_2$ from $x_2$ to $y_2$, each of length at least $2$, 
by adding edges from $x_i$ to every vertex of $P_{3-i}$ for all $i=1,2$
and adding an edge $y_1y_2$, see Figure~\ref{fig: TMwheel}(left).
Note that 
if $M$ is the cycle matroid of a twisted wheel,
then its fundamental graph 
with respect to the base marked by thick edges in Figure~\ref{fig: TMwheel}(left)
is isomorphic to $\theta(1,a,b)$ with odd integers $a,b\ge 5$.

A \emph{multidimensional wheel} is a graph obtained from $\theta(\ell_1,\ldots,\ell_k)$ with $k\ge 3$ and $\ell_1,\ldots,\ell_k \ge 2$ by adding a new vertex $c$ adjacent to all the other vertices, see Figure~\ref{fig: TMwheel}(right).
If $M$ is the cycle matroid of such a multidimensional wheel and $B$ is its base consisting of the edges incident with $c$, then its fundamental graph with respect to $B$ is isomorphic to $\theta(2\ell_1,\ldots,2\ell_k)$ and 
the two vertices of degree~$k$ belong to $B$.

For a binary matroids $M_1$ and a matroid $M_2$, if $T:=E(M_1) \cap E(M_2)$ is a triangle in both $M_1$ and $M_2$, 
then 
there is a unique matroid $P_T(M_1,M_2)$, called the \emph{generalized parallel connection of $M_1$ and $M_2$ across $T$}, on $E(M_1) \cup E(M_2)$ such that 
for every subset $X$ of $E(M_1) \cup E(M_2)$, 
it is a flat in $P_T(M_1,M_2)$ if and only if $X \cap E(M_1)$ is a flat in $M_1$
and $X\cap E(M_2)$ is a flat in $M_2$.
Observe that if $M_1$ and $M_2$ have binary representations
\[ 
  M_1: \quad \left(\begin{array}{@{}c|c@{}}
    \begin{matrix} 1 & 0 & 1 \\
      0 & 1 & 1
    \end{matrix} 
    & A_1 \\
    \hline 
    \mathbf{0}
    & B_1
  \end{array}\right), 
  \qquad\qquad
  M_2: \quad 
  \left(\begin{array}{@{}c|c@{}}
    \begin{matrix} 1 & 0 & 1 \\
      0 & 1 & 1
    \end{matrix} 
    & A_2 \\
    \hline 
    \mathbf{0}
    & B_2
  \end{array}\right)
\] 
where the first three columns are indexed by elements of $T$ in the same order
and $\mathbf{0}$ represents a zero matrix,
then $P_{T}(M_1,M_2)$ has the following binary representation
\[ 
  \left(\begin{array}{@{}c|c|c@{}}
    \begin{matrix}
      1&0&1\\
      0&1&1 
    \end{matrix} 
    & A_1 & A_2 \\    \hline 
    \mathbf{0}  & B_1 & \mathbf{0} \\\hline
    \mathbf{0}  & \mathbf{0} & B_2 
  \end{array}\right).
\] 

For a positive integer $k$, 
a binary matroid $\Lambda_k$ is the vector matroid associated with the following $(k+1) \times (2k+2)$ binary matrix
\[ 
  \begin{pmatrix}
    I_k & \mathbf{0} & \mathbf{1} - I_k & \mathbf{1} \\
    \mathbf{0} & 1 & \mathbf{1} & 1 \\
  \end{pmatrix},
\] 
where $\mathbf{1}$ is a matrix whose every entry is $1$.
If $k$ is odd, then let $b_1,\ldots,b_k,x,a_1,\ldots,a_k,y$ be the indices of the columns in order, and if $k$ is even, then let $b_1,\ldots,b_{k-1},a_k,x,a_1,\ldots,a_{k-1},b_k,y$ be the indices of the columns in order; see~{\cite[page~225]{Oxley2000b}}.
Then by rearranging the columns $b_1,\ldots,b_k,y,a_1,\ldots,a_k,x$ in order and taking elementary row operations, we obtain another binary representation of $\Lambda_k$ as follows:
\begin{align*}
  \begin{pmatrix}
    I_k & \mathbf{0} & I_k & \mathbf{1} \\
    \mathbf{0} & 1 & \mathbf{1} & \delta_k \\
  \end{pmatrix}
  \text{ where }
  \delta_k :=
  \begin{cases}
    1 & \text{if $k$ is odd}, \\
    0 & \text{otherwise}.
  \end{cases}
\end{align*}
Note that 
for each $i\in\{1,2,\ldots,k\}$, 
$T_i := \{y,a_i,b_i\}$ is a triangle of $\Lambda_k$.
We are going to define matroids $M_1(n_1,n_2,\ldots,n_k)$ and $M_2(n_1,n_2,\ldots,n_k)$ for positive integers $n_1,n_2,\ldots,n_k$.
Let $W_1,W_2,\ldots,W_k$ be matroids such that for each $i\in\{1,2,\ldots,k\}$,
$W_i\cong M(\mathcal{W}_{n_{i}+2})$, 
$E(W_i)\cap E(\Lambda_k)=T_i$, and $T_i$ is a triangle of $W_i$.  
For each $i\in\{0,1,2\ldots,k\}$, let 
\[ 
  N_i=\begin{cases}
    \Lambda_k & \text{if }i=0,\\
    P_{T_i}(W_i,N_{i-1})\setminus x_i &\text{otherwise}
  \end{cases}
\]
where $x_i$ is the unique element of $T_i$ in a rim of $W_i$. 
If $x_i=a_i$ for each $i=1,2,\ldots,k$, then 
we define $M_1(n_1,\ldots,n_k):=N_k$. 
If $x_i=a_i$ for each $i=1,2,\ldots,k-1$ and $x_k=b_k$, then
we define $M_2(n_1,\ldots,n_k):=N_k$.
Then $M_1(n_1,\ldots,n_k) $ has the following standard binary representation:
\begin{align*}
  \begin{pmatrix}
    I_{n_1} & \cdots & \mathbf{0} & \mathbf{0} & 
    A_{n_1} & \cdots & \mathbf{0} & \mathbf{1} \\
    \vdots & \ddots & \vdots & \vdots & 
    \vdots & \ddots & \vdots & \vdots \\
    \mathbf{0} & \cdots & I_{n_k} & \mathbf{0} & 
    \mathbf{0} & \cdots & A_{n_k} & \mathbf{1} \\
    \mathbf{0} & \cdots & \mathbf{0} & 1 & 
    R_{n_1} & \cdots & R_{n_k} & \delta_k \\
  \end{pmatrix}
\end{align*}
where $A_n$ is an $n\times n$ matrix such that all diagonal entries and $(i+1,i)$-entries with $1\le i\le n-1$ are one and the other entries are zero, and $R_n$ is a $1\times n$ matrix such that $(1,n)$-entry is one and the other entries are zero.
Similarly, $M_2(n_1,\ldots,n_k)$ has the following standard binary representation:
\begin{align*}
  \begin{pmatrix}
    I_{n_1} & \cdots & \mathbf{0} & \mathbf{0} & 
    A_{n_1} & \cdots & \mathbf{0} & \mathbf{1} \\
    \vdots & \ddots & \vdots & \vdots & 
    \vdots & \ddots & \vdots & \vdots \\
    \mathbf{0} & \cdots & I_{n_k} & \mathbf{0} & 
    \mathbf{0} & \cdots & A_{n_k} & \mathbf{1} \\
    \mathbf{0} & \cdots & \mathbf{0} & 1 & 
    R_{n_1} & \cdots & R_{n_k} & 1-\delta_k \\
  \end{pmatrix}.
\end{align*}
Then we deduce that $M_1(n_1,\ldots,n_k)$ has a fundamental graph isomorphic to $\theta(1,2n_{1}+1,\ldots,2n_{k}+1)$ if $k$ is odd and $\theta(2n_{1}+1,\ldots,2n_{k}+1)$ otherwise.
Similarly $M_2(n_1,\ldots,n_k)$ has a fundamental graph isomorphic to $\theta(2n_{1}+1,\ldots,2n_{k}+1)$ if $k$ is odd and $\theta(1,2n_{1}+1,\ldots,2n_{k}+1)$ otherwise.

Now we state two results of Oxley and Wu~\cite{Oxley2000b} and their corollary for binary matroids.

\begin{theorem}[Oxley and Wu~{\cite[Theorem~1.3]{Oxley2000b}}]\label{thm:OW2}
  The class of $3$-connected matroids that have exactly two non-essential elements, each of which is deletable, coincides with the class of matroids that are constructed as described in (i)--(vi) of~{\cite[Theorem~1.3]{Oxley2000b}}.
\end{theorem}

\begin{theorem}[Oxley and Wu~{\cite[Theorem~1.4]{Oxley2000b}}]\label{thm:OW3}
  The class of $3$-connected matroids that have exactly two non-essential elements, one of which is deletable and one of which is contractible, coincides with the non-wheels and non-whirls that are in the class of matroids constructed as described in (i)--(vi) of~{\cite[Theorem~1.4]{Oxley2000b}}.
\end{theorem}

Because deletable elements of a matroid are precisely contractable elements of its dual matroid, the above two theorems give a full characterization of $3$-connected matroids that have exactly two non-essential elements.
The construction of matroids described in Theorems~1.3 and~1.4 in~\cite{Oxley2000b} is complicated, but it is simpler for binary matroids as described below.
The following corollaries can be deduced directly from Corollary~5.4 (or Theorems~5.1 and~5.2) in~\cite{Oxley2000b}.
\begin{corollary}\label{cor:binmat2}
  The class of $3$-connected binary matroids that have exactly two non-essential elements, each of which is deletable, coincides with the class of the cycle matroids of multidimensional wheels.
\end{corollary}

\begin{corollary}\label{cor:binmat3}
  The class of $3$-connected binary matroids that have exactly two non-essential elements, one of which is deletable and one of which is contractible, coincides the class of
  \begin{enumerate}[label=\rm(\roman*)]
    \item the cycle matroids of twisted wheels and
    \item $M_1(n_1,\ldots,n_k)$ and $M_2(n_1,\ldots,n_k)$ for $k\ge 3$ and $n_1,\ldots,n_k\ge 1$.
  \end{enumerate}
\end{corollary}

If a matroid $M$ is isomorphic to $U_{1,3}$, $U_{2,3}$, or the cycle matroid of a triangle-sum of $n$ wheels with $n\ge 2$, then $M$ has exactly three non-essential elements; see the first two lines of page 225 and the first sentence of the proof of Theorem~1.2 in~\cite{Oxley2000b}.
Thus, we deduce the following proposition from~\cite[Corollary 5.4]{Oxley2000b}.

\begin{proposition}\label{prop:binmat}
  A $3$-connected binary matroid with at least four elements has exactly two non-essential elements if and only if it is isomorphic to one of the following matroids:
  \begin{enumerate}[label=\rm(\roman*)]
    \item\label{item:twisted wheel} the cycle matroid of a twisted wheel,
    \item\label{item:multidimensional wheel} the cycle or cocycle matroid of a multidimensional wheel, and
    \item\label{item:M1M2} $M_1(n_1,\ldots,n_k)$ or $M_2(n_1,\ldots,n_k)$ for some $k\ge 3$ and $n_1,\ldots,n_k \ge 1$.
  \end{enumerate}
\end{proposition}

Then Corollaries~\ref{cor:binmat2} and~\ref{cor:binmat3} follow Proposition~\ref{prop:binmat} by checking whether non-essential elements in matroids from \ref{item:twisted wheel}--\ref{item:M1M2} in Proposition~\ref{prop:binmat} are deletable or contractible.
This can be easily done by observing their fundamental graphs.

As an example, let us check non-essential elements in the cycle matroid $M$ of a twisted wheel.
  Let $B$ be its base inducing a fundamental graph $G$ isomorphic to $\theta(1,a,b)$ with odd $a,b\ge 5$.
  Let $x$ and $y$ be vertices of degree~$3$ in $G$ such that $x\in B$ and $y\in V(G)-B$, respectively.
  Then as $G/x$ and $G/y$ are prime, by Lemmas~\ref{lem:fundamental graphs of binary matroids} and~\ref{lem:fundamental graphs of binary matroids 2}, $x$ is deletable and $y$ is contractible in $M$.
  Moreover, since neither $G\setminus x$ nor $G\setminus y$ is prime, $x$ is not contractible and $y$ is not deletable in $M$.
  This strategy of deciding deletable and contractible elements similarly works for other matroids in the previous proposition.

  Now, the following lemma immediately shows that 
Proposition~\ref{prop:binmat} 
and Corollary~\ref{cor:bippiv3} 
are equivalent,
assuming Corollary~\ref{cor:bippiv}.
We omit its easy proof.

\begin{lemma}
  The following are equivalent for a binary matroid $M$.
  \begin{enumerate}[label=\rm(\roman*)]
    \item  $M$ is isomorphic to one of \ref{item:twisted wheel}--\ref{item:M1M2} in Proposition~\ref{prop:binmat}.
    \item $M$ has a fundamental graph isomorphic to one of the following.
    \begin{itemize}
      \item $\theta(1,a,b)$ with odd integers $a,b\ge 5$,
      \item $\theta(2\ell_1,\ldots,2\ell_k)$ with $k\ge 3$ and $\ell_1,\ldots,\ell_k \ge 2$, and
      \item $\theta(2n_{1}+1,\ldots,2n_{k}+1)$ or $\theta(1,2n_{1}+1,\ldots,2n_{k}+1)$ with $k\ge 3$ and $n_1,\ldots,n_k \ge 1$.
    \end{itemize}
  \item $M$ has a fundamental graph that is isomorphic to a graph in $\Theta$
  and is not pivot-equivalent to an even cycle.\qed
  \end{enumerate}
\end{lemma}

\end{document}